\documentclass[onefignum,onetabnum]{siamart190516}

\usepackage{braket,amsfonts}

\usepackage{array}

\usepackage[caption=false]{subfig}

\usepackage{pgfplots}

\usepackage{cite,alltt}

\newsiamthm{claim}{Claim}
\newsiamremark{remark}{Remark}
\newsiamremark{hypothesis}{Hypothesis}
\crefname{hypothesis}{Hypothesis}{Hypotheses}
\usepackage{multirow}
\usepackage{booktabs}
\usepackage{algorithm,algorithmicx,algpseudocode}

\usepackage{graphicx,epstopdf}

\Crefname{ALC@unique}{Line}{Lines}

\usepackage{amsopn}

\usepackage{xspace}
\usepackage{bold-extra}
\usepackage{bm}
\usepackage[most]{tcolorbox}

\colorlet{texcscolor}{blue!50!black}
\colorlet{texemcolor}{red!70!black}
\colorlet{texpreamble}{red!70!black}
\colorlet{codebackground}{black!25!white!25}

\lstdefinestyle{siamlatex}{%
  style=tcblatex,
  texcsstyle=*\color{texcscolor},
  texcsstyle=[2]\color{texemcolor},
  keywordstyle=[2]\color{texemcolor},
  moretexcs={cref,Cref,maketitle,mathcal,text,headers,email,url},
}

\tcbset{%
  colframe=black!75!white!75,
  coltitle=white,
  colback=codebackground, % bottom/left side
  colbacklower=white, % top/right side
  fonttitle=\bfseries,
  arc=0pt,outer arc=0pt,
  top=1pt,bottom=1pt,left=1mm,right=1mm,middle=1mm,boxsep=1mm,
  leftrule=0.3mm,rightrule=0.3mm,toprule=0.3mm,bottomrule=0.3mm,
  listing options={style=siamlatex}
}

\newtcblisting[use counter=example]{example}[2][]{%
  title={Example~\thetcbcounter: #2},#1}

\newtcbinputlisting[use counter=example]{\examplefile}[3][]{%
  title={Example~\thetcbcounter: #2},listing file={#3},#1}

\DeclareTotalTCBox{\code}{ v O{} }
{ %fontupper=\ttfamily\color{texemcolor},
  fontupper=\ttfamily\color{black},
  nobeforeafter,
  tcbox raise base,
  colback=codebackground,colframe=white,
  top=0pt,bottom=0pt,left=0mm,right=0mm,
  leftrule=0pt,rightrule=0pt,toprule=0mm,bottomrule=0mm,
  boxsep=0.5mm,
  #2}{#1}

\patchcmd\newpage{\vfil}{}{}{}
\newcommand{\R}{\mathbb{R}}
\newcommand{\M}{\mathcal{M}}

\newcommand{\ds}{\displaystyle}

\newcommand{\vect}[1]{\mathbf{#1}}

\newcommand{\vx}{\vect{x}}
\newcommand{\xh}{\hat{x}}
\newcommand{\yh}{\hat{y}}
\newcommand{\vxh}{\vect{\xh}}
\newcommand{\vy}{\vect{y}}

\newcommand{\vQ}{R}

\newcommand{\vn}{\vect{n}}

\newcommand{\vt}{\vect{t}}

\newcommand{\vxi}{\boldsymbol{\xi}}

\newcommand{\X}{\mathbf{X}}
\newcommand{\Xh}{\hat{\X}}
\newcommand{\uc}{\underline{c}}
\newcommand{\ud}{\underline{d}}
\newcommand{\ub}{\underline{b}}
\newcommand{\ua}{\underline{a}}
\newcommand{\uu}{\underline{u}}
\newcommand{\us}{\underline{s}}
\newcommand{\up}{\underline{p}}

\newcommand{\ulam}{\underline{\lambda}}

\newcommand{\lap}{\Delta}
\newcommand{\laps}{\lap_{\M}}
\newcommand{\laph}{\hat{\lap}}

\newcommand{\diag}{\text{diag}}
\DeclareMathOperator*{\argmin}{argmin}

\DeclareMathOperator{\rank}{\text{rank}}

\usepackage{alltt}
\usepackage{color}
\definecolor{string}{rgb}{0.7,0.0,0.0}
\definecolor{comment}{rgb}{0.13,0.54,0.13}
\definecolor{keyword}{rgb}{0.0,0.0,1.0}

\begin{tcbverbatimwrite}{tmp_\jobname_header.tex}
\title{MGM: A meshfree geometric multilevel method for systems arising from elliptic equations on point cloud surfaces}

\author{Grady B. Wright\thanks{Boise State University (\email{gradywright@boisestate.edu}).}
\and Andrew M. Jones\thanks{Boise State University (\email{andrewjones237@u.boisestate.edu}).}
\and Varun Shankar\thanks{University of Utah (\email{shankar@cs.utah.edu}).}}

\headers{MGM: A method for systems arising from elliptic PDEs on surfaces}{G.~B. Wright, A.~M. Jones, V.~Shankar}
\end{tcbverbatimwrite}
\input{tmp_\jobname_header.tex}

\ifpdf
\hypersetup{ pdftitle={Meshfree Geometric Multilevel Methods} }
\fi

\begin{document}
\maketitle

%% ------------------------------------------------------------------
%% ABSTRACT
%% ------------------------------------------------------------------
%\begin{tcbverbatimwrite}{tmp_\jobname_abstract.tex}
\begin{abstract}
We develop a new meshfree geometric multilevel (MGM) method for solving linear systems that arise from discretizing elliptic PDEs on surfaces represented by point clouds.  The method uses a Poisson disk sampling-type technique for coarsening the point clouds and new meshfree restriction/interpolation operators based on polyharmonic splines for transferring information between the coarsened point clouds.  These are then combined with standard smoothing and operator coarsening methods in a V-cycle iteration.  MGM is applicable to discretizations of elliptic PDEs based on various localized meshfree methods, including RBF finite differences (RBF-FD) and generalized finite differences (GFD).  We test MGM both as a standalone solver and preconditioner for Krylov subspace methods on several test problems using RBF-FD and GFD, and numerically analyze convergence rates, efficiency, and scaling with increasing point cloud sizes.  We also perform a side-by-side comparison to algebraic multigrid (AMG) methods for solving the same systems.  Finally, we further demonstrate the effectiveness of MGM  by applying it to three challenging applications on complicated surfaces: pattern formation, surface harmonics, and geodesic distance.

\end{abstract}


%\end{tcbverbatimwrite}
%\input{tmp_\jobname_abstract.tex}
%% ------------------------------------------------------------------
%% END HEADER
%% ------------------------------------------------------------------

%%%%%%%%%%%%%%%%%%%%%%%%%%%%%%%%%%%%%%%%%%%%%%%%%%%%%%%%%%%%%%%%%%%%%%%%%
%INTRODUCTION
%%%%%%%%%%%%%%%%%%%%%%%%%%%%%%%%%%%%%%%%%%%%%%%%%%%%%%%%%%%%%%%%%%%%%%%%%
\section{Introduction}\label{sec:intro}
%Partial differential equations (PDEs) defined on surfaces (or manifolds) arise in many areas of science and engineering, with the most prominent being the atmospheric sciences, where spherical surfaces play a key role.  However, PDEs on more general surfaces also arise, for example, in models for chemical diffusion on cell membranes, morphogenesis, nonlinear chemical oscillators in excitable media, and textures in computer graphics.  Solutions of these models can rarely be achieved by analytical means and must instead be approximated using numerical techniques.  

%Partial differential equations (PDEs) defined on surfaces (or manifolds) arise in many areas of science and engineering, where, for example, they are used to model atmospheric flows,  chemical diffusion on cell membranes, morphogenesis, nonlinear chemical oscillators in excitable media, and textures for computer graphics.  Solutions of these models can rarely be achieved by analytical means and must instead be approximated using numerical techniques

Partial differential equations (PDEs) defined on surfaces (or manifolds) arise in many areas of science and engineering, where they are used to model, for example, atmospheric flows~\cite{Williamson:2007DynamicalCores}, chemical signaling on cell membranes~\cite{Liue2104191118}, morphogenesis~\cite{stoop2015curvature}, and textures for computer graphics~\cite{turk1991generating}. Solutions of these models can rarely be achieved by analytical means and must instead be approximated using numerical techniques. While numerical methods for PDEs on the sphere have been developed since the 1960s~\cite{Williamson:2007DynamicalCores}, development of methods for PDEs on more general surfaces only began in the late 1980s~\cite{Dziuk88}, with interest growing considerably in the early 2000s~\cite{dziuk_elliott_2013}.  These techniques include surface finite element (SFE)~\cite{dziuk_elliott_2013}, embedded finite element (EFE)~\cite{BertalmioEtAl2001,olshanskii2017trace}, and closest point (CP)~\cite{MacDondaldRuuth2009} methods.  More recently, various meshfree (or meshless) methods have also been developed for PDEs on general surfaces that use a local stencil approach, including radial basis function-finite differences (RBF-FD)~\cite{Alvarez2021,LSW2016,SHANKAR2014JSC,PIRET2016,SHANKAR2018722,Wendland2020}, generalized finite differences (GFD)~\cite{SUCHDE20192789}, and generalized moving least squares (GMLS)~\cite{LiangZhao13,TraskKuberry20,GrossEtAl20}.  These methods can be applied for surfaces represented only by point clouds and do not require a surface triangulation like SFE methods or a level-set representation of the surface like EFE  methods.  Additionally, these meshfree methods approximate the solutions directly on the point cloud and do not extend the PDEs into the embedding space like the EFE and CP methods.

In this paper, we concentrate on local meshfree methods for elliptic PDEs on surfaces, which are challenging to solve with iterative methods because of the poor conditioning of the systems.  We specifically focus on the surface Poisson and shifted (or screened) surface Poisson problems, which, for example, in surface hydrodynamics~\cite{GrossEtAl20}, computer graphics~\cite{reuter2009discrete}, and time-implicit discretizations of surface reaction diffusion equations~\cite{SHANKAR2014JSC}.  We focus on two methods for these PDEs: polyharmonic spline-based RBF-FD with polynomials and GFD.  These meshfree discretizations result in large, sparse, non-symmetric, linear systems of equations that need to be solved.  Direct solvers for these systems have most commonly been used, but these do not scale well to large point clouds and high-orders of accuracy, motivating the need for efficient and robust iterative methods.
%with a few studies using Krylov iterative methods preconditioned with algebraic multigrid (AMG)~\cite{GrossEtAl20} or incomplete LU~\cite{LSW2016}.  To make such methods practical for large scale problems, effective iterative methods need to be developed for solving these systems.

Multigrid methods are known to be effective solvers and preconditioners for linear systems that arise from discretizing elliptic PDEs (e.g.,~\cite{mgtrott}).  These methods can be classified into two types: geometric and algebraic.  While algebraic multigrid (AMG) methods are general purpose solvers/preconditioners, geometric methods, when they can be developed, generally converge faster and work as better preconditioners for Krylov subspace methods.  Geometric multigrid methods have been developed for SFE and CP discretizations (e.g.,~\cite{landsberg2010multigrid} and~\cite{chen2015closest}) and it is the aim of this paper to also develop these methods for meshfree discretizations.

The basic components of geometric multigrid methods that need to be developed are (1) techniques for coarsening the grid or mesh, (2) constructing interpolation/restriction operators for transferring the information between levels, (3) discretizing the differential operator on the coarser levels, (4) smoothing the approximate solution, and (5) solving the system of the coarsest level.  The first component presents a challenge for meshfree surface PDEs as there is only a point cloud available and no grid or mesh to create a hierarchy of coarser levels. To overcome this challenge we use the weighted sample elimination (WSE) method from~\cite{Yuksel2015}, which is a general purpose method for selecting quasi-uniformly spaced subsets of points from a point cloud and falls into the general category of Poisson disk sampling methods~\cite{Bridson07}.  The lack of a grid or mesh also presents a challenge for component (2) as standard transfer operators cannot be used.  To overcome this challenge, we use RBFs to construct the interpolation operators for transferring the defect from coarser to finer levels.  For the restriction operators, we simply use the transpose of the interpolation operators, which is a standard choice~\cite{mgtrott}.  With these transfer operators, we generate component (3) using a Galerkin projection, often referred to as the Galerkin coarse grid operator.  Finally, for component (4) we use standard Gauss-Seidel smoothing and for (5) we use a direct solver.  We combine all of these components in a V-cycle iteration and apply it both as a solver and preconditioner.  The resulting method is entirely meshfree and we refer to it as the meshfree geometric multilevel\footnote{We use the term multilevel rather than multigrid, since this method does not depend on a grid.} (MGM) method.

The new MGM method has some similarities to the meshfree \textit{multicloud} methods~\cite{katz2009multicloud,zamolo2019novel}, but also some key differences.  The first major difference is that  multicloud methods have been developed for PDEs posed in planar domains, whereas MGM is for surface PDEs.  Another difference is with the choice of transfer operators. The method of~\cite{zamolo2019novel} uses one-point, piecewise constant operators, while the method of~\cite{katz2009multicloud} uses two-point, inverse-distance weighted interpolation and restriction operators.  It is not clear how these latter transfer operators should be generalized to surfaces. MGM instead uses transfer operators based on RBFs, which are well suited for interpolation on surfaces~\cite{FuselierWright:2010}.  A second difference is the strategy for geometric coarsening of the given point cloud.  Multicloud methods use a graph coloring-type scheme for finding maximally-independent subsets of vertices to determine the coarser levels.  This does not allow the size of the point clouds on the coarser levels to be controlled precisely, and it limits the coarsening factors to approximately four (for 2D problems).  MGM instead uses WSE~\cite{Yuksel2015}, which allows for arbitrary coarsening factors and for the sizes of the points in the coarser levels to be controlled exactly.  A third difference is that MGM can handle degenerate PDEs (e.g., the surface Poisson equation), while the multicloud methods have been tailored to non-degenerate PDEs (e.g. planar Poisson equation with mixed Dirichlet-Neumann boundary conditions).  Finally, multicloud methods have only been tested on second order accurate discretizations of PDEs; these discretizations typically only use small stencils.  We demonstrate that MGM works for discretizations at least up to sixth order accurate with large stencil sizes.   

%There are geometric multilevel methods that have previously been developed for meshfree discretizations of PDEs, which are termed multicloud methods~\cite{katz2009multicloud,zamolo2019novel}.  MGM has some similarities to these methods, but also some key differences.  However, there are some key differences....
%\begin{itemize}
%\item \cite{zamolo2019novel} uses piecewise constant restriction and interpolation operators and   They develop their method for MQ RBFs and small stencil sizes (n=7 in 2D and n=14 in 3D), which means they are of low accuracy.  Their coarsening strategy only uses one coarse level node for the interpolation, this does keep their stencil sizes roughly constant through the levels.   They use aggregation based interpolation and restriction operators which introduce more tuning parameters. Does not handle degenerate operators like surface Poisson problem.  Uses SOR for the smoother.
%\item ~\cite{katz2009multicloud} is based on RBF-FD using multiquadric with only linear polynomials appended -- low order.  They apply their method to the Euler equations which are discretized on a fine mesh using a conservative FV scheme.  They only use multicloud on the coarse levels to accelerate convergence.  Coarsen using a graph coloring-type scheme for finding maximally independent set of vertices.  Does not allow one to have precise control on the number of points on the coarsest level.
%\end{itemize}
  
The remainder of the paper is organized as follows.  In Section \ref{sec:meshfree_review}, an overview is given of the two meshfree methods for surface PDEs that the MGM algorithm is used to solve.  The next section presents the RBF-based transfer operators.  Section \ref{sec:mgm} describes the remaining components of the MGM algorithm, including a discussion of the changes necessary to solve the degenerate surface Poisson problem.  Section \ref{sec:results} presents an extensive array of numerical results for the MGM method, including a comparison with algebraic multigrid (AMG) methods.  Section \ref{sec:applications} then uses the method in three challenging applications to further demonstrate its effectiveness. Finally, the paper concludes with some final remarks on the method and some future directions in Section \ref{sec:discussion}.

\subsection{Assumptions and notation}\label{sec:assumptions}
Throughout the manuscript we let $\M$ be a smooth embedded manifold of co-dimension one in $\R^3$ with no boundary and let $\laps$ denote the Laplace-Beltrami operator (LBO) (or surface Laplacian) on $\M$.  When referencing points on $\M$, we assume they are represented as coordinates in $\R^3$, e.g., for $\vx\in\M$, and we write $\vx = (x,y,z)$. For a point $\vx\in\M$, we let $T_{\vx}\M$ denote the tangent plane to $\M$ at $\vx$. We denote normal vectors to $\M$ as $\vn$ and assume that they are available either analytically or using some approximation technique (see for example~\cite{Klas_09}). 

We use sub/superscripts $h$ and $H$ on variables to indicate whether they are associated with the fine or coarse level point clouds, respectively.  For example, $X_h$ and $X_H$ denote the set of points in the fine level and coarse level point clouds, respectively.  This is meant to mimic the notation that is used in traditional grid based geometric multigrid methods~\cite{mgtrott}, but these parameters do not relate to anything specific about the spacing of the points and do not need to be computed.

The focus of this study is on the elliptic equation
\begin{align}
\mathcal{L} u = f,
\label{eq:poisson}
\end{align}
where $u:\M\rightarrow\mathbb{R}$ is unknown and $f:\M\rightarrow\mathbb{R}$ is known.  Here $\mathcal{L} = \laps$ or $\mathcal{L} = \mathcal{I} - \mu \laps$, where $\mathcal{I}$ is the identity operator and $\mu > 0$, which correspond to the surface Poisson and shifted (or screened) surface Poisson equation \eqref{eq:poisson}, respectively.  
%The shifted surface Poisson equation appears, for example, in the time-implicit discretization of surface diffusion problems.  
For the surface Poisson problem, we assume $f$ satisfies the \textit{compatibility condition} $\int_{\M} f\,dA=0$, which is a necessary and sufficient condition for \eqref{eq:poisson} to have a solution.  In this case, any solution is unique up to the addition of a constant since constants satisfy the homogeneous equation.

%%%%%%%%%%%%%%%%%%%%%%%%%%%%%%%%%%%%%%%%%%%%%%%%%%%%%%%%%%%%%%%%%%%%%%%%%
%METHODOLOGY RBF METHODS
%%%%%%%%%%%%%%%%%%%%%%%%%%%%%%%%%%%%%%%%%%%%%%%%%%%%%%%%%%%%%%%%%%%%%%%%%

\section{Localized meshfree discretizations\label{sec:meshfree_review}}
Several localized meshfree methods have been developed for approximating the solution of~\eqref{eq:poisson}, e.g.,~\cite{LiangZhao13,TraskKuberry20,Alvarez2021,LSW2016,SHANKAR2014JSC,PIRET2016,SHANKAR2018722}.  For the sake of brevity, we limit the focus of this study to two localized meshfree methods: polyharmonic spline (PHS)-based RBF-FD with polynomials and GFD. Both of these methods use the so-called tangent plane approach, but differ in the approximation spaces used.  They should be sufficient to demonstrate the general applicability of the MGM method. 

The RBF-FD and GFD methods are based on approximating the strong form of the equation, and amount to discretizing the LBO $\laps$ over a local stencil of points on the surface. This stencil based approach can generally be described as follows.  Let $\X_h=\{\vx_i\}_{i=1}^{N_h}$ denote the global point cloud (node set) discretizing $\M$.  For each $\vx_i\in \X_h$, $i=1,\ldots,N_h$, let $\sigma^h_{i}$ denote the set of indices of the $n_i>1$ nearest neighbor nodes in $\X_h$ to $\vx_j$.  Here we use the Euclidean distance in $\R^3$ to define the nearest neighbor distances.  The points $\X_h^{i} = \{\vx_j\}_{j\in \sigma_i^h}$ are the \textit{stencil} for $\vx_i$, and $\vx_i$ is the \textit{stencil center}. The stencil based approximation to $\laps u$ at $\vx_j$ then takes the form
\begin{align}
%\laps u\bigr|_{\vx_j} \approx \sum_{i=1}^{n_j} c_i^{(j)} u(\vx_{i}^{(j)}),
\laps u\bigr|_{\vx_i} \approx \sum_{j\in\sigma_h^i} c_{ij} u(\vx_{j}),
\label{eq:general_fd_approx}
\end{align}
%where $c_i^{(j)}$ 
where $c_{ij}$ are some set of weights determined by the RBF-FD or GFD methods discussed below.  These weights can then be assembled into a global $N_h$-by-$N_h$ (sparse) differentiation matrix $D_h$ and an approximate solution to \eqref{eq:poisson} is given as a solution to the linear system 
\begin{align}
L_h u^h = f^h,
\label{eq:discrete_system}
\end{align}
where $u^h,f^h\in\mathbb{R}^{N_h}$ contain the unknown solution and known right hand side of \eqref{eq:poisson}, respectively, sampled at $\X_h$.  The matrix $L_h$ is given by $L_h = D_h$ or $L_h = I_h - \mu D_h$, where $I_h$ is the $N_h$-by-$N_h$ identity matrix.  Note that the latter $L_h$ arises from time-implicit discretizations of surface diffusion-type problems.  The MGM method will be used for solving the system \eqref{eq:discrete_system}.

In the remainder of this section, we give specific details on determining the stencil weights in \eqref{eq:general_fd_approx} for the LBO using RBF-FD and GFD methods. Since both of these methods use the tangent plane technique, we review it first.

\subsection{Tangent plane method\label{sec:tangent_plane}}

The tangent plane idea was introduced by Demanet~\cite{DEMANET2006} and recently further refined by Suchde \& Kuhnert~\cite{SUCHDE20192789}, Shaw~\cite{ShawThesis}, and, in the case of the unit two-sphere, by Gunderman et. al.~\cite{GUNDERMAN2020109256}.  The central idea of the method is to approximate the LBO at the center of each stencil $\X_h^{i}$ using an approximation to the standard Laplacian on the plane tangent to the surface at the stencil center, $T_{\vx_i}\M$.  The approximation is constructed from a projection of the stencil points to $T_{\vx_i}\M$.  In this work, we use the projection advocated in~\cite{SUCHDE20192789}, which is known as an orthographic projection when $\M$ is the unit sphere.
\begin{figure}[htb]
\centering
\begin{tabular}{cc}
\includegraphics[width=0.4\textwidth]{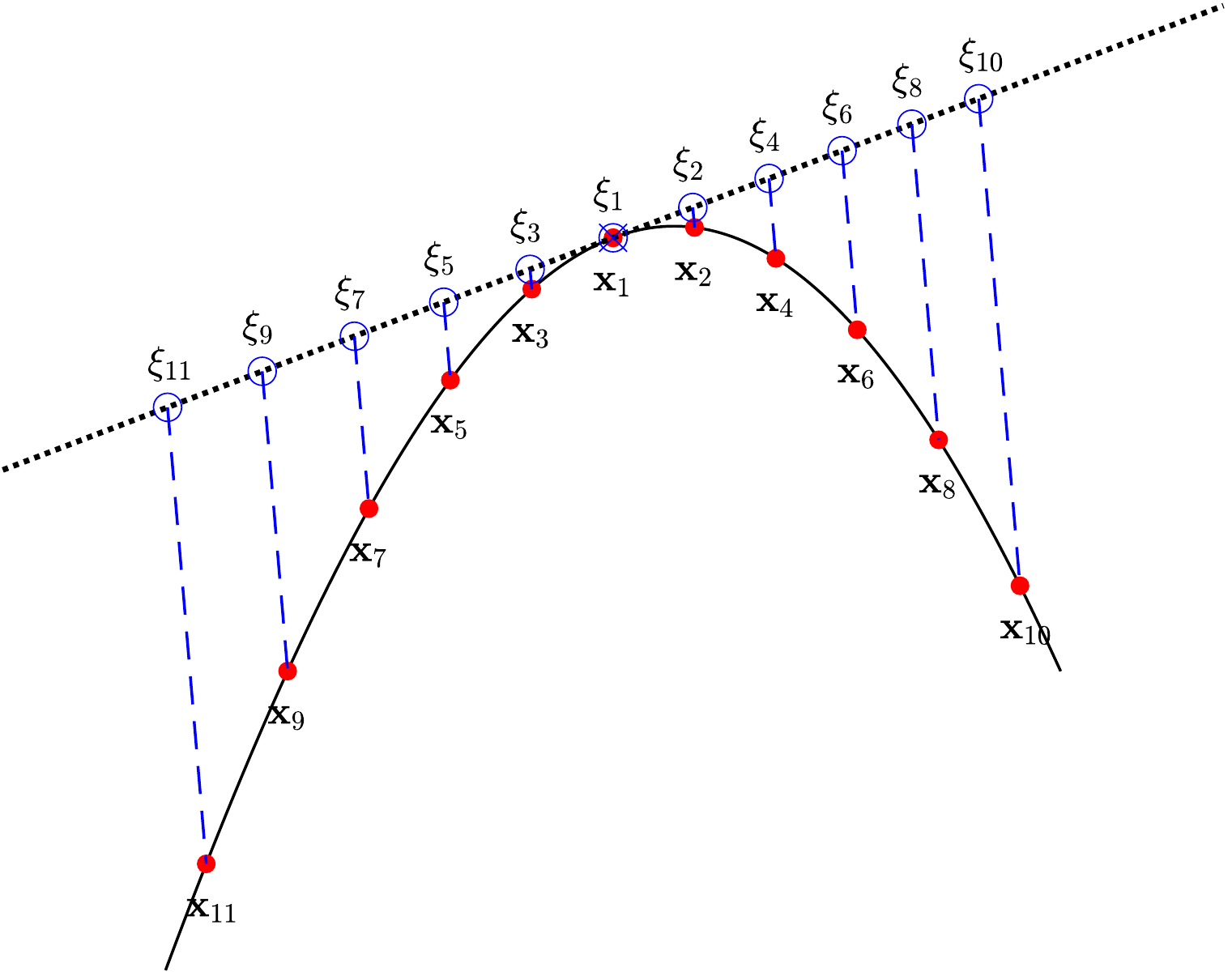} & \includegraphics[width=0.4\textwidth]{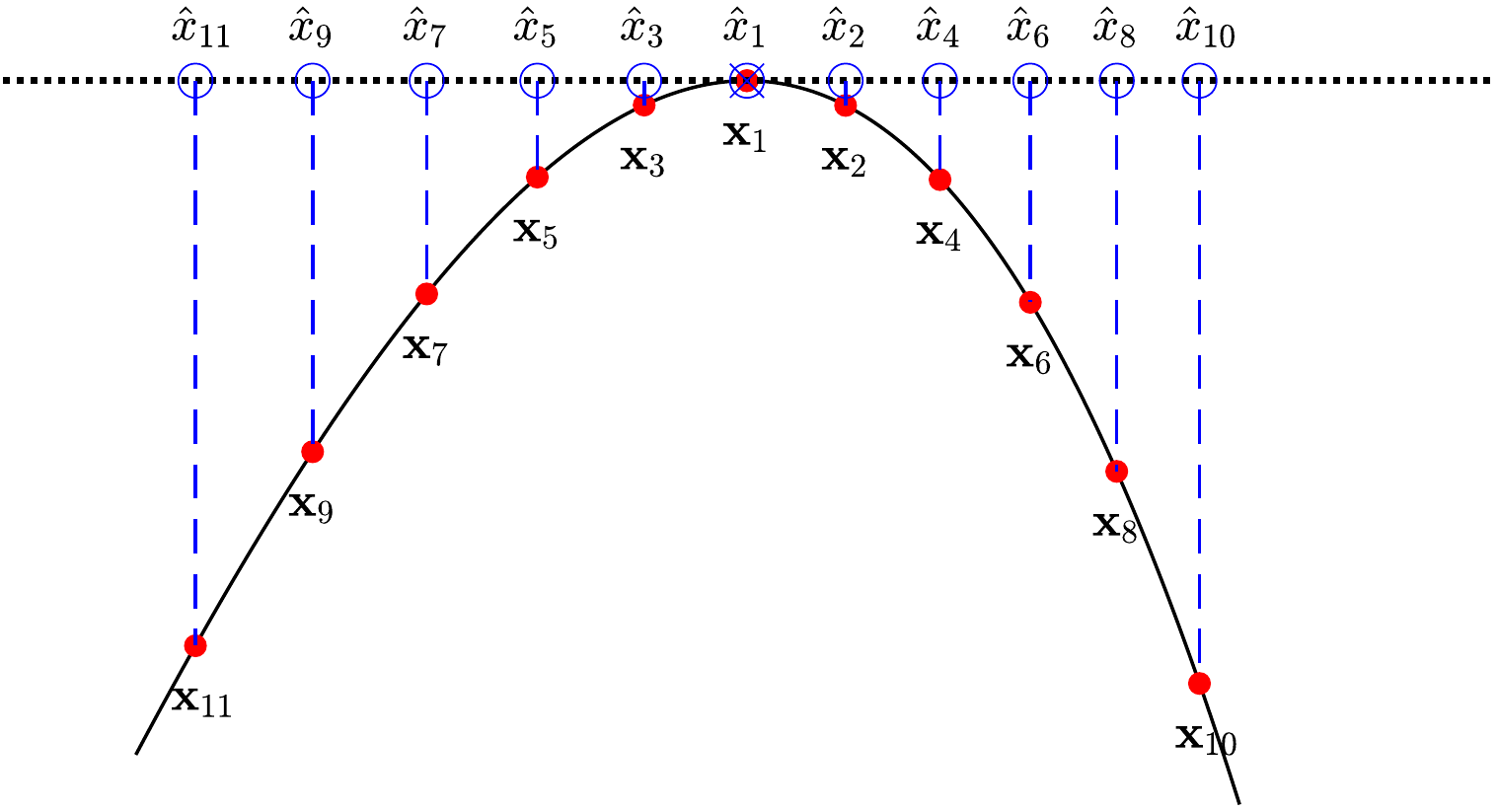} \\
(a) & (b)
\end{tabular}
\caption{Illustration of the tangent plane method for a 1D surface (curve). The solid black lines indicates the surface, the solid red circles mark the $n=11$ stencil nodes, the open blue circles mark the projected nodes, and the $\times$'s marks the stencil center.  (a) Direct projection of the stencil points according to \eqref{eq:projection1}. (b) Rotation and projection of the stencil points according to \eqref{eq:tp_pts}. \label{fig:projection}}
\end{figure}

%With out loss of generality, consider the first stencil with $n_1=n$ points, and to simplify notation, drop the superscripts on the stencil points so that $\X^{(1)}_h=\{\vx_1,\vx_2,\ldots,\vx_n\}$.  The tangent plane method from~\cite{SUCHDE20192789} projects these points onto $T_{\vx_1}\M$ along the normal vector $\vn_1$ (the normal to the surface at $\vx_1$).  This projection is illustrated in Figure \ref{fig:projection} (a) for a one dimensional surface (curve) in $\R^2$.  The projected points can be computed explicitly by
With out loss of generality, we describe the method for the first stencil $\X^1_h$ with index set $\sigma_h^1$.  To simplify notation, we assume $\sigma_h^1 = \{1,\ldots,n\}$, so that the stencil points are simply $\X_h^{1}=\{\vx_1,\ldots,\vx_n\}$.  The tangent plane method from~\cite{SUCHDE20192789} projects these points onto $T_{\vx_1}\M$ along the normal vector $\vn_1$ (the normal to the surface at $\vx_1$).  This projection is illustrated in Figure \ref{fig:projection} (a) for a one dimensional surface (curve) in $\R^2$.  The projected points can be computed explicitly by
\begin{align}
%\vxi_j = \vx_j + \vn_1^{T}(\vx_1-\vx_j)\vn_1 = \underbrace{(I-\vn_1\vn_1^T)\}_{\ds P}vx_j + \vn_1^{T}\vx_1\vn_1,\; j=1,\ldots,n.
\vxi_j = (I-\vn_1\vn_1^T)(\vx_j-\vx_1),\; j=1,\ldots,n,
\label{eq:projection1}
\end{align}
where we have shifted the projected points so that $\vxi_1$ is at the origin.  For a two dimensional surface, the projected points can be expressed in terms of orthonormal vectors $\vt_1$ and $\vt_2$ that span $T_{\vx_1}\M$ as follows
\begin{align}
\vxi_j  = \underbrace{\begin{bmatrix} \vt_1 & \vt_2 \end{bmatrix}}_{\ds \vQ} \underbrace{\begin{bmatrix} \xh_j \\ \yh_j \end{bmatrix}}_{\ds \vxh_j},\; j=1,\ldots,n.
\label{eq:projection2}
\end{align}
The 2D coordinates $\vxh_j$ in the tangent plane for the projected points are what will be used for constructing the approximations to the Laplacian.  These can be computed directly from the relationship \eqref{eq:projection1} and \eqref{eq:projection2} as 
\begin{align}
\vxh_j = \vQ^T (\vx_j - \vx_1),\; j=1,\ldots,n,
\label{eq:tp_pts}
\end{align}
where we have used $\vQ^T (I - \vn_1\vn_1^T) = \vQ^T$.  We denote the the projected stencil as $\Xh_h^{1} = \{\vxh_1,\ldots,\vxh_n\}$.  The above procedure is repeated for every stencil $\X_h^{i}$,  to obtain the projected stencils $\Xh_h^{i}$, $i=1,\ldots,N_h$.

We note that, geometrically speaking,  \eqref{eq:tp_pts} amounts to first shifting the stencil points so the center is at the origin, rotating them so the normal $\vn_1$ is orthogonal to the $xy$-plane, and then dropping the third component.  This is illustrated in Figure \ref{fig:projection} (b) for the case of a 1D curve.

\subsection{PHS-based RBF-FD with polynomials}\label{sec:rbffd_lap}
The RBF-FD method for determining the weights in \eqref{eq:general_fd_approx} can be derived by constructing an RBF interpolant over each of the projected stencil points to the tanget plane, applying the standard 2D Laplacian to the interpolants, and then evaluating them at the stencil center.  In this study we focus on interpolants constructed from PHS kernels and polynomials~\cite{bayona2017role,Du77,ShawThesis}.  Without loss of generality, we again describe the method for the first stencil, with $\X_h^1 = \{\vx_1,\ldots,\vx_n\}$, to simplify notation.

For the stencil $\X_h^{1}$, the PHS interpolant to the projected stencil $\Xh_h^{1}$ takes the form
\begin{align}
s(\vxh) = \sum_{i=1}^{n} a_i \|\vxh - \vxh_i\|^{2k+1} + \sum_{j=1}^{L} b_j p_j(\vxh),
\label{eq:phs_interp}
\end{align}
where $\vxh=\begin{bmatrix} \xh & \yh \end{bmatrix}^T \in T_{\vx_1}\M$, $\|\cdot\|$ denotes the Euclidean norm, $k$ is the order of the PHS kernel, and $\{p_1,\ldots,p_L\}$ is a basis for bivariate polynomials in $T_{\vx_1}\M$ of degree $\ell$ (so that $L = (\ell+1)(\ell+2)/2$).  The order $k$ controls the smoothness of the PHS and is chosen such that $0\leq k\leq \ell$. We note that the polynomials can be chosen to be the standard bivariate monomials in the components of $\vxh$.  For samples $\{u_1,\ldots,u_n\}$ of an arbitrary function at the stencil points $\X_h^1$, the coefficients for the interpolant in the tangent plane are determined by the conditions
\begin{align}
s(\vxh_i) = u_i,\; i=1,\ldots,n \quad \text{and} \quad \sum_{i=1}^n a_i p_j(\vxh_i) = 0,\; j=1,\ldots,L,
\label{eq:phs_interp_conditions}
\end{align}
which can be written as the following linear system
\begin{align}
\begin{bmatrix}
A & P \\
P^T & \vect{0}
\end{bmatrix}
\begin{bmatrix}
\ua \\
\ub
\end{bmatrix}
=
\begin{bmatrix}
\uu \\
\underline{0}
\end{bmatrix},
\label{eq:rbf_fd_interp_system}
\end{align}
where $A_{ij} = \|\vxh_i - \vxh_j\|^{2k+1}$ ($i,j=1,\ldots,n$), $P_{ij} = p_j(\vxh_i)$ ($i=1,\ldots,n$, $j=1,\ldots,L$), and underlined terms denote vectors containing the corresponding terms in \eqref{eq:phs_interp_conditions}.  If the stencil nodes $\Xh_h^{1}$ are unisolvent with respect to the space of bivariate polynomials of degree $\ell$ (i.e., $\text{rank}(P)=L$), then \eqref{eq:rbf_fd_interp_system} has a unique solution, so that \eqref{eq:phs_interp} is well-posed~\cite{Wendland:2004}.  This is a mild condition on the stencil nodes, especially for ``scattered'' nodes on the tangent plane.

The stencil weights $c_{1j}$ in \eqref{eq:general_fd_approx} are determined from the approximation
\begin{align*}
\laps u\bigr|_{\vx_1} \approx \laph s \bigr|_{\vxh_1} &= \sum_{i=1}^{n} a_i \laph(\|\vxh - \vxh_i\|^{2k+1})\bigr|_{\vxh_1} + \sum_{j=1}^{L} b_j \laph(p_j(\vxh))\bigr|_{\vxh_1} \\
& = \begin{bmatrix} \laph\us & \laph\up \end{bmatrix}^T\begin{bmatrix} \ua \\ \ub \end{bmatrix},
\end{align*}
where $\laph = \partial_{\xh\xh} + \partial_{\yh\yh}$,  $\laph\us$  and $\laph\up$ are vectors containing the entries $\laph(\|\vxh - \vxh_i\|^{2k+1})\bigr|_{\vxh_1}$, $i=1,\ldots,n$, and  $\laph(p_j(\vxh))\bigr|_{\vxh_1}$, $j=1,\ldots,L$, respectively. Using \eqref{eq:rbf_fd_interp_system} in the above expression, the stencil weights are given as the solution to the following linear system
\begin{align}
\begin{bmatrix}
A & P \\
P^T & \vect{0}
\end{bmatrix}
\begin{bmatrix}
\uc \\
\ulam
\end{bmatrix}
=
\begin{bmatrix}
\laph\us \\
\laph\underline{p}
\end{bmatrix},
\label{eq:rbf_fd_lap}
\end{align}
where $\uc$ contains $c_{1j}$ and $\ulam$ are unused.  

We note that one can interpret \eqref{eq:rbf_fd_lap} as the solution to an equality constrained optimization problem where the weights are determined by enforcing they are exact for $\laph(\|\vxh - \vxh_i\|^{2k+1})\bigr|_{\vxh_1}$, $j=1,\ldots,n$, subject to the constraint that they are also exact for $\laph(p_j(\vxh))\bigr|_{\vxh_1}$, $j=1,\ldots,L$.  Under this interpretation, $\ulam$ is the vector of Lagrange multipliers~\cite{bayona2017role}.

In this work, we choose the order of the PHS as $k=\ell$ and fix the stencil size $n_j=n$, $j=1,\ldots,N$ as $n = \lceil2 L\rceil$, which is a common choice for RBF-FD methods~\cite{bayona2017role}.  The degree $\ell$ of the appended polynomial can then be used to control the approximation order of the method, with larger $\ell$ leading to higher orders~\cite{davydov2019optimal}.

\subsection{GFD}

This method is similar to the RBF-FD method, but instead of using an interpolant, the method is based on a (weighted) polynomial least squares approximant.  Using the same notation and assumptions as the previous section and again focusing only on the first stencil $\X_h^1$, the approximant for the projected stencil $\Xh_h^{1}$ takes the form
\begin{align}
q(\vxh) = \sum_{j=1}^L b_j p_j(\vxh).
\label{eq:gmls_interp}
\end{align}
%where $\{p_1,\ldots,p_L\}$ are again a basis for bivariate polynomials of degree $\ell$ in $\T$.  
The coefficients of the approximant are determined from the samples $\{u_1,\ldots,u_n\}$ according to the the following weighted least squares problem:
\begin{align}
\ub^{*} = \argmin_{\ub\in\mathbb{R}^n} \sum_{i=1}^n w(\vxh_i) (q(\vxh_i) - u_i)^2 =   \argmin_{\ub\in\mathbb{R}^n} \|W^{1/2}(P\ub - \uu)\|_2^2,
\label{eq:wls_fd_approx_system}
\end{align}
where $W = \diag(w(\vxh_i))$.  Here we again assume the stencil nodes $\Xh_h^{1}$ are unisolvent with respect to bivariate polynomials of degree $L$ so that \eqref{eq:wls_fd_approx_system} has a unique solution.  There are many different options for selecting the weight function $w$ in the literature.  In this work, we follow~\cite{SUCHDE20192789} and use the following Gaussian function:
\begin{align*}
w(\vxh_i) = \exp\left(-\alpha\frac{\|\vxh_1 - \vxh_i\|^2}{\rho_1^2 + \rho_i^2}\right),
\end{align*}
where $\rho_k$ is the support of the $k$th projected stencil $\Xh_h^i$, i.e.\ the radius of the minimum ball centered at $\vxh_i$ that encloses all the points in $\Xh_h^i$.  The parameter $\alpha > 0$ is used for controlling the shape of the weight function and is typically chosen in an ad hoc manner~\cite{SUCHDE20192789}.

The stencil weights in \eqref{eq:general_fd_approx} are determined from the approximation
\begin{align*}
\laps u\bigr|_{\vx_1} \approx \laph q \bigr|_{\vxh_1} &= (\laph \up)^T \ub^*.
\end{align*}
Using the normal equation solution of \eqref{eq:wls_fd_approx_system} in the above expression, the vector of stencil weights $\uc$ is given as
\begin{align}
\uc = WP(P^TWP)^{-1} (\laph \up).
\label{eq:wls_fd_lap}
\end{align}
In practice, a QR factorization of $WP$ is used instead of the normal equations to improve the numerical conditioning of \eqref{eq:wls_fd_lap}.

In this work, we choose the number of points in the stencils $\X_h^{i}$ for this method in the same manner as the RBF-FD technique.  Note that, similar to RBF-FD, increasing the polynomials degree $\ell$ also increases the order of accuracy of the GFD method.

\section{Multilevel transfer operators using RBFs \label{sec:transer_ops}}
Operators for transferring information between coarse and fine levels are one of the key components of multilevel methods.  The interpolation transfer operators are used to transfer information from coarse level to a finer level, while the restriction operators are used for the reverse.  Let $\X_h=\{\vx_i\}_{i=1}^{N_h}\subset\M$ denote the fine set of nodes and $\X_H=\{\vy_j\}_{j=1}^{N_H}\subset\M$ the coarse set, where $N_H < N_h$.  We denote the interpolation operator by $I_{H}^h$ and the restriction operator by $I_h^H$.  These can be represented of as (sparse) matrices, so that for a vector $u_H$ of data on the coarse nodes $\X_H$, the vector containing the interpolation of $u^H$ to $\X_h$ is given as $u^h=I_H^{h}u_H$.  In this section, we discuss a novel meshfree method for constructing $I_{H}^h$ based on RBF interpolation.  For the restriction operator, we use $I_h^H = (I_H^h)^T$, which is a standard choice, especially in AMG methods~\cite[Appendix A]{mgtrott}.

Similar to the discrete LBO, we compute the interpolation operator using a stencil based approach.  Letting $\sigma_H^i$ be the indicies of the $m_i$ nearest neighbors in $\X_H$ to $\vx_i$, the interpolation of $\{u^H_j\}_{j\in\sigma^H_i}$ to $u^{h}_i$, the entry in $u_h$ corresponding to $\vx_i$, is given as
\begin{align}
u^h_i = \sum_{j\in\sigma^i_{H}} d_{ij} u^H_{j}.
\label{eq:interp_fd_approx}
\end{align}
We again use the Euclidean distance in $\mathbb{R}^3$ to define the index set $\sigma_H^i$.  The weights for each stencil can be assembled to form the (sparse) interpolation matrix $I_{H}^{h}$.

We use local RBF interpolants about each stencil $\X_H^i=\{\vy_j\}_{j\in\sigma^H_i}$ to determine the weights in \eqref{eq:interp_fd_approx} and form these interpolants in the embedding space $\mathbb{R}^3$.  This is considerably simpler than using intrinsic coordinates to $\M$ and the resulting interpolants have good approximation properties~\cite{FuselierWright:2010}.  One could alternatively use interpolants in the tangent plane about each stencil center similar to Section \ref{sec:rbffd_lap}, but these are only accurate when the surface is well discretized by the underlying point cloud.  This will not necessarily be the case for the nodes $\X_{H}$ as we coarsen the the finer levels.  We again use a PHS kernel to form the interpolants and describe the method for the first stencil $\X_H^1$, which, to simplify the notation, we assume to consist of the nodes $\{\vy_1,\ldots,\vy_m\}$.

For the stencil $\X_H^1$, the PHS interpolant takes the same form as \eqref{eq:phs_interp}, but with $\vxh$ replaced by $\vx$ and $\{\vxh_i\}$ is replaced by $\{\vy_i\}$.  Additionally, we only consider the PHS kernel with $k=0$ and a constant term appended to the interpolant (i.e. $\ell = 0$).  The weights $d_{1j}$ in \eqref{eq:interp_fd_approx} are determined by evaluating this interpolant at $\vx_1$ and can be computed in a similar  procedure to that used in deriving the system \eqref{eq:rbf_fd_lap}.  The linear system for the interpolation weights takes the form
\begin{align}
\begin{bmatrix}
A & \underline{1} \\
\underline{1}^T & \vect{0}
\end{bmatrix}
\begin{bmatrix}
\ud \\
\lambda
\end{bmatrix}
=
\begin{bmatrix}
\us \\
1
\end{bmatrix},
\label{eq:rbf_fd_interp}
\end{align}  
where $A_{ij}= \|\vy_i-\vy_j\|$ ($i,j=1,\ldots,m$), $\underline{1}$ is the vector of length $n$ with all ones, and $\underline{s}$ has entries $\|\vx_1 - \vy_i\|$, $i=1,\ldots,m$.  Again, $\lambda$ is unused.

While higher order PHS kernels ($k > 0$) and higher degree polynomials ($\ell > 0$) could be used in constructing the interpolation weights, we found that the simple formulation above gave good results, while also being efficient, for the range of problems we considered. This formulation also has the added benefit that the system \eqref{eq:rbf_fd_interp} has a unique solution, provided the points are distinct~\cite{Wendland:2004}. When using larger $k$ and $\ell$ this may not be the case as the points must be unisolvent with respect to the space of trivariate polynomials of degree $\ell$, i.e., $\rank(P)=L$.  Since the interpolation is done in the embedding space, this can be an issue for certain algebraic surfaces (e.g., the sphere with $\ell \geq 2$).

\section{Meshfree geometric multilevel (MGM) method}\label{sec:mgm}
In this section we present the MGM method for solving the discrete problem \eqref{eq:discrete_system}.  We first present the MGM method in terms of a two-level cycle, which is summarized in Algorithm \ref{alg:two_level}, describing its primary components: coarsening the point cloud $\X_h\rightarrow \X_H$, forming the coarse level operator $L_H$, smoothing the approximation, and solving for the defect on the coarse level.  The interpolation/restriction operators are described in the previous section.  We then focus on some modifications to the algorithm that are necessary when \eqref{eq:discrete_system} corresponds to the surface Poisson problem.  This is followed by a description of the multilevel extension of the method.  Finally, we comment on using the method as a preconditioner for Krylov subspace methods.
%The primary components of the method that we still have left to discuss are techniques for coarsening the node sets, forming the coarse level operator, smoothing the defect, and solving the coarse system.

%=================================
\begin{algorithm}[t]
\begin{algorithmic}[1]
\State  Pre-smooth initial guess: \phantom{} $u^h \leftarrow \text{presmooth}(L_h,u^h,f^h,\nu_1)$
\State  Compute residual: \phantom{} $r^h = f^h - L_h u^h$
\State  Restrict the residual to $\X_H$: \phantom{}$r^H = I_h^H r^h$
\State  Solve for the defect: \phantom{}$L_H e^H = r^H$
\State  Interpolate the defect to $\X_h$: \phantom{}$e^h = I_H^h e^H$
\State  Correct the approximation: \phantom{}$u^h \leftarrow u^h + e^h$
\State  Post-smooth the approximation: \phantom{}$u^h \leftarrow \text{postsmooth}(L_h,u^h,f^h,\nu_2)$
%\State   Repeat until $\|r^h\|_2 \leq \ tol\|f^h\|$
\caption{Two-Level Cycle\label{alg:two_level}}
\end{algorithmic}
\end{algorithm}

\subsection{Node coarsening}

The technique we propose for generating the coarser point clouds on general surfaces is based on (WSE) method from~\cite{Yuksel2015}.   This algorithm falls into the category of Poisson disk sampling methods, which produce quasiuniformly spaced point sets~\cite{Bridson07}.  The WSE method approximates the solution to the following optimization problem:  Given a point cloud $\X_h$ with $N_h$ samples, determine a subset $\X_H$ of $\X_h$ with $N_H$ samples that has maximal Poisson disk radius.  The Poisson disk radius is defined as one half the minimum distance between neighboring points in the set (which is called the separation radius in the meshfree methods literature~\cite{Wendland:2004}).  This optimization problem is $NP$ complete, but the WSE algorithm approximates the solution in $N_h-N_H$ steps with a theoretical complexity of $\mathcal{O}(N_h \log N_h)$ operations~\cite{Yuksel2015}.  The method works for point clouds defined on many different sampling domains, including arbitrary manifolds, where it uses the Euclidean norm in $\mathbb{R}^3$ to define nearest neighbor distances.  We use the implementation by the author of the WSE method, called \texttt{cySampleElimination}, that is provided in the \texttt{cyCodeBase}~\cite{YukselSoftware}.

\begin{figure}[ht]
\centering
\begin{tabular}{cc}
\includegraphics[width=0.22\textwidth]{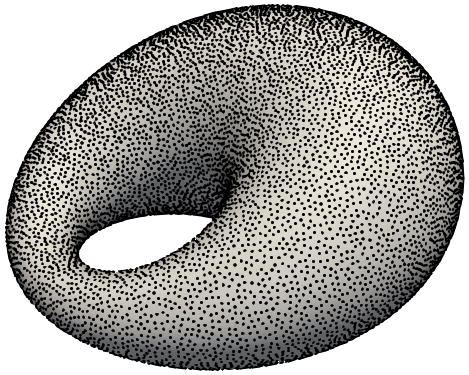}  \includegraphics[width=0.22\textwidth]{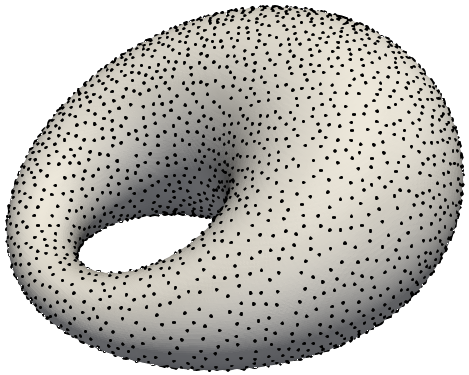} & \includegraphics[width=0.22\textwidth]{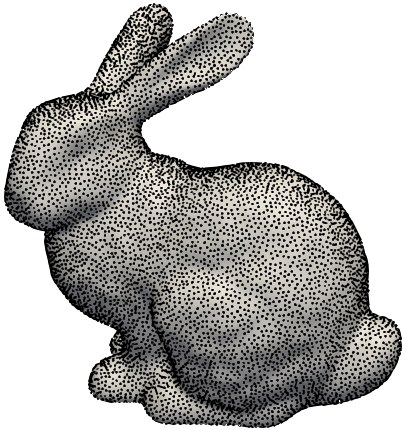}  \includegraphics[width=0.22\textwidth]{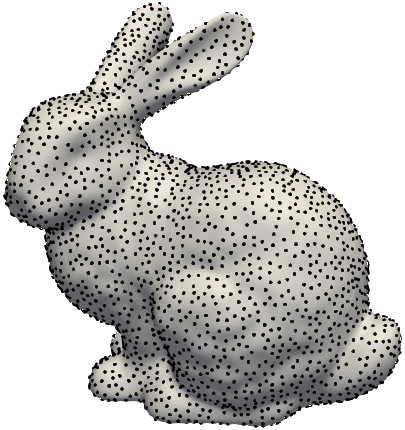} \\
(a) $\X_h \rightarrow \X_H$ & (b) $\X_h \rightarrow \X_H$
\end{tabular}
\caption{Illustration of the WSE algorithm for generating a coarse level set $\X_H$ of $N_H=\lfloor N_h/4 \rfloor$ points from a fine level set $\X_h$ of $N_h$ points.  Here $N_h=14561$ \& $N_H=3640$ for the cyclide (a) and $N_h=14634$ \& $N_H= 3658$ for the Stanford Bunny (b) .\label{fig:WSECoarsening}}
\end{figure}

In this work, we coarsen the point cloud $\X_h$ by a fixed factor of 4, so that $\X_H$ has $N_H=\lfloor N_h/4 \rfloor$ points.  This mimics the standard coarsening of geometric multigrid for two-dimensional domains.  We tested other coarsening factors, but found that coarsening by 4 generally gave the best results in terms of iteration count and wall clock time for the multilevel method. Figure \ref{fig:WSECoarsening} illustrates the coarse point clouds $\X_H$ with this coarsening factor computed from the WSE algorithm for two example surfaces.  

%\subsubsection{Point sets for the unit sphere}
%In the important case of the unit sphere $\Sph^2$, one can use structured, quasi-uniform points for the coarse node sets $\X_H$ as an alternative to WSE.  Several options for these point sets are available such as icosahedral, cubed sphere, HEALPix, and Fibonacci~\cite{HMS16}; see~\cite{SpherePts} for implementations in MATLAB.  An important distinction between using these node sets instead of those from WSE, is that they do not necessarily form nested subsets of the original fine level point cloud $\X_h$.  The meshfree interpolation schemes we propose in Section \ref{sec:transer_ops} do not require such a nesting and, as we demonstrate in the results section, excellent results can be obtained for these classes of coarse level point sets, even when the fine level point cloud is not of the same class.  The benefits of using this procedure is that it eliminates setup costs of the method associated with running the WSE algorithm.  For some of these coarser point sets, it may not be possible to coarsen by an exact factor of 4.  In these cases, we choose the value of $N_H$ for these node sets to be as close to $N_h/4$ as possible.
%
%This type of coarsening can also be used for more general surface where a hierarchy of points sets may be available, either from analytical means or from some node generator (e.g.,~\cite{ShankarKirbyFogelson}).

\subsection{Coarse level operator}
There are two main approaches to constructing the coarse level operator in multilevel methods.  The first is \textit{direct discretization}, where the differential operator is discretized directly on the coarse level points $\X_H$.  The second is based on a Galerkin projection involving the interpolation $I_H^h$ and restriction $I_h^H$ operators, and is defined as follows:
\begin{align}
L_H = I_h^H L_h I_H^h.
\label{eq:galerkin}
\end{align}
This latter operator, referred to as the \textit{Galerkin coarse grid operator}, provides a simple means of coarsening $L_h$ and has been shown to be robust for a large class of problems, especially those where a direct discretization on the coarse grid does not adequately represent the approximation on the fine grid~\cite[\S 7.7.4]{mgtrott}.  It also gives rise to a variational principle that is exploited in the analysis of algebraic multigrid (AMG)~\cite[\S A.2.4]{mgtrott} methods.  While this latter result relies on the matrix being symmetric positive define, modifications to this theory have also been developed for non-symmetric problems, which involve choosing the restriction operator differently than the transpose of the interpolation operator~\cite{ManteuffelSouthworth19}.  
We use the Galerkin approach for forming $L_H$, as we have found that it approximates the fine grid operator on the coarser levels better than the direct discretization technique and it makes for a more robust solver/preconditioner.  While $L_h$ is not symmetric for our discretizations, we have nonetheless found that simply choosing $I_h^H = (I_H^h)^T$ works well over a large array of test problems.  

One disadvantage of the Galerkin approach is that $L_H$ has to be formed explicitly through sparse matrix-matrix multiplication, which is more computationally expensive in terms of time and memory than the direct discretization approach. Several researchers have developed methods to reduce this cost on parallel architectures (e.g.,\cite{BellEtAll2012,BallardEtAl2016}), but we have not used these methods in our implementation.  We simply construct $L_H$ as part of a set-up phase using a sparse matrix library.  However, we do some minor alterations to improve the computational performance.  These include reordering the rows and columns of $L_h$ to decrease its bandwidth using the reverse Cuthill-McKee (RCM) algorithm prior to forming $L_H$.  This essentially leads to a reordering of the nodes $\X_h$, which in turn leads to a reordering of the interpolation operator $I_H^h$.  We also use RCM to reorder the rows and columns of $L_H$ after it is formed, which leads to a re-ordering of the nodes $\X_h$ and of the columns of $I_H^h$.  We have found that these matrix reorderings not only reduce the wall-clock time of MGM, but also the number of iterations to reach convergence.

\subsection{Smoother and coarse level solver}
For the smoothing operator we use classical Gauss-Seidel (GS) method.  One application of the smoother can be written as
\begin{align}
u^h \leftarrow u^h + B_h^{-1}(f^h - L_h u^h),
\label{eq:smoother}
\end{align}
where $B_h$ is the lower triangular part (called forward GS) or upper triangular part (called backward GS) of $L_h$.  In some cases we vary the version of the smoother for the pre- and post- smoothing operations (e.g., forward GS for pre-smooth and backward GS for the post-smooth).  We denote the number of applications of the smoother for the pre- and post-phases of the cycle as $\nu_1$ and $\nu_2$, respectively.  

To solve for the defect on the coarse level, we use a direct solver based on a sparse $LU$ factorization of $L_H$ (e.g., SuiteSparse or SuperLU).

\subsection{Modifications to the two-level cycle for the surface Poisson problem\label{sec:poisson_mods}}
When \eqref{eq:discrete_system} corresponds to the discretization of a surface Poisson problem ($L_h = D_h$), the system is singular and some modifications to the two level cycle in Algorithm \ref{alg:two_level} are necessary.  To understand the nature of the singularity, we can look at the continuous problem \eqref{eq:poisson}.  As discussed in Section \ref{sec:assumptions}, this problem has a solution if and only if the right hand side satisfies the compatibility condition.  Furthermore, the solution is only unique up to the addition of a constant.  The degeneracy in the continuous problem manifests in the discrete problem as a one dimensional null space of $L_h$ corresponding to constant vectors.  The discrete analog of the consistency condition is that \eqref{eq:discrete_system} has a solution if and only if $f^h$ is orthogonal to the left null vector of $L_h$ (i.e., $f^h$ is in the range of $L_h$). Also, similar to the continuous case, any solution of \eqref{eq:discrete_system} is only unique up to the addition of a constant vector. 

%%The analog for the discrete system is that \eqref{eq:discrete_system} has a solution if and only if $f^h$  is orthogonal to the left null space of $L_h$ (i.e., $f^h$ is in the range of $L_h$), and that a solution is only unique up to the addition of a constant vector (since constant vectors satisfy the discrete homogeneous equation). 

The primary issue that arises with using multilevel methods (and other iterative methods) for these types of singular systems stems from the fact that, in practice, $f_h$ is rarely in the range of $L_h$.  This can cause the iterations to fail to converge to a suitable approximation.  Three standard approaches to bypass this issue include the following.  First, one can project $f^h$ into the range of $L_h$.  However, this requires computing the left null vector, which can be computationally expensive\footnote{Note that $L_h$ is not symmetric, so the constant vector is not necessarily the left null vector}.  It also requires modifying the coarse level solver to use the pseudoinverse (or some approximate inverse).  A second approach is to impose that the solution is zero at one point.  This fixes the non-uniqueness issue and transforms the problem into solving a non-singular system of one dimension smaller.  However, this can lead to a deterioration of the convergence of the multilevel method since the pointwise condition is not well approximated on coarser levels~\cite{wesseling2004}.  Additionally, the solution to this approach can be less accurate and less smooth than the projection approach~\cite{YoonEtAl2016}.  The third approach is to enforce a global constraint on the solution, such as the discrete mean of $u_h$ is zero~\cite[\S 5.6.4]{mgtrott}.  This constraint can be enforced using a Lagrange multiplier, which transforms the linear system into the constrained system
\begin{align}
\underbrace{
\begin{bmatrix} 
L_h & b_h^T \\
b_h & 0
\end{bmatrix}}_{\ds \tilde{L}_h}
\underbrace{
\begin{bmatrix}
u^h \\
\lambda^h
\end{bmatrix}}_{\ds \tilde{u}^h}
=
\underbrace{
\begin{bmatrix}
f^h \\
0
\end{bmatrix}}_{\ds \tilde{f}^h},
\label{eq:saddle_point}
\end{align}
where $b_h$ is a row vector of length $N$ with all of its components set to 1 (i.e., the summation operator), and $\lambda^h$ is the Lagrange multiplier.  Provided $b_h$ is not orthogonal to the left null space of $L_h$ (which is likely to be true because of the compatibility condition for the continuous problem), this constrained system will have a unique solution~\cite[Lemma 5.6.1]{mgtrott}.  Furthermore, if this condition holds, the solution will be the same (up to a constant) as the projection approach, since  $f^h-\lambda^h b_h^T$ is then necessarily in the range of $L_h$.   We use the third approach in the MGM method.

Some modifications to the two-level cycle are required to handle the constrained system \eqref{eq:saddle_point}.  First, the transfer operators have to also transfer the Lagrange multipler through the fine and coarse levels and the Galerkin coarse grid operator has to include the transferred constraint.  We follow the approach from~\cite{Adams2004} and modify these operators  according to the following definitions:
\begin{align}
\tilde{L}_H =
\underbrace{
\begin{bmatrix}
I_h^H & 0 \\
0 & 1
\end{bmatrix}}_{\ds \tilde{I}_h^H}
\underbrace{
\begin{bmatrix} 
L_h & b_h^T \\
b_h & 0
\end{bmatrix}}_{\ds \tilde{L}_h}
\underbrace{
\begin{bmatrix}
I_H^h & 0 \\
0 & 1
\end{bmatrix}}_{\ds \tilde{I}_H^h}
=
\begin{bmatrix}
I_h^H L_h I_H^h & I_h^H b_h^T \\
b_h I_H^h & 0
\end{bmatrix},
\label{eq:gco_saddle_point}
\end{align}
where $\tilde{I}_H^h$ and $\tilde{I}_h^H$ are the modified interpolation and restriction operators, respectively, and $\tilde{L}_H$ is the modified Galerkin operator.  These modified transfer operators simply pass the Lagrange multiplier between levels without alteration.

For the smoother of the constrained system \eqref{eq:saddle_point}, we use the approach discussed in~\cite[\S 5.6.5]{mgtrott}, where only the solution $u_h$ is smoothed and the constraint is left alone.  We again use GS for smoothing $u_h$ and one application of the modified smoother takes the form
\begin{align*}
u^h \leftarrow u^h + B_h^{-1}(f^h - L_h u^h - b_h \lambda^h),
\end{align*}
where $B_h$ is the same as \eqref{eq:smoother}.  This smoother is equivalent to one iteration of the undamped inexact Uzawa method with the Schur complement  set equal to zero~\cite{benzi2005}.

Finally, we use a direct solve to compute the defect $e^h$ and Lagrange multiplier $\lambda^H$ on the coarse level.  This system takes the form $\tilde{L}_H \tilde{e}^H = \tilde{r}^H$,
%\begin{align*}
%\tilde{L}_H \tilde{e}^H = \tilde{r}^H,
%\end{align*}
where $\tilde{e}^H = \begin{bmatrix} e^H & \lambda^H\end{bmatrix}^T$ and $\tilde{r}^H$ is the restricted residual for the modified system: $\tilde{r}^H = \tilde{I}_h^H (\tilde{f}^h - \tilde{L}_h \tilde{u}^h)$. When solving a Poisson problem, we use the modifications described above in Algorithm \ref{alg:two_level}.

\subsection{Multilevel extension}
\begin{algorithm}[t!]
\begin{algorithmic}[1]
%\State Initialize node set $X^h$ of size $N$
%\State Initialize coarsening factor $CF$
\State \textbf{Input:} Fine level nodes $\X_{1}$ and operator $L_{1}$; minimum number of coarse level points $N_{\rm min}$
\State Re-order rows and columns of $L_{1}$ using RCM
\State Re-order $\X_{1}$ according to the RCM ordering
\State Compute number of levels: $p=\lfloor \log(N_{1}/N_{\rm min})/\log(4)\rfloor+1$ \label{alg:p}
\For{$j=1\ldots p-1$}
\State Generate coarse point cloud $\X_{{j+1}}$ with $N_{{j+1}} = \lfloor N_1/4^j\rfloor$ points \label{alg:wse}
\State Generate interpolation operator $I_{{j+1}}^{{j}}$ from $\X_{{j+1}}$ to $\X_{{j}}$
\State Set the restriction operator to $I_{j}^{{j+1}} = (I_{{j+1}}^{{j}})^{T}$
\State Generate Galerkin coarse level operator $L_{{j+1}} = I_{{j+1}}^{{j}} L_{{j}} I_{j}^{{j+1}}$
\State Re-order rows and columns of $L_{{j+1}}$ using RCM
\State Re-order rows of $I_{{j+1}}^{{j}}$ and columns of $I_{j}^{{j+1}}$ according to the RCM ordering
\EndFor
\State Compute sparse LU decomposition of $L_{p}$
\caption{MGM preprocessing phase\label{alg:preprocessing}}
\end{algorithmic}
\end{algorithm}
The multilevel extension of the two-level cycle can be obtained by applying it recursively until a sufficiently coarse level is reached to make a direct solver practical.  To simplify the notation in describing the multilevel cycle, we replace the $h/H$ superscript/subscript notation with a number corresponding to the level, with $j=1$ being the finest level.  For example, for the $j$th level, $\X_j$ denotes the point cloud, $N_j$ denotes its size, $L_j$ denotes the operator, $r^j$ denotes the residual, and $I_{j}^{j+1}$ is the restriction to level $j+1$.  
%=================================

Before the multilevel cycle begins, we compute all the coarse point clouds, transfer operators, and Galerkin coarse level operators in a preprocessing step, which is outlined in Algorithm \ref{alg:preprocessing}.  The number of levels, $p$, depends on the number of fine level nodes and minimum number of nodes on the coarsest level, $N_{\rm min}$, and is determined on line \ref{alg:p} of this algorithm.  This guarantees that the number of nodes on the coarsest level satisfies $N_{\rm min} \leq N_p < 4N_{\rm min}$.  We note that when using WSE to generate the coarse point cloud $\X_{j}$ on line \ref{alg:wse} of the preprocessing algorithm, we use the finer point cloud $\X_{j+1}$, rather than the finest node $\X_{1}$.  This reduces the cost in performing this step.

The multilevel cycle is outlined in Algorithm \ref{alg:mgm} in non-recursive form.  This algorithm is what we call the MGM method and corresponds to a traditional V-cycle in multigrid methods, which is typically denoted $\rm V(\nu_1,\nu_2)$ corresponding to the number of pre-/post smoothing operations.  Other cycling methods can also be used (e.g.,\ F- or W-cycle~\cite[\S 2.4]{mgtrott}), but we limit our focus to the V-cycle.  While this algorithm is described for a shifted Poisson problem, it can be easily modified for solving a Poisson problem following the modifications discussed in Section \ref{sec:poisson_mods}.
 
%As discussed in Section \ref{sec:intro} we discuss the M3 algorithm and its primary components:
%node set selection,  node coarsening with weighted sample elimination (WSE) from \cite{Yuksel2015}, transfer and operator generation and storage, smoothing, restriction, coarse level correction (CLC), interpolation, correct. We first detail the preprocessing stage of M3, then we show the two and multilevel V-cycle algorithms \cite{mgtrott}.

%=================================
\begin{algorithm}[tbh]
\begin{algorithmic}[1]
%\State Initialize node set $X^h$ of size $N$
%\State Initialize coarsening factor $CF$
\State \textbf{Input:} Right hand side $f^{1}$; Initial guess $u^1$; Number levels $p$; $\{L_{j}\}_{j=1}^{p-1}$; $\{I_{{j+1}}^{{j}}\}_{j=1}^{p-1}$; $\{I_{j}^{{j+1}}\}_{j=1}^{p-1}$;
\State \phantom{\textbf{Input:}}  RCM re-orderings; Sparse LU factorization of $L_p$; 
\State Re-order $f^{1}$ and $u^1$ according to RCM re-ordering of $L_{1}$
\State Presmooth initial guess: $u^1 \leftarrow \text{presmooth}(L_{1},u^1,f^1,\nu_1)$
\State Compute/restrict residual: $r^{1} =I_{1}^{{2}}(f^1 - L_{1} u^h)$
\For{$j=2\ldots p-1$}
\State Presmooth defect: $e^{j} = \text{presmooth}(L_{j},0,r^{j},\nu_1)$
\State Compute/restrict residual: $r^{{j+1}} = I_{j}^{{j+1}} (r^{j} - L_{j} e^{j})$
\EndFor
\State Compute defect: Solve $L_{p} e^{p} = r^{p}$ using sparse LU decomposition of $L_{p}$
\For{$j=p-1,\ldots,2$}
\State Interpolate/correct defect: $e^{j} \leftarrow  e^{j} + I_{{j+1}}^{{j}} e^{{j+1}}$
\State Post smooth defect: $e^{j} \leftarrow \text{postsmooth}(L_{j},e^{j},r^{j},\nu_2)$
\EndFor
\State Interpolate defect/correct approximation: $u^{1} \leftarrow  u^{1} + I_{{2}}^{{1}} e^{{2}}$
\State Post smooth approximation: $u^{1} \leftarrow \text{postsmooth}(L_{1},u^{1},f^{1},\nu_2)$
\State Undo re-ordering of $u^1$ from RCM of re-ordering of $L_1$
\caption{MGM $\rm V(\nu_1,\nu_2)$-cycle \label{alg:mgm}}
\end{algorithmic}
\end{algorithm}

%%=============================
%\begin{algorithm}[H]
%\begin{algorithmic}[1]
%\State PreSmooth $L_F\hat{u}_F = f_F$
%\State Compute $r_F = f_F - L_Fu_F$
%\State Restrict  $r_{k-1} = R_{k} r_F$
%
%\WHILE{$k \geq k_{C}$ }
%\State Smooth $L_{k-1}\hat{v}_{k-1} = r_{k-1}$}\;
%\State Compute Defect $d_{k-1} = r_{k-1} - L_{k-1}  \hat{v}_{k-1}$}\;
%\State Restrict  $r_{k-1} = R_{k-1} d_{k-1}$}\;
%\ENDWHILE
%
%\State Solve Correction $L_{C} \hat{v}_{C} = r_{C}$}
%\State $\hat{v}_{k-1} = \hat{v}_{C}$}
%\WHILE {$k \leq k_{F}$}
%  \State Interpolate and Add Smoothed Correction $v_k = I_k \hat{v}_{k-1} + \hat{v}_k$}\;
%  \State Smooth $L_{k}v_{k} = r_{k}$}\;    
%\ENDWHILE
%
%\State Correct Approximation $u_k = \hat{u}_k + v_k$
%\State PostSmooth $L_Fu_F = f_F$
%\State Continue Until Relative $||r_F||_2 \leq \ tol$ }
%\caption{Multilevel V-Cycle}
%\end{algorithmic}
%\end{algorithm}
%%=================================

\subsection{Preconditioner for Krylov subspace methods}
The MGM method has the benefit of being relatively straightforward to implement.  However, as shown in the numerical experiments in the next section, it may converge slowly when using it as a standalone solver, especially for higher order discretizations of the LBO on more irregular point clouds.  A common approach to bypassing these issues for standard geometric and algebraic multigrid methods is to combine them with a Krylov subspace method (e.g.,~\cite[\S 7.8]{mgtrott} or \cite{OoWiWaGa00,GhaiEtAl2018}).  In this case, multigrid is viewed as preconditioner for the Krylov method.  We also take this approach with MGM, using it a preconditioner for two Krylov methods: generalized minimum residual (GMRES) and bi-conjugate gradient stabilized (BiCGSTAB)~\cite{Saad}.  This combination appears to result in an efficient and robust method for solving the discretized surface Poisson and shifted Poisson equations on quite complicated surface.

%%%%%%%%%%%%%%%%%%%%%%%%%%%%%%%%%%%%%%%%%%%%%%%%%%%%%%%%%%%%%%%%%%%%%%%%%
%RESULTS
%%%%%%%%%%%%%%%%%%%%%%%%%%%%%%%%%%%%%%%%%%%%%%%%%%%%%%%%%%%%%%%%%%%%%%%%%

\section{Numerical results}\label{sec:results}
In this section, we analyze the MGM method as a solver and preconditioner for the Poisson and shifted Poisson problem on two surfaces: the unit sphere and the cyclide.  The latter is shown in Figure \ref{fig:WSECoarsening} and the implicit equation describing the surface is given in~\cite{LSW2016}.  We test the method on both RBF-FD and GFD discretizations using the parameters given in the first part of Table \ref{tbl:parameters}.  In all the tests, we are interested in how the method scales to higher order discretizations, and thus give results for polynomial degrees $\ell=3$, 5, and 7, which correspond to approximately second, fourth, and sixth order accuracy, respectively~\cite{ShawThesis,SUCHDE20192789}.  For the sphere tests, we generate the point clouds from the vertices of icosahedral node sets, which are used extensively in numerical weather prediction~\cite{majewski2002operational}. For the cyclide, we use point clouds produced from Poisson disk sampling of the surface.  This latter approach results in much more unstructured point clouds than the sphere case (see Figure \ref{fig:WSECoarsening} for an illustration).  
\begin{table}[htb]
\begin{tabular}{| c | p{0.67\textwidth} | c | }
\hline
Variable & \multicolumn{1}{c}{Description} & Value(s) \\
\hline
\hline
\multicolumn{3}{|c|}{Parameters for the discretization the LBO} \\
\hline
$\ell$ & Poly.\ degree for discretizing the LBO with RBF-FD or GFD & $3$, $5$, or $7$ \\
\hline
$k$ & Order of the PHS kernel for discretizing the LBO with RBF-FD & $\ell$ \\
\hline
$\alpha$ & Weighting parameter for the Gaussian kernel in GFD & $4$ or $5$ \\
\hline
$n$ & Stencil size for discretizing the LBO with RBF-FD or GFD & $\lceil (\ell+1)(\ell+2) \rceil$ \\
\hline
\hline
\multicolumn{3}{|c|}{Parameters for MGM} \\
\hline
$N_{\rm min}$ & Minimum number of nodes on the coarsest level & 250 \\ 
\hline
$B_h$ & Pre- and post-smoother (see \eqref{eq:smoother}) & Forward GS \\ 
\hline
$\nu_1, \nu_2$ & Number of applications of the pre- and post-smoother & $1$ \\
\hline
$m$ & Stencil size of the interpolation/restriction operators & 3\\
\hline
\end{tabular}
\caption{Description of parameters and their values used in the numerical results.\label{tbl:parameters}}
\end{table}

Unless otherwise specified, the parameters of the MGM method are set according to those given in the second part of Table \ref{tbl:parameters}.  We tested the method with different combinations of these parameters and found that the ones listed in the table generally gave the best results in terms of iteration count and wall-clock time.  Additionally, when using MGM with Krylov methods, we use it as a right preconditioner, which is generally recommended~\cite{GhaiEtAl2018}.  Finally, all the MGM results presented were obtained from a MATLAB implementation of the method, with a MEX interface to the WSE method, which is implemented in C++. 

In the first several experiments, we compare MGM to AMG, as implemented in the Python package PyAMG~\cite{OlSc2018}.  In addition to being very popular blackbox solvers and preconditioners for a wide range of problems, AMG methods have been used previously for solving linear systems associated with meshfree discretizations of elliptic PDEs in the plane~\cite{seibold2010performance} and on surfaces~\cite{GrossEtAl20}.  We use the smoothed aggregation version of AMG, as we found it performed better than classical AMG.  Additionally, we use one application of symmetric GS as the pre- and post-smoother, a V-cycle for the multilevel cycle, and sparse LU for the coarse level solver.  We experimented with other combinations of parameters and again found these generally gave the best results in terms of iteration count and wall-clock time.  Additionally, when using PyAMG with GMRES, we use it as a right preconditioner (with the fgmres option), while for BiCGSTAB we use it as a left preconditioner (as this is the only option).  Finally, in the comparisons with AMG, we focus on the shifted Poisson problem as PyAMG does not offer a specialized way to deal with the constrained system \eqref{eq:saddle_point}.

\begin{figure}[tbh]
\centering
\begin{tabular}{ccc}
 & {\small Sphere} & {\small Cyclide} \\
\rotatebox{90}{\hspace{0.1\textwidth}\small $\ell=3$} & \includegraphics[width=0.45\textwidth]{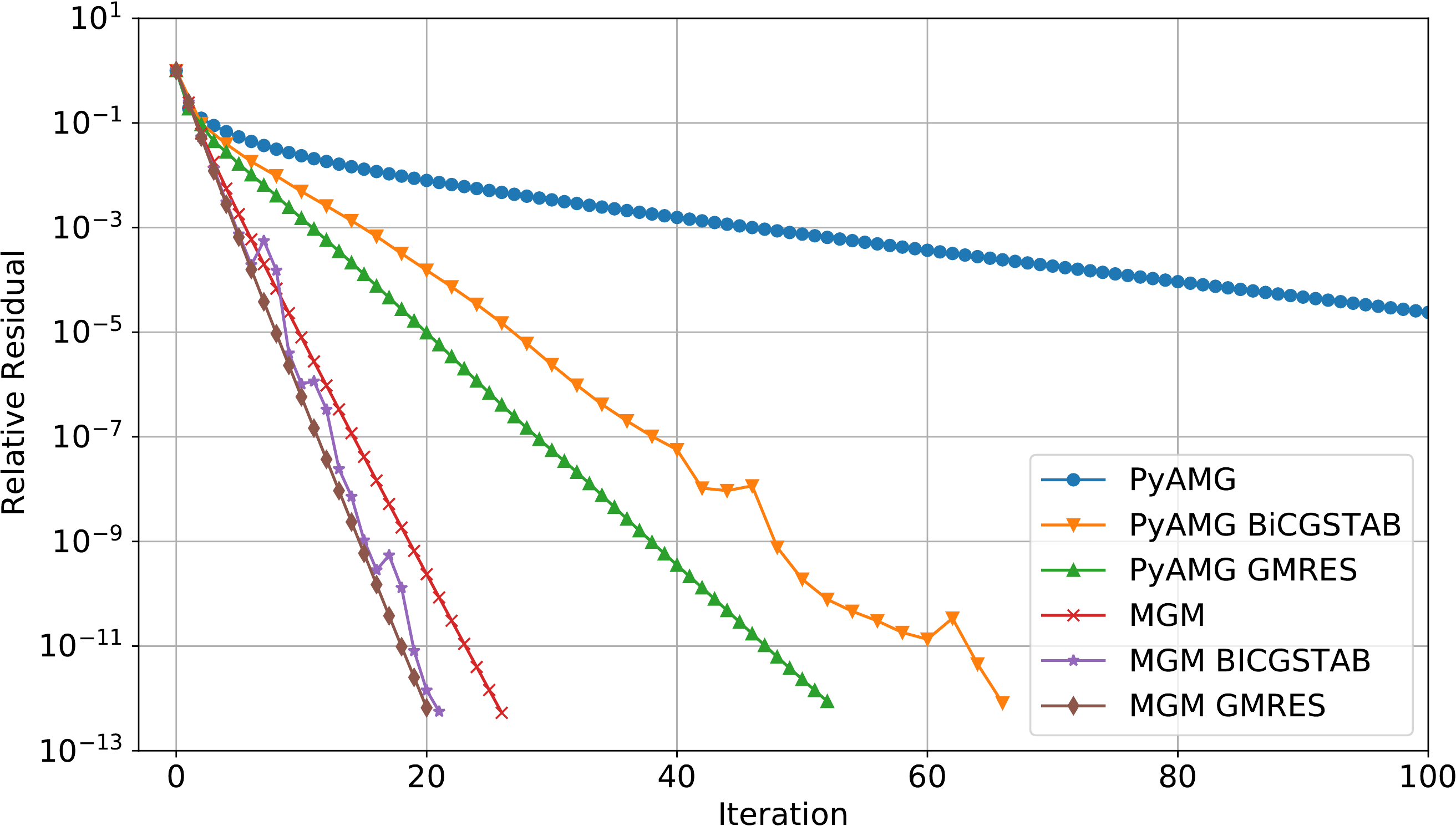} &
\includegraphics[width=0.45\textwidth]{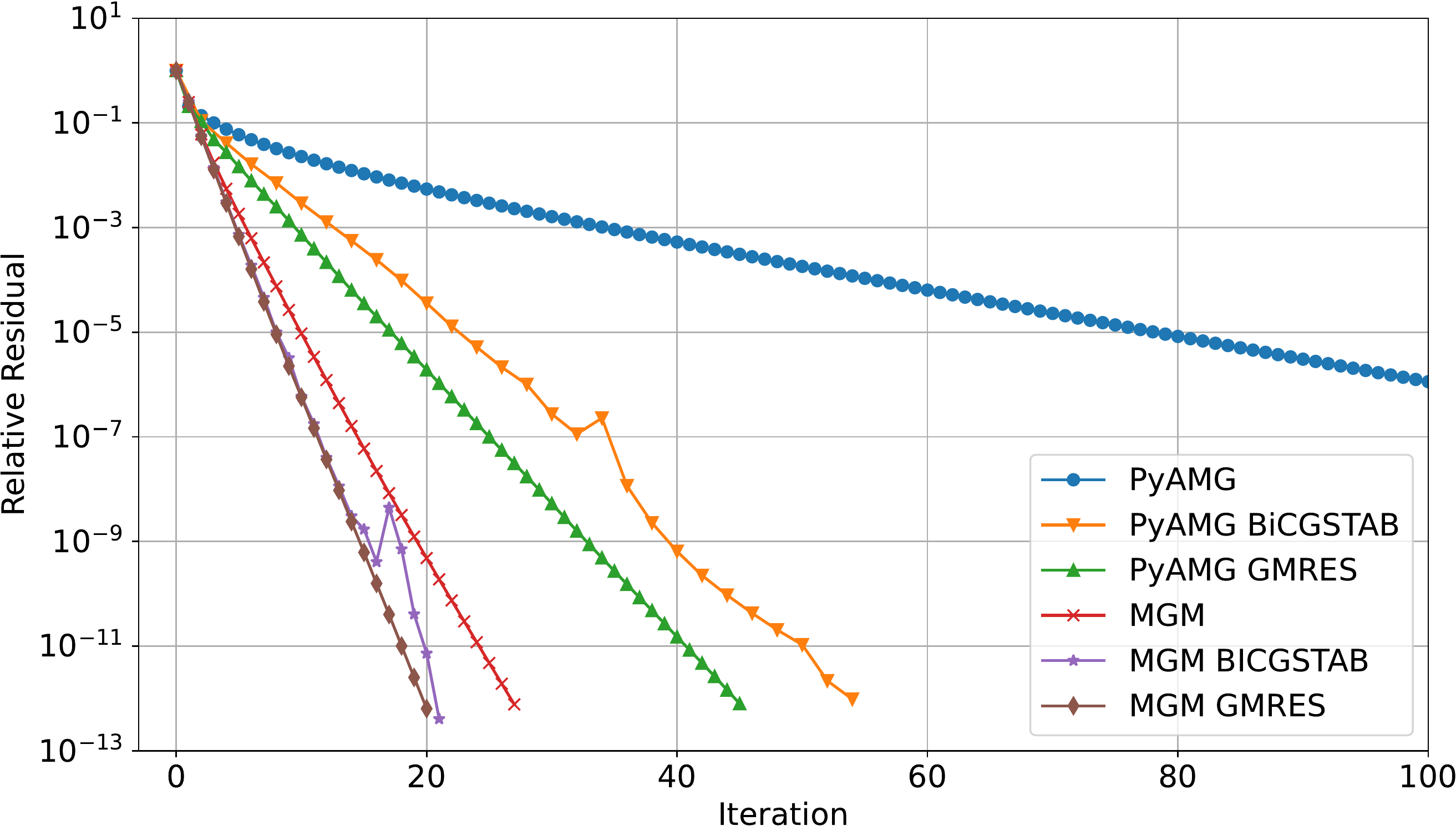} \\
\rotatebox{90}{\hspace{0.1\textwidth} \small $\ell=5$} & \includegraphics[width=0.45\textwidth]{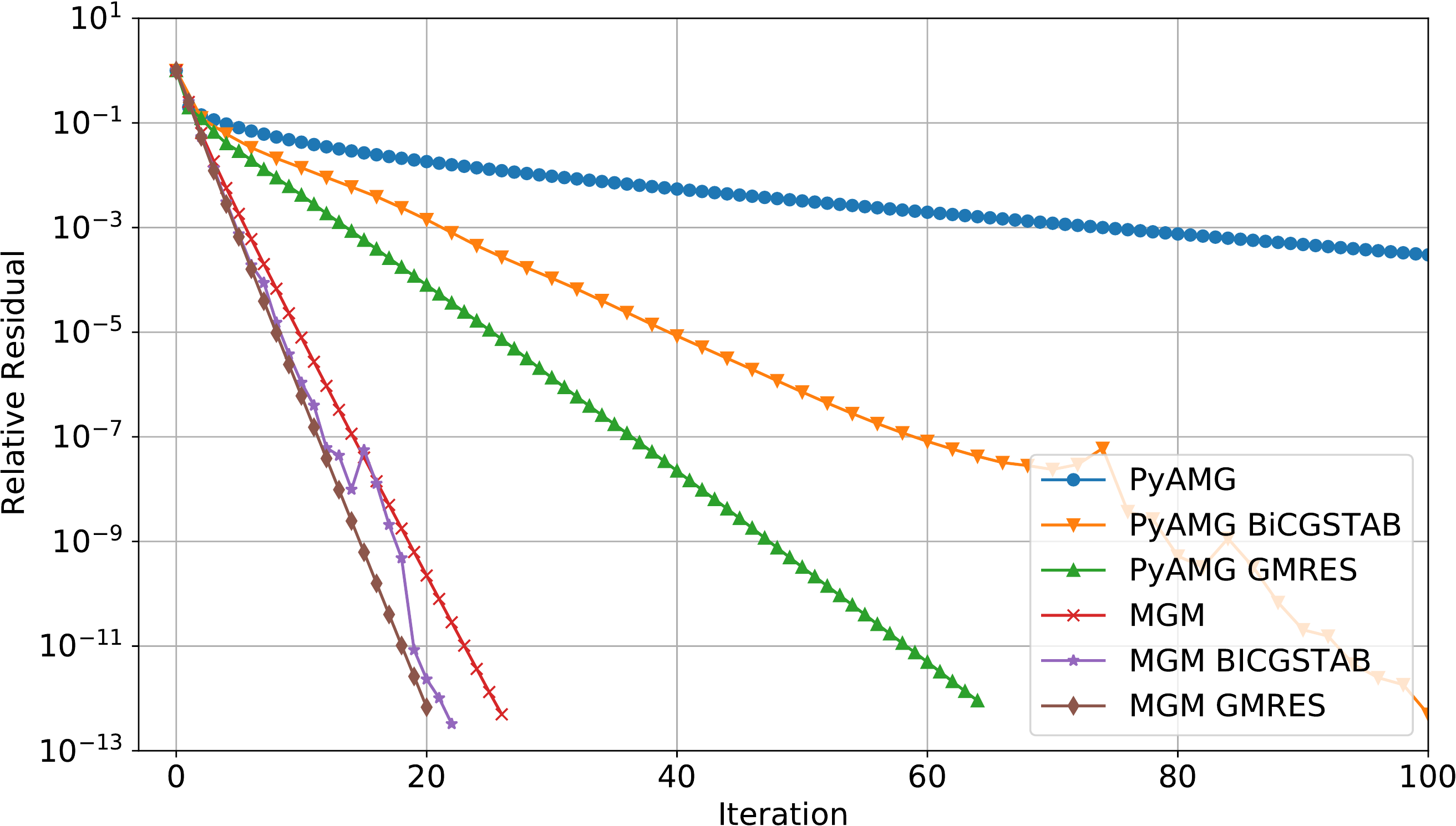} & 
\includegraphics[width=0.45\textwidth]{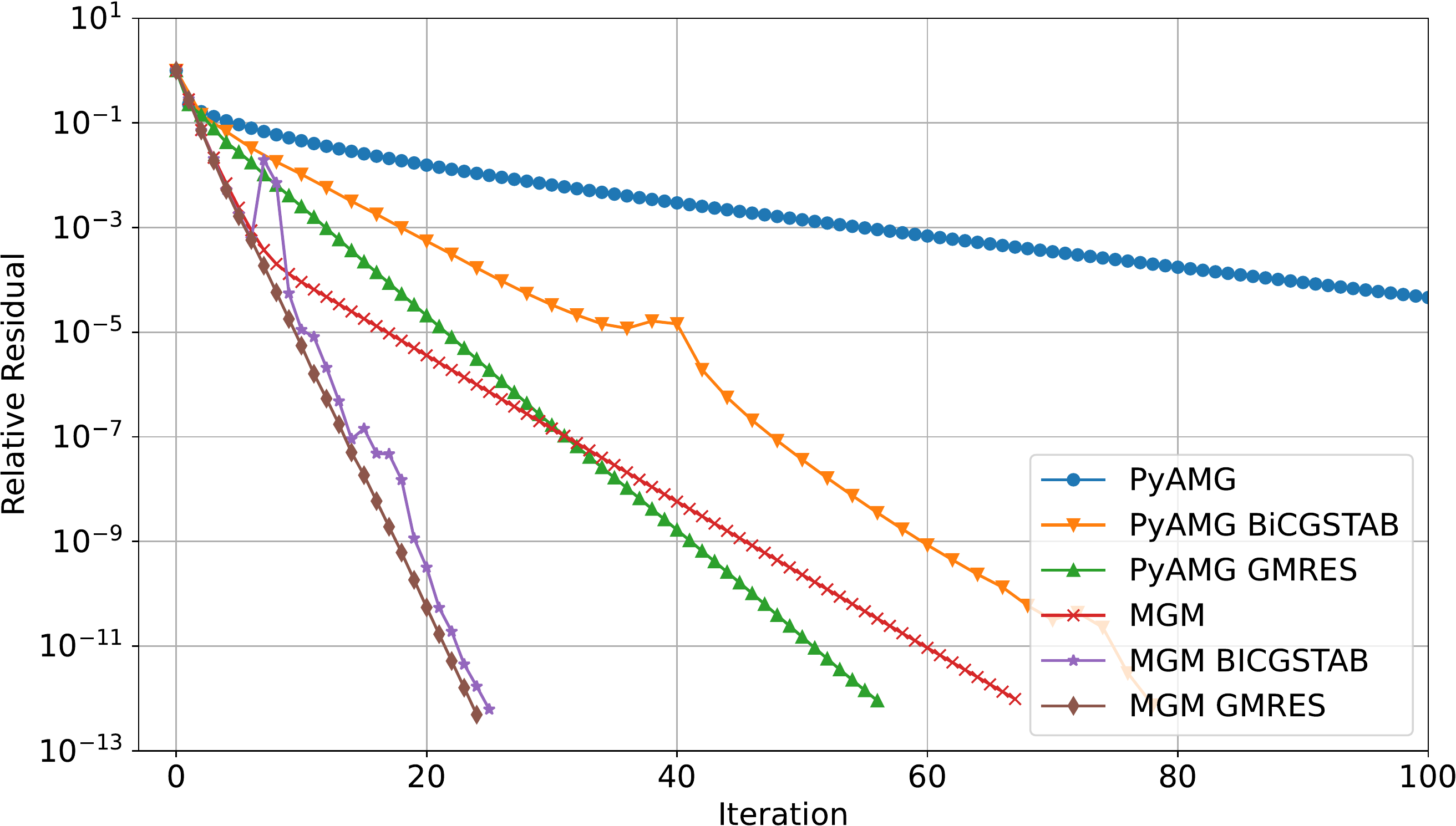} \\ 
\rotatebox{90}{\hspace{0.1\textwidth} \small $\ell=7$} & \includegraphics[width=0.45\textwidth]{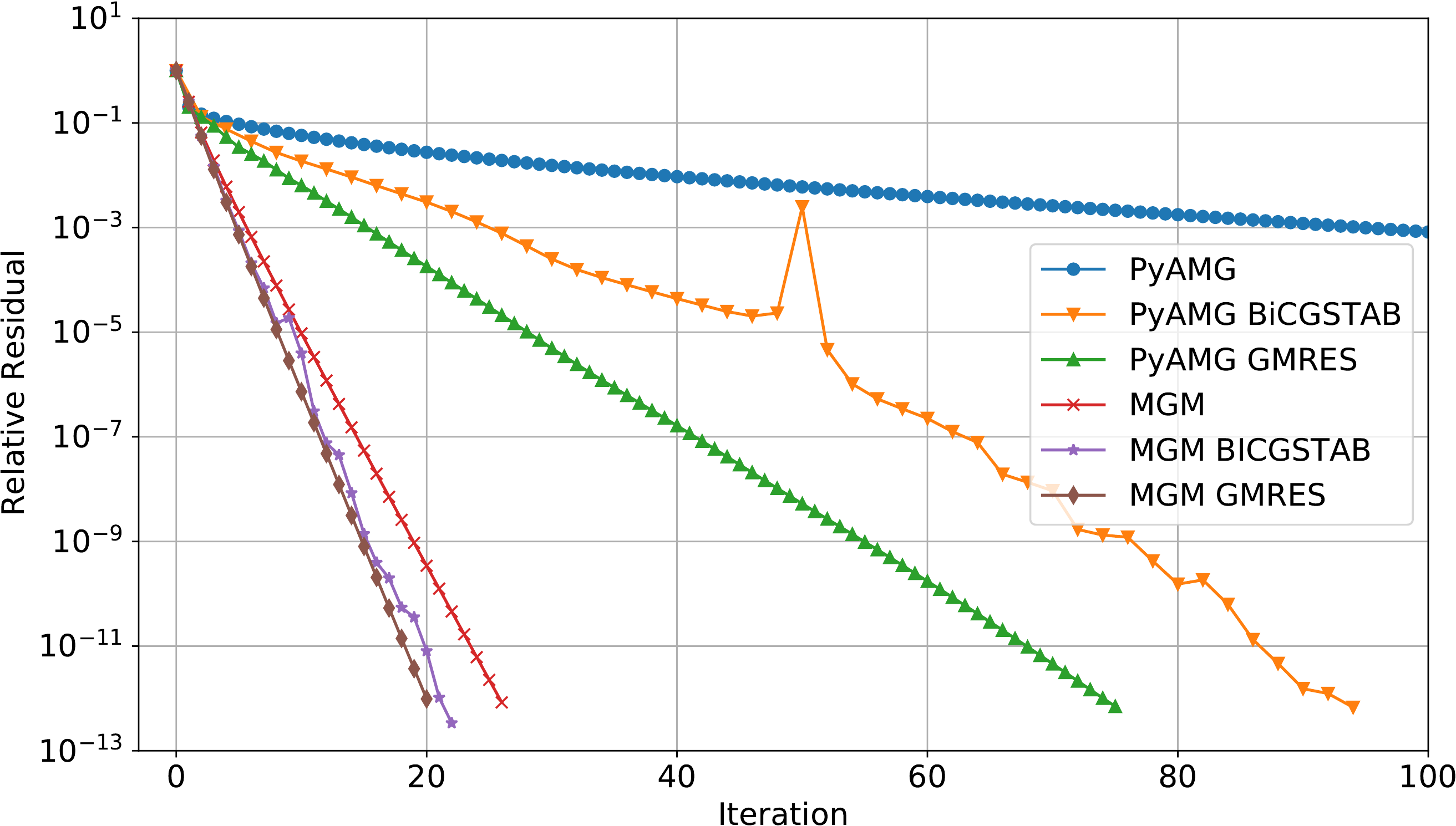} & 
\includegraphics[width=0.45\textwidth]{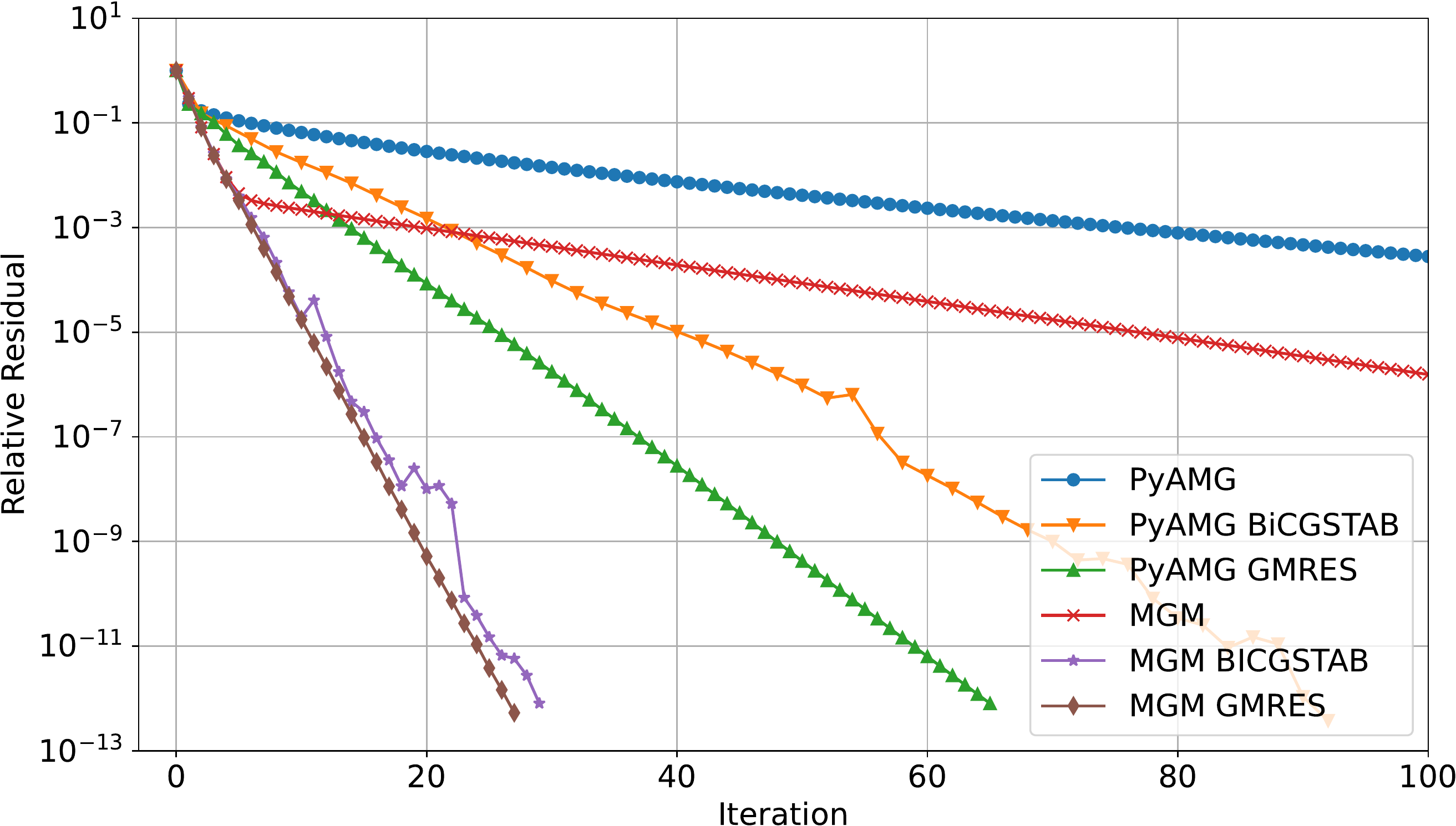} 
\end{tabular}
\caption{Convergence results for MGM and PyAMG based solvers for RBF-FD discretizations of a shifted Poisson problem with random right hand side.  The sphere results are for $N_h=2621442$, while for the cyclide $N_h = 2097152$.\label{fig:convergence_rbffd}}
\end{figure}

\subsection{Standalone solver vs.\ preconditioner}

In the first set of tests, we compare MGM and PyAMG both as standalone solvers and preconditioners.  For the latter approaches we refer to these solvers as MGM GMRES, MGM BiCGSTAB, PyAMG GMRES, and PyAMG BiCGSTAB, to indicate the type of Krylov method employed.  We use these solvers on the shifted Poisson problem on the unit sphere and cyclide with $N_h$=2,621,422 and $N_h$=2,097,152 nodes, respectively.  For the BiCGSTAB results, we count the number of applications of the preconditioner as the iterations since each step of this method applies the preconditioner twice, whereas GMRES applies it once.

\begin{figure}[htb]
\centering
\begin{tabular}{ccc}
& {\small Sphere} & {\small Cyclide} \\
\rotatebox{90}{\hspace{0.1\textwidth}\small $\ell=3$} & \includegraphics[width=0.45\textwidth]{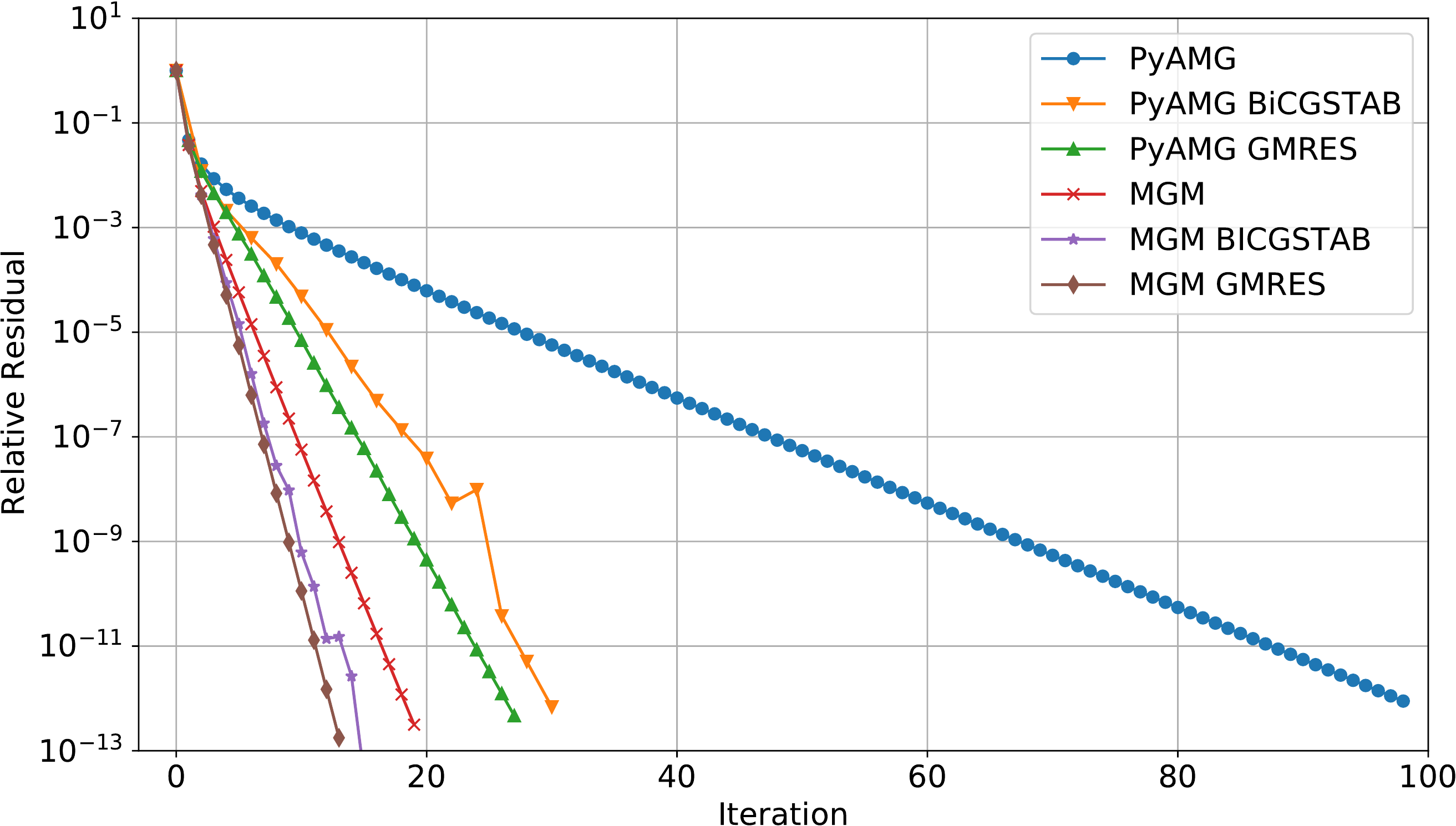} &
\includegraphics[width=0.45\textwidth]{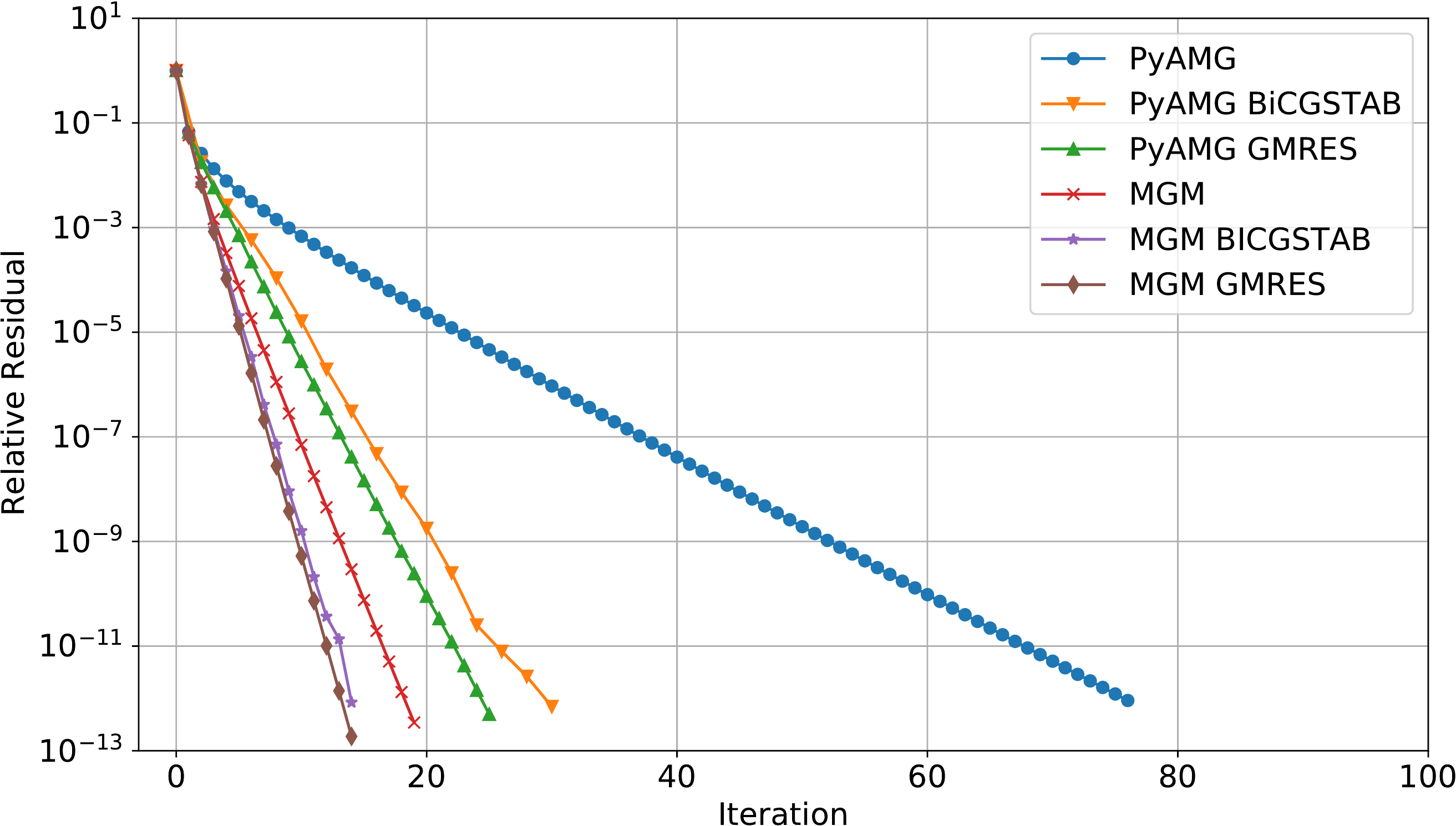} \\
\rotatebox{90}{\hspace{0.1\textwidth}\small $\ell=5$} & 
\includegraphics[width=0.45\textwidth]{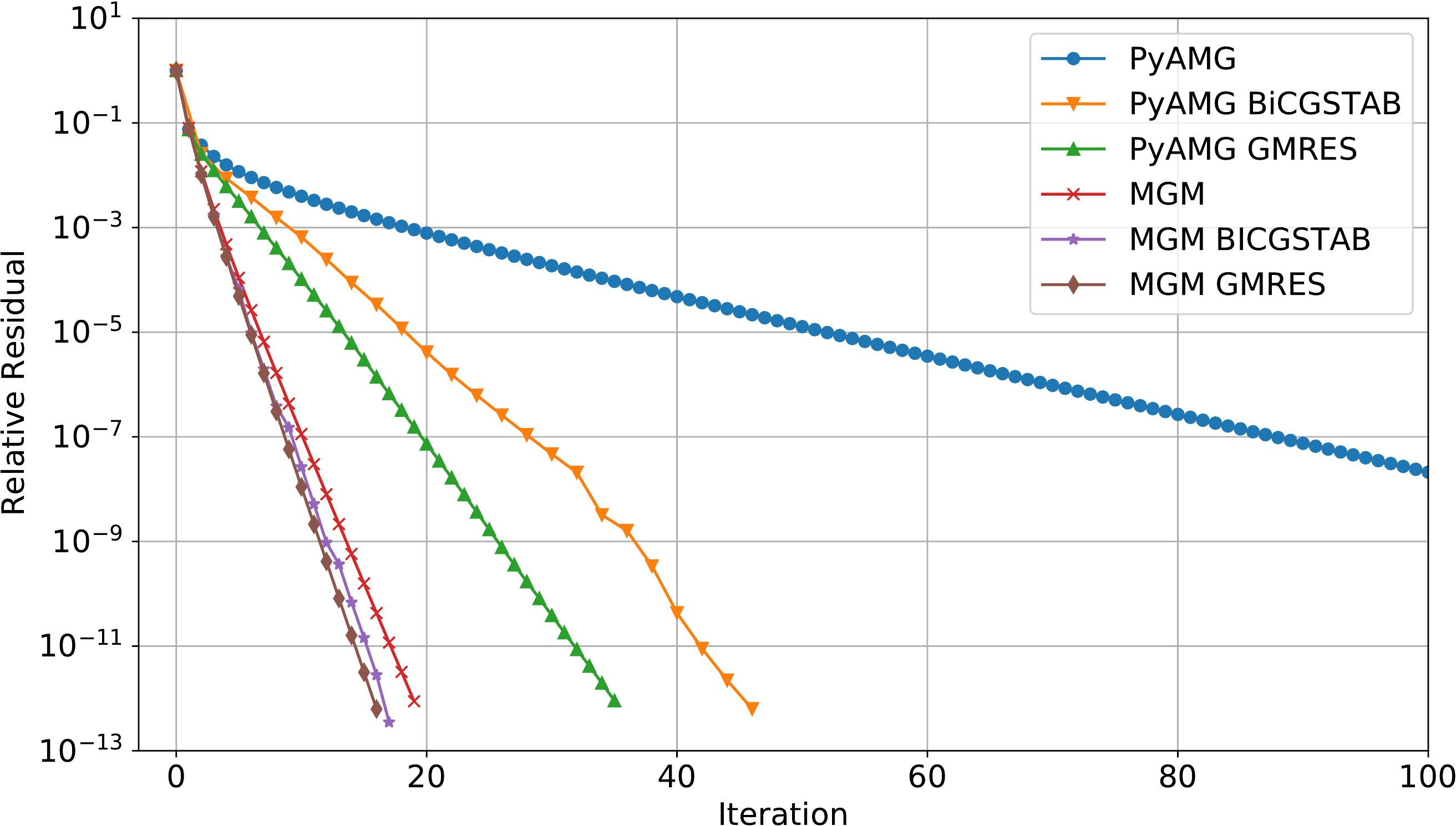} & 
\includegraphics[width=0.45\textwidth]{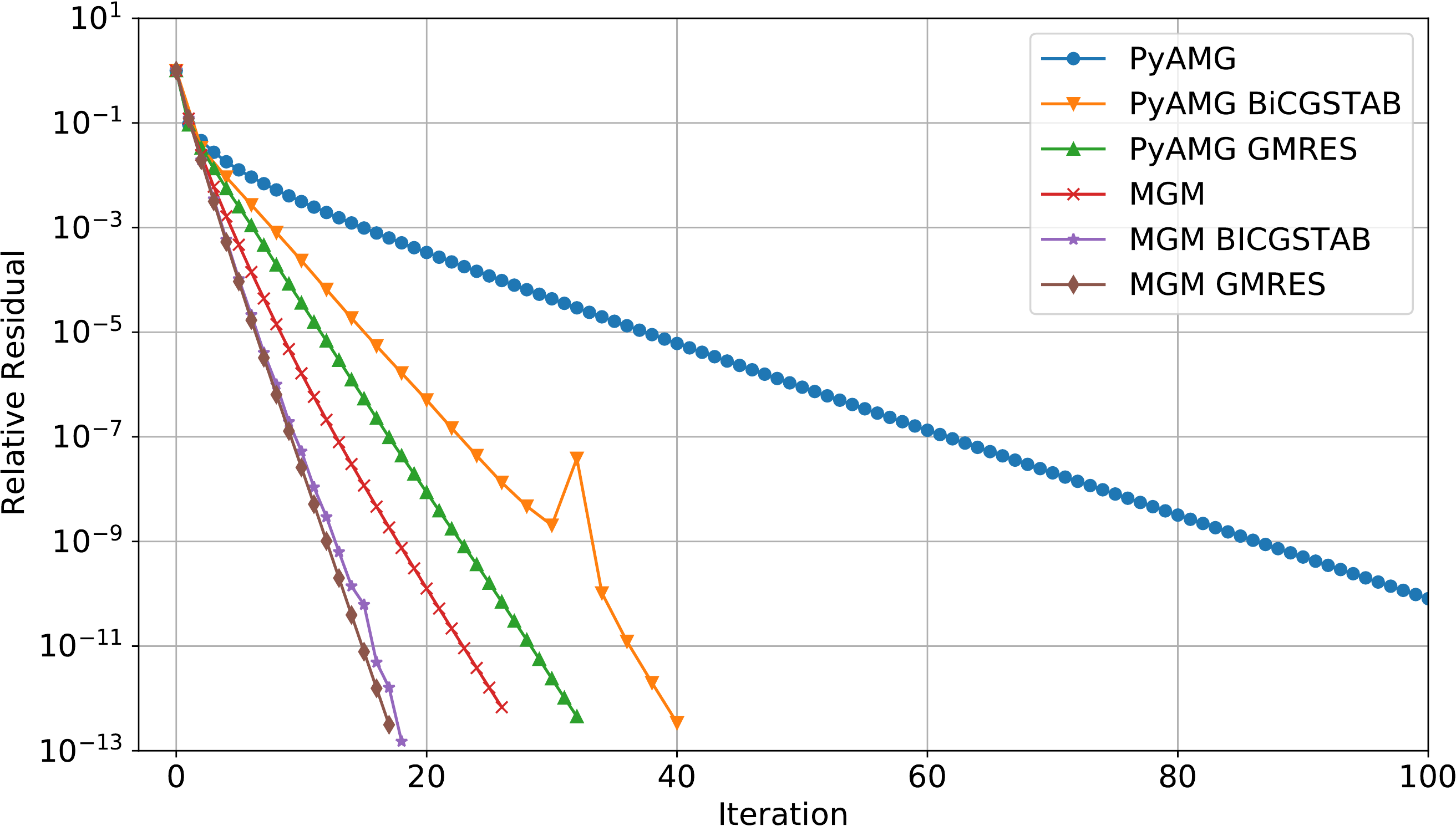} \\ 
\rotatebox{90}{\hspace{0.1\textwidth}\small $\ell=7$} & 
\includegraphics[width=0.45\textwidth]{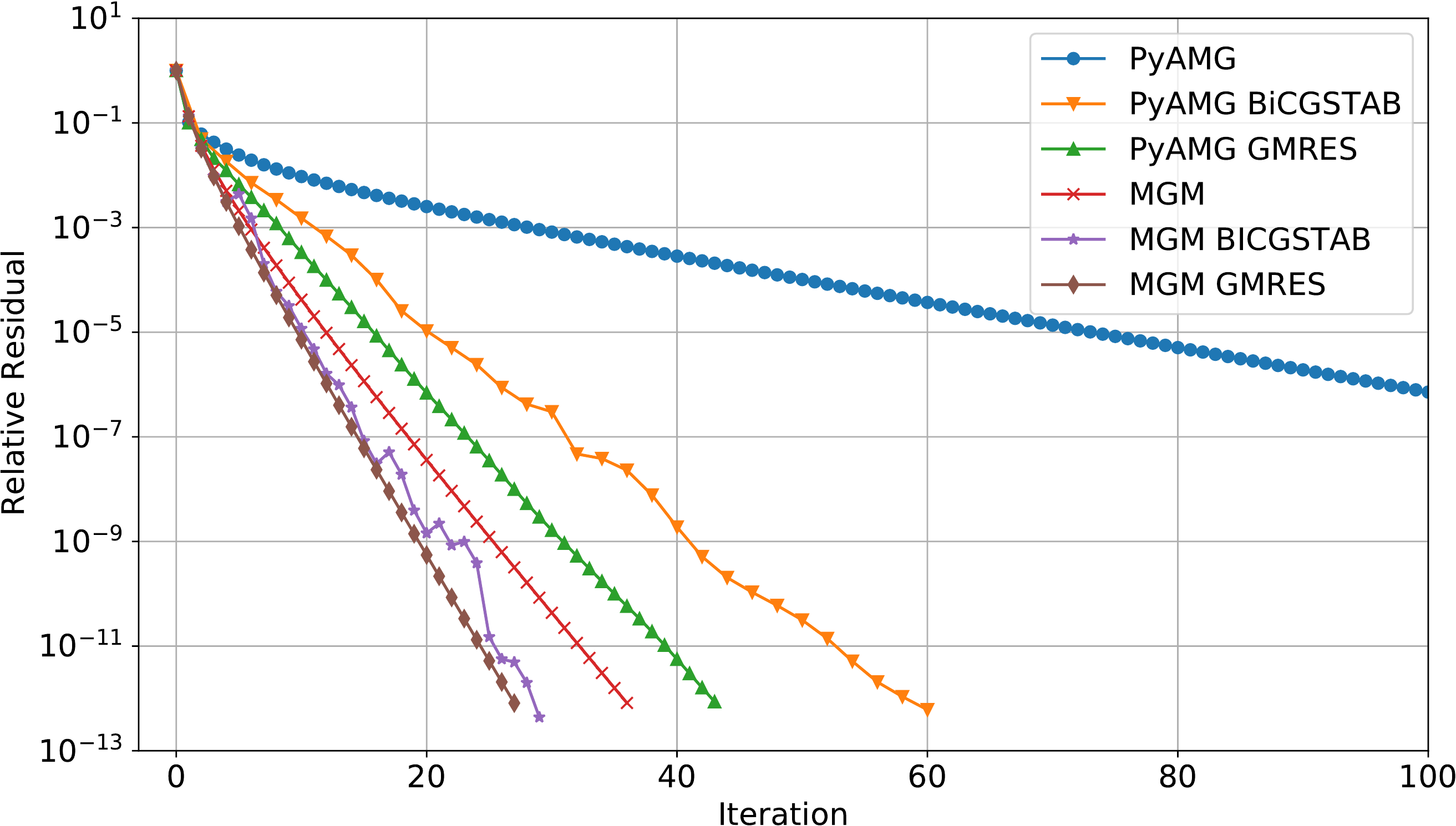} & 
\includegraphics[width=0.45\textwidth]{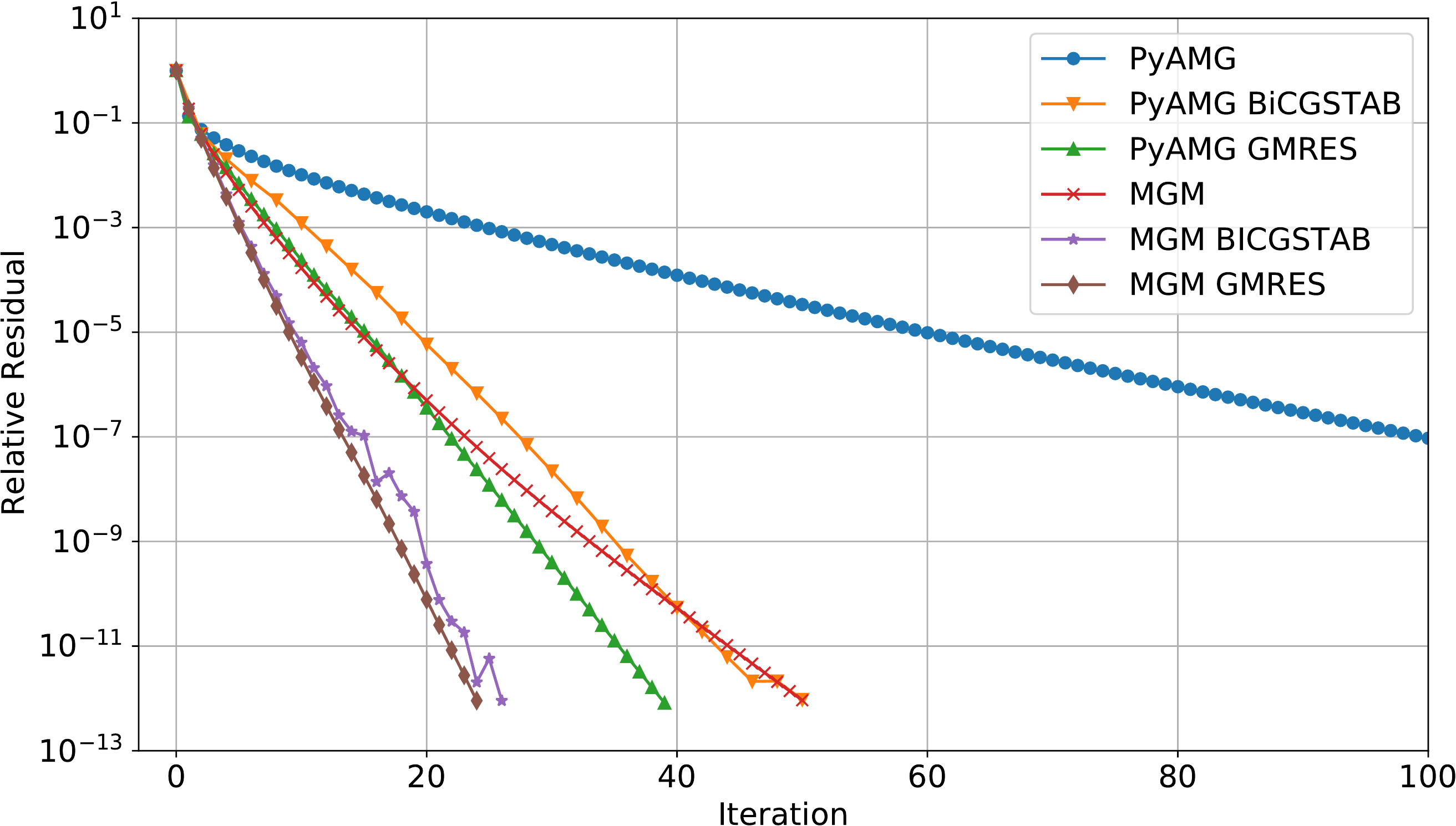}
\end{tabular}
\caption{Same as Figure \ref{fig:convergence_rbffd}, but for GFD discretizations of the shifted surface Poisson problem.\label{fig:convergence_gfd}}
\end{figure}

Figure \ref{fig:convergence_rbffd} displays the results in terms of relative residual vs.\ iteration count for RBF-FD, while Figure \ref{fig:convergence_gfd} displays the results for GFD.  For the RBF-FD results, we see that the methods using MGM converge more rapidly than the methods based on PyAMG for both surfaces.  For the sphere, MGM works very well as a standalone solver and preconditioner even as $\ell$ increases, but for the cyclide the convergence rates of MGM as a standalone solver decrease considerably.  This may be due to the more irregular nature of the cyclide point cloud.  We note, however, that the preconditioned versions of MGM only have a very mild decrease in convergence rates for the cyclide.  The figures also show that the methods using PyAMG do not converge as rapidly as the corresponding MGM methods, with the fastest converging PyAMG method taking more than double the number of iterations as the fastest MGM method when $\ell=3$ and triple when $\ell=5$ and $7$.  We see similar patterns in the GFD results, but the methods based on both MGM and PyAMG generally converge faster in this case and the gap between the fastest converging MGM and PyAMG methods is not as wide.  Finally, we note that MGM BiCGSTAB seems to converge at a very similar rate to MGM GMRES, whereas this does not hold for PyAMG.  This is a promising result for large systems since the storage requirements of BiCGSTAB are fixed, whereas they grow with the size of the Krylov subspace for GMRES~\cite{GhaiEtAl2018}.

These experiments also indicate that, while MGM can be an effective standalone solver for small $\ell$ (lower order discretizations), it is more robust for larger $\ell$ (higher order discretizations) and when used as a preconditioner.  This also seems to be the case when applying it to different surfaces and point clouds based on regular nodes (like the sphere) and irregular nodes (like the cyclide).  From the PyAMG results, it is clear that it should be used as a preconditioner to get the most robust results, which is generally the case for AMG methods applied to nonsymmetric systems~\cite{GhaiEtAl2018}.

\subsection{Scaling with problem size}

\begin{table}[tbh]
\centering
\begin{tabular}{ |c||c|c|c|c|c|c|} 
\hline
\multicolumn{7}{|c|}{Sphere} \\
 \hline
 & \multicolumn{3}{c|}{PyAMG} & \multicolumn{3}{c|}{MGM} \\
 $N$ & $\ell=3$ & $\ell=5$ & $\ell=7$ & $\ell=3$ & $\ell=5$ & $\ell=7$ \\
\hline
10242 & 30 (38) & 39 (56) & 46 (70) & 18 (19) & 18 (18) & 18 (20) \\ 
\hline 
40962 & 35 (42) & 43 (54) & 53 (78) & 19 (19) & 19 (20) & 19 (20) \\ 
\hline 
163842 & 40 (50) & 49 (64) & 59 (84) & 19 (21) & 19 (20) & 19 (20) \\ 
\hline 
655362 & 45 (62) & 55 (80) & 68 (90) & 20 (20) & 20 (20) & 20 (21) \\ 
\hline 
2621442 & 52 (66) & 64 (102) & 75 (94) & 20 (21) & 20 (22) & 20 (22) \\
\hline 
\hline
\multicolumn{7}{|c|}{Cyclide} \\
 \hline
 & \multicolumn{3}{c|}{PyAMG} & \multicolumn{3}{c|}{MGM} \\
 $N$ & $\ell=3$ & $\ell=5$ & $\ell=7$ & $\ell=3$ & $\ell=5$ & $\ell=7$ \\
\hline
8192 & 25 (30) & 34 (48) & 42 (62) & 17 (19) & 19 (20) & 20 (22) \\ 
\hline 
32768 & 32 (40) & 39 (48) & 46 (64) & 19 (20) & 22 (22) & 24 (26) \\ 
\hline 
131072 & 35 (46) & 45 (66) & 54 (74) & 20 (21) & 23 (26) & 28 (31) \\ 
\hline 
524288 & 41 (50) & 49 (68) & 60 (78) & 21 (22) & 25 (27) & 30 (39) \\ 
\hline 
2097152 & 45 (54) & 56 (78) & 65 (90) & 20 (21) & 24 (25) & 27 (29) \\ 
\hline
\end{tabular}
\caption{Comparison of the number of PyAMG and MGM preconditioned GMRES/BiCGSTAB iterations required to reach a relative residual tolerance of $10^{-12}$ for solving the shifted Poisson problem with RBF-FD discretizations.  The numbers not in parenthesis are for GMRES, while the numbers in parenthesis are for BiCGSTAB.\label{tbl:iters_vs_N_rbffd}}
\end{table}
In the next set of tests, we examine how both the MGM and PyAMG methods scale as the size of the point clouds $N_h$ increases.  We focus on the preconditioned versions of these methods and test them again on the sphere and cyclide.  
Tables \ref{tbl:iters_vs_N_rbffd} and \ref{tbl:iters_vs_N_gfd} display the results for the RBF-FD and GFD methods, respectively, in terms of number of iterations required to reach a relative residual of $10^{-12}$.  We see from these tables that the preconditioned MGM methods appear to scale much better than the PyAMG methods, both in terms of $N_h$ and $\ell$.  For RBF-FD discretizations, the increase in the iteration count for MGM is more mild with increasing $\ell$ than for GFD.  However, in all cases but $\ell=7$ on the sphere, the iteration count is lower for the GFD discretizations;  we examine this further in Section \ref{sec:accuracy}.

\begin{table}[tbh]
\centering
\begin{tabular}{ |c||c|c|c|c|c|c|} 
\hline
\multicolumn{7}{|c|}{Sphere} \\
 \hline
 & \multicolumn{3}{c|}{PyAMG} & \multicolumn{3}{c|}{MGM} \\
 $N$ & $\ell=3$ & $\ell=5$ & $\ell=7$ & $\ell=3$ & $\ell=5$ & $\ell=7$ \\
\hline
10242 & 15 (18) & 20 (26) & 26 (34) & 10 (10) & 14 (14) & 20 (21) \\ 
\hline 
40962 & 17 (20) & 23 (26) & 30 (42) & 10 (11) & 15 (16) & 23 (25) \\ 
\hline 
163842 & 20 (22) & 27 (36) & 34 (42) & 11 (12) & 15 (17) & 26 (29) \\ 
\hline 
655362 & 23 (26) & 30 (42) & 39 (48) & 12 (13) & 16 (17) & 27 (28) \\ 
\hline 
2621442 & 27 (30) & 35 (46) & 43 (60) & 13 (15) & 16 (17) & 27 (29) \\ 
\hline 
\hline
\multicolumn{7}{|c|}{Cyclide} \\
 \hline
 & \multicolumn{3}{c|}{PyAMG} & \multicolumn{3}{c|}{MGM} \\
 $N$ & $\ell=3$ & $\ell=5$ & $\ell=7$ & $\ell=3$ & $\ell=5$ & $\ell=7$ \\
\hline
8192 & 13 (16) & 17 (20) & 23 (28) & 10 (11) & 14 (14) & 18 (19) \\ 
\hline 
32768 & 16 (20) & 21 (24) & 26 (30) & 12 (12) & 17 (18) & 24 (26) \\ 
\hline 
131072 & 19 (24) & 24 (32) & 30 (36) & 12 (13) & 18 (21) & 25 (27) \\ 
\hline 
524288 & 22 (26) & 27 (34) & 34 (40) & 13 (14) & 18 (18) & 26 (29) \\ 
\hline 
2097152 & 25 (30) & 32 (40) & 39 (50) & 14 (14) & 17 (18) & 24 (26) \\ 
\hline
\end{tabular}
\caption{Same as Table \ref{tbl:iters_vs_N_rbffd}, but for GFD discretizations.  For these results, we set $\alpha=4$ for all $N_h$, but the largest, where we set $\alpha=5$.\label{tbl:iters_vs_N_gfd}}
\end{table}

%\begin{table}[htb]\centering
%\begin{tabular}{ |c||c|c|c|c|c|c|} 
%\hline
% \multicolumn{7}{|c|}{Sphere} \\
%  \hline
% & \multicolumn{3}{c|}{PyAMG} & \multicolumn{3}{c|}{MGM} \\
% $N$ & $\ell=3$ & $\ell=5$ & $\ell=7$ & $\ell=3$ & $\ell=5$ & $\ell=7$ \\
%\hline
%10242  & 30 (15) &  39 (20) &  46 (26) & 18 (10) &  18 (14)  &  18 (20) \\
%\hline
% 40962  & 35 (17) &  43 (23) &  53 (30) &  19 (11) &  19 (15) &  19 (23) \\ 
%\hline
% 163842  & 40 (20) &  49 (27) &  59 (34) &  19 (11) &  20 (15) &  20 (25) \\ 
%\hline
% 655362  & 45 (23) &  55 (30) &  68 (39) &  20 (12) &  20 (16) &  20 (27) \\ 
%\hline
% 2621442  & 53 (26) &  63 (35) &  75 (44)  & 20 (13) & 20 (16) & 21 (27)\\
%\hline
% \multicolumn{7}{|c|}{Cyclide} \\
% \hline
% & \multicolumn{3}{c|}{PyAMG} & \multicolumn{3}{c|}{MGM} \\
% $N$ & $\ell=3$ & $\ell=5$ & $\ell=7$ & $\ell=3$ & $\ell=5$ & $\ell=7$ \\
%\hline
%8192  & 25 (13) &  34 (17) &  42 (23) &  17 (10) &  19 (14) &  20 (18) \\ 
%\hline
% 32768  & 31 (16) &  39 (21) &  46 (26) &  19 (12) &  21 (17) &  24 (24) \\ 
%\hline
% 131072  & 35 (19) &  45 (24) & 54 (30) &  20 (13) &  23 (18) &  28 (25) \\ 
%\hline
% 524288  & 40 (22) &  49 (27) &  60 (34) &  21 (13)  &  26 (18) &  30 (26) \\ 
%\hline
%2097152 & 45 (25) & 56 (33) & 65 (39) & 20 (14) & 23 (17) & 27 (24) \\
%\hline
%\end{tabular}
%\caption{Comparison of the number of PyAMG and MGM preconditioned GMRES iterations required to reach a relative residual tolerance of $10^{-12}$ for solving the screened Poisson equation on the sphere and cyclide.  }
%\end{table}

\begin{figure}[htb]
\centering
\begin{tabular}{ccc} 
& {\small Sphere} & {\small Cyclide} \\
\rotatebox{90}{\hspace{0.1\textwidth}\small $\ell=3$} & 
\includegraphics[width=0.4\textwidth,scale=0.7]{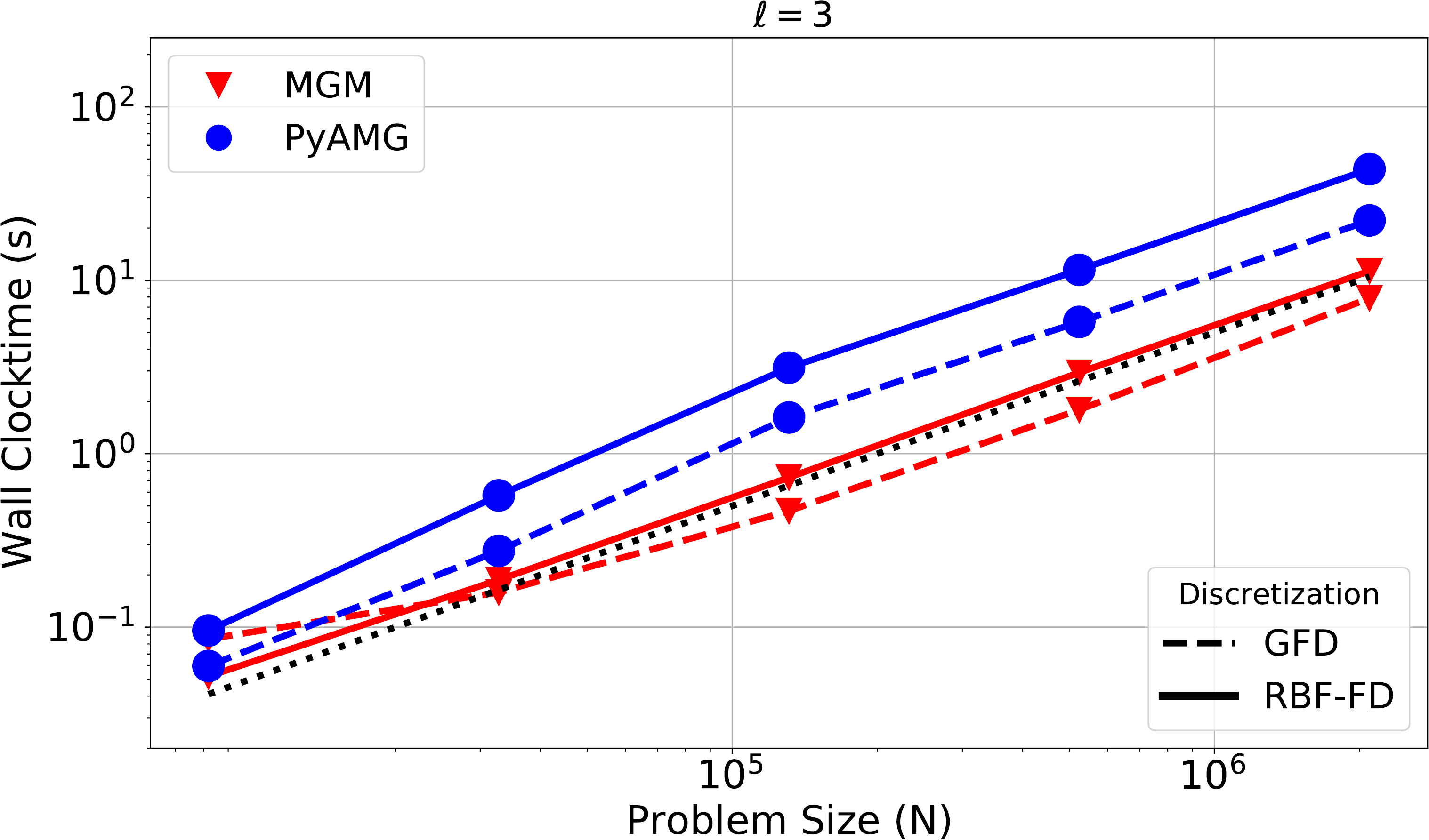}  & 
\includegraphics[width=0.4\textwidth,scale=0.7]{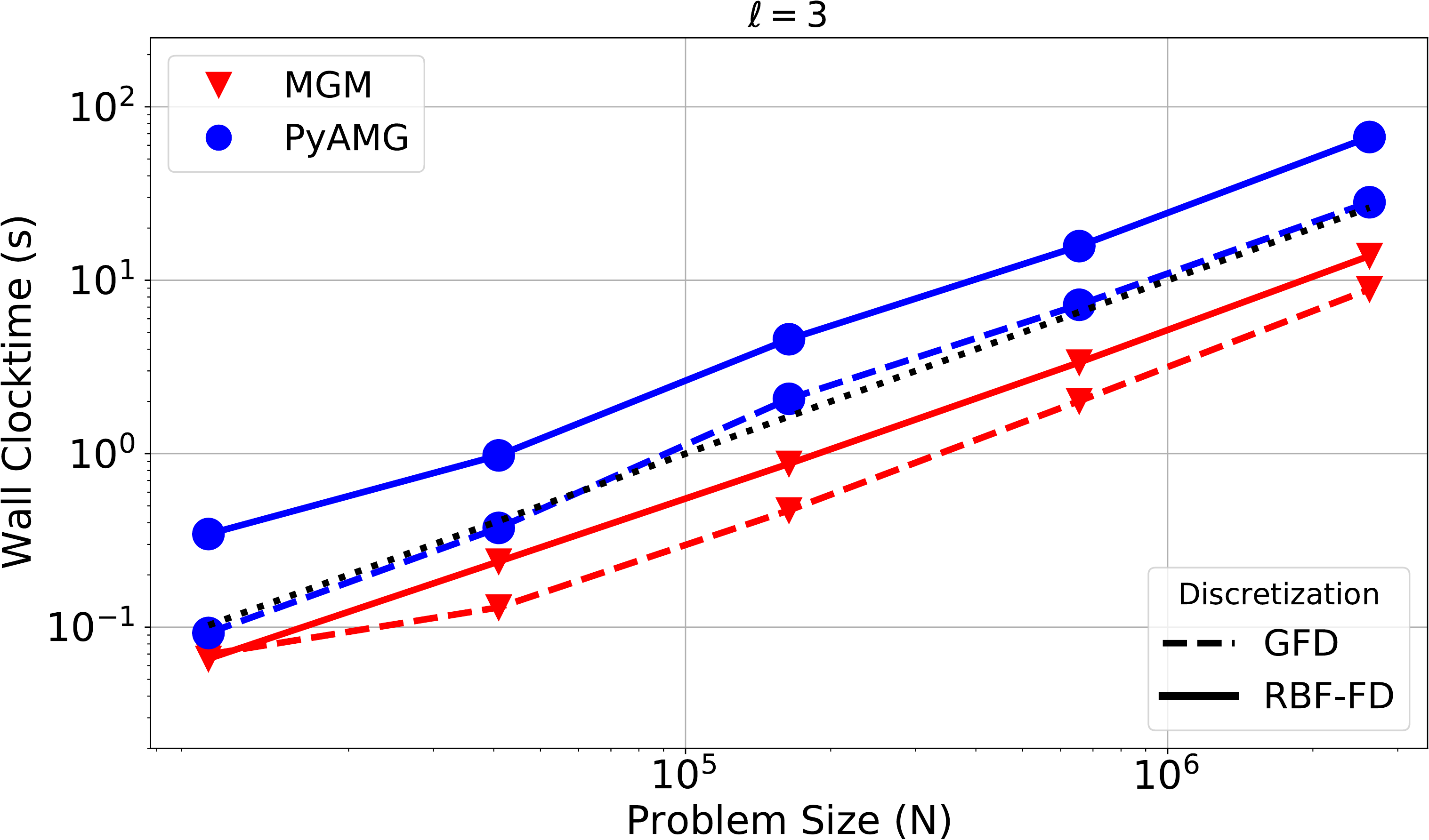} \\ 
\rotatebox{90}{\hspace{0.1\textwidth}\small $\ell=5$} & 
\includegraphics[width=0.4\textwidth,scale=0.7]{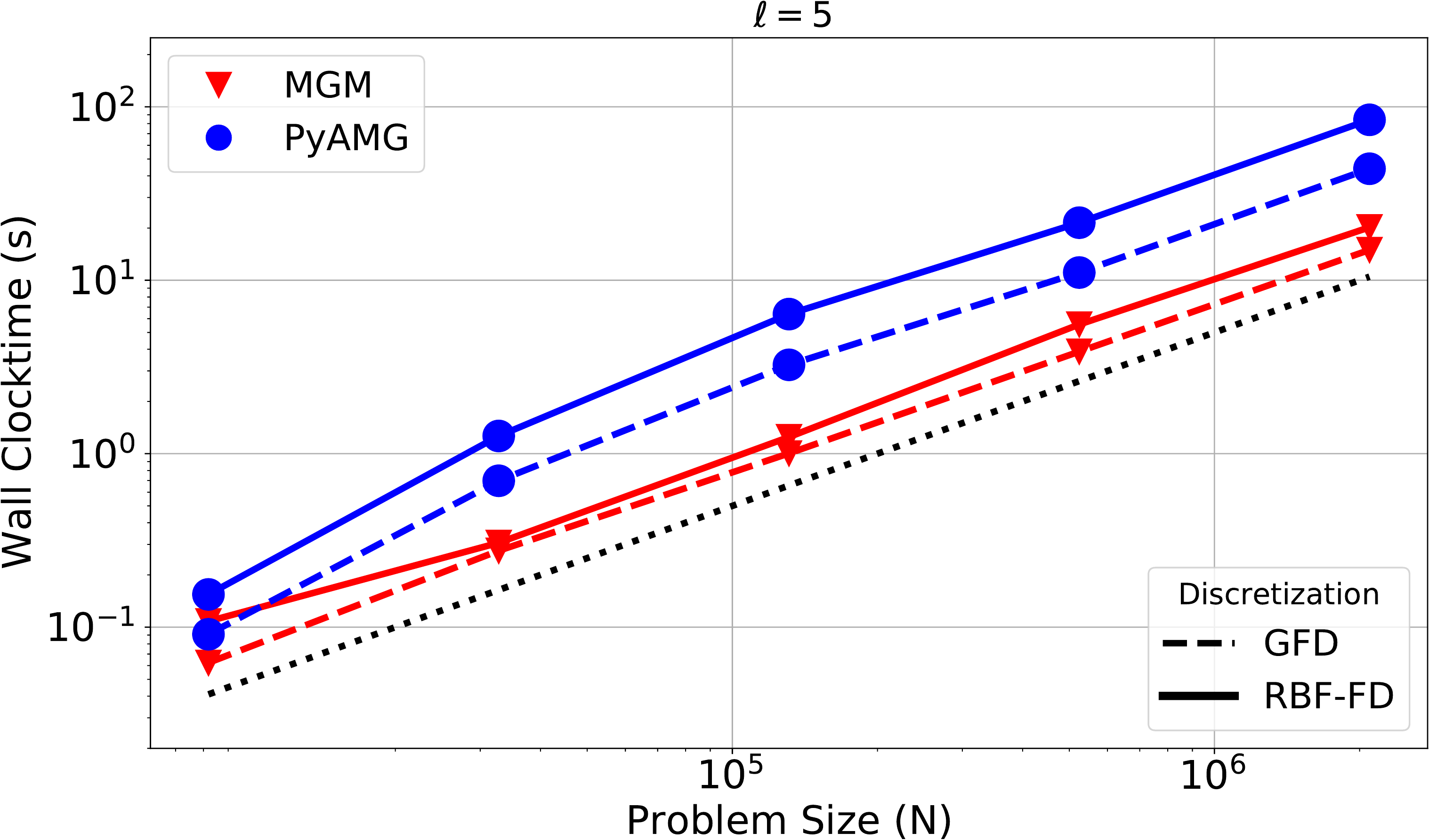}  & 
\includegraphics[width=0.4\textwidth,scale=0.7]{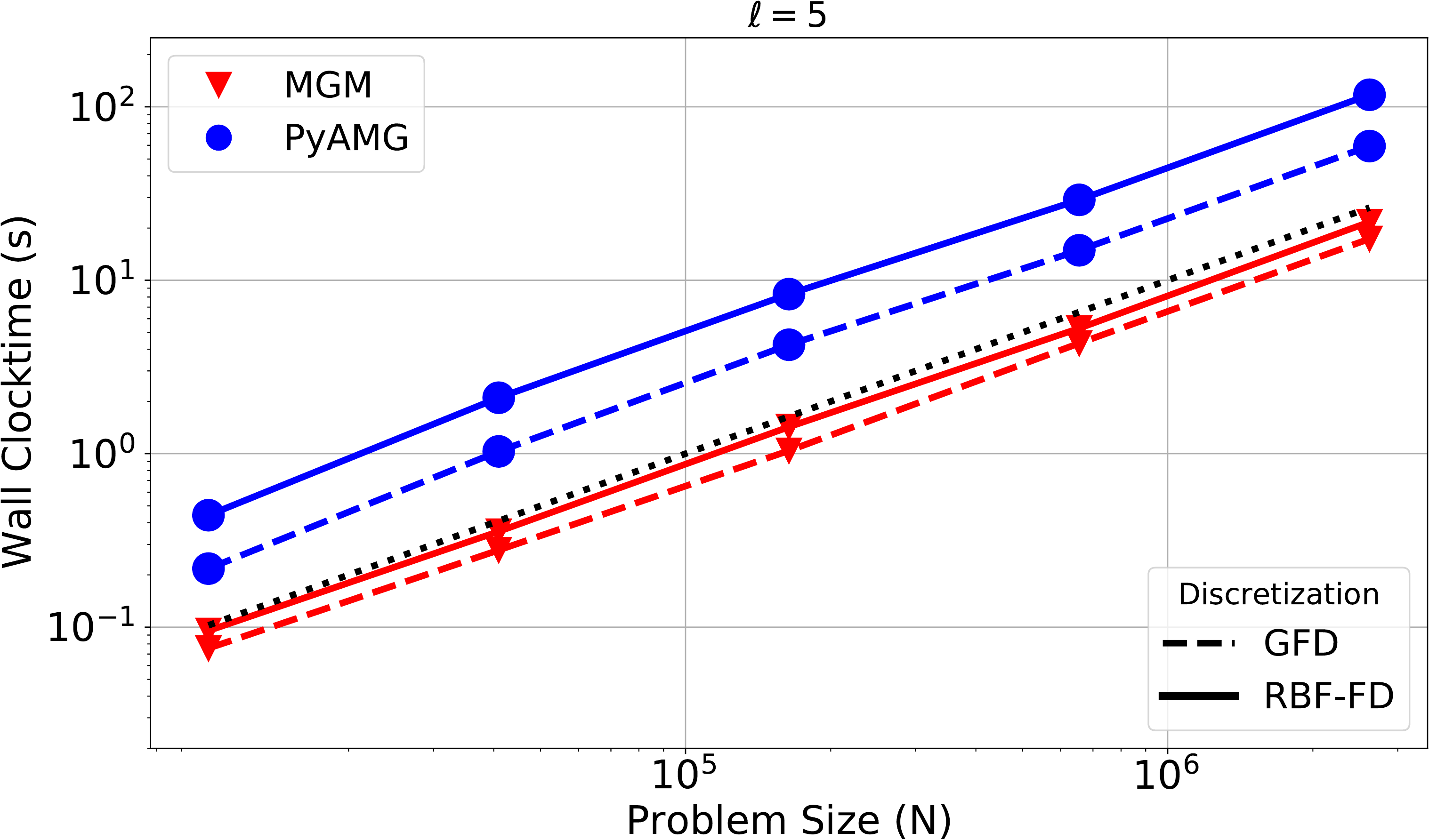} \\ 
\rotatebox{90}{\hspace{0.1\textwidth}\small $\ell=7$} & 
\includegraphics[width=0.4\textwidth,scale=0.7]{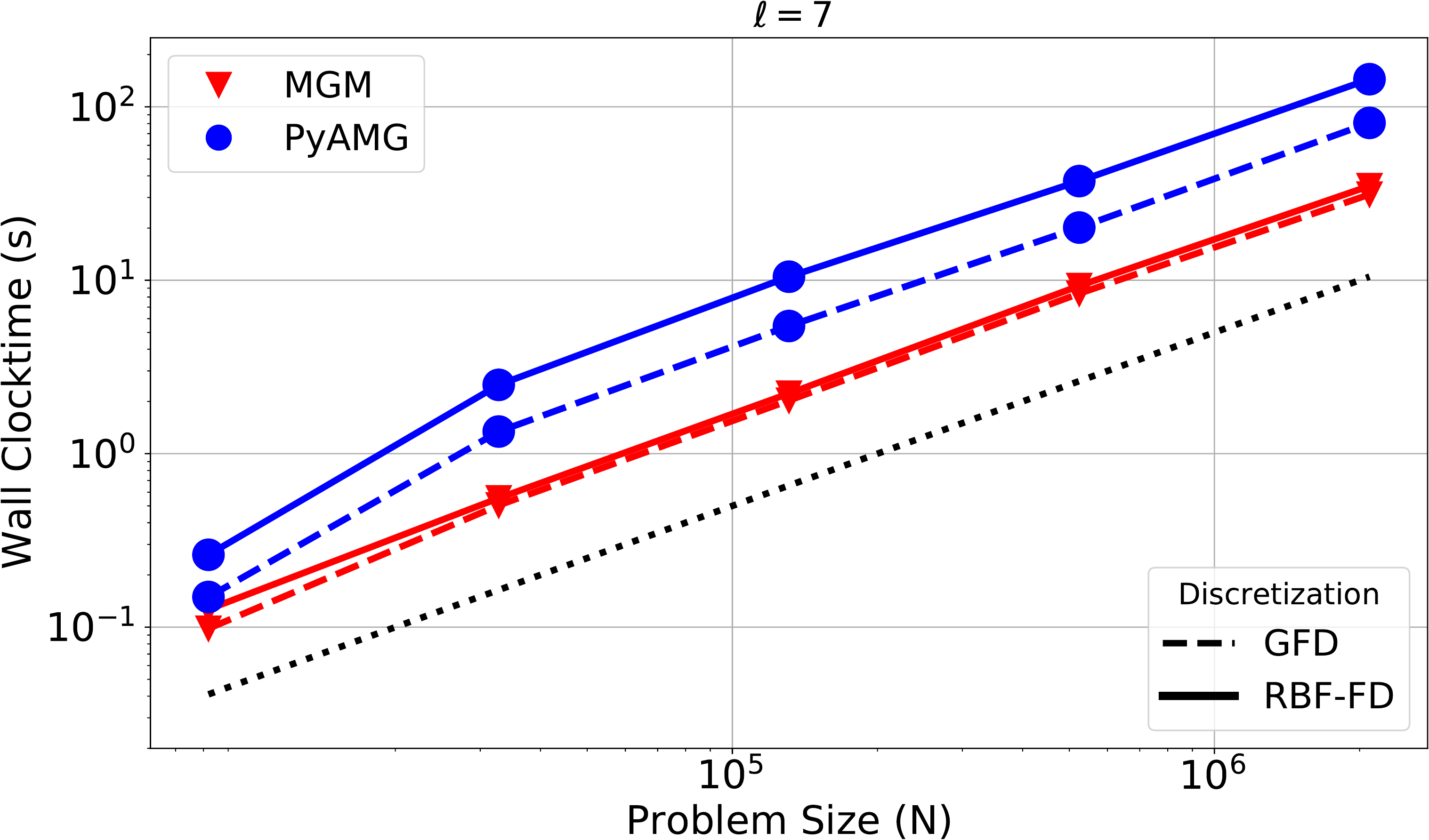}  & 
\includegraphics[width=0.4\textwidth,scale=0.7]{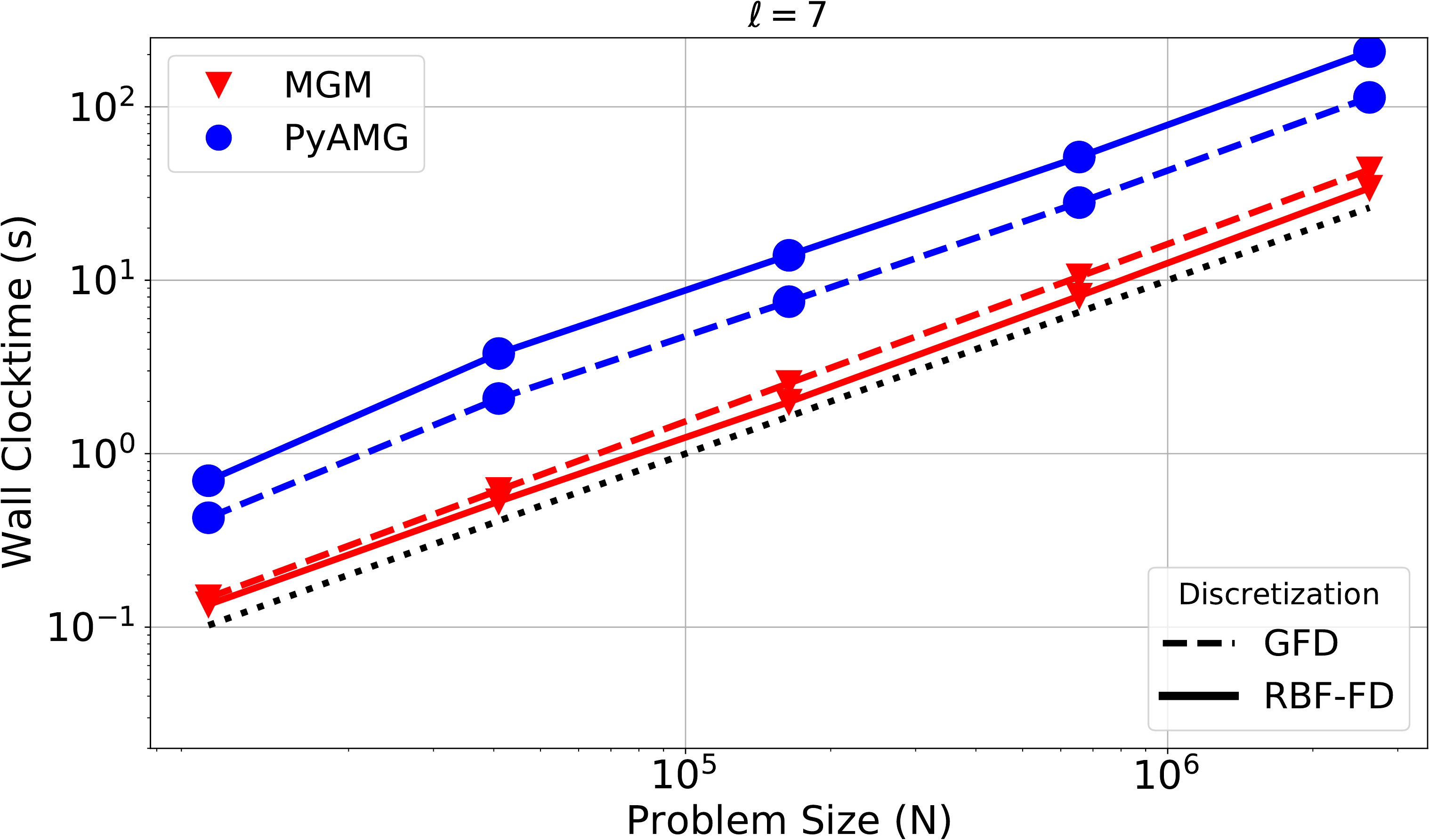} 
\end{tabular}
\caption{Wall-clock time (in seconds) for PyAMG GMRES and MGM GMRES to converge to a relative residual of $10^{-12}$ problem size ($N_n$) increases for RBF-FD (solid line) and GFD (dashed line) discretizations.  The black dotted line marks linear scaling, $\mathcal{O}(N_h)$, for reference.\label{fig:cost_geo}}
\end{figure}

In Figure \ref{fig:cost_geo} we display the wall-clock times for the results in Tables \ref{tbl:iters_vs_N_rbffd} and \ref{tbl:iters_vs_N_gfd} for GMRES PyAMG and MGM.  These results were run on a Linux Workstation with Intel i9-9900X 3.5 GHz processor (with no explicit parallelization) and do not include the set-up times.  We see from Figure \ref{fig:cost_geo} that MGM has a lower wall-clock time than PyAMG  for all but the first $N_h$ in the case of the sphere.  Furthermore, for the largest $N_h$, MGM is between 3 and 5 times faster.  Additionally, the dotted line in these scaling plots marks perfect linear scaling and we see that the results for both methods have a very similar slope to this line.  Finally, we note that the timing results for BiCGSTAB follow a similar trend to GMRES, so we omitted displaying the results.  However, the gap between the MGM and PyAMG timings were larger in this case. 

\subsection{Spectrum analysis}
\begin{figure}
\centering
\begin{tabular}{ccc}
& {\small Sphere} & {\small Cyclide} \\
\rotatebox{90}{\hspace{0.08\textwidth}\small $\ell=3$} & 
\includegraphics[width=0.45\textwidth]{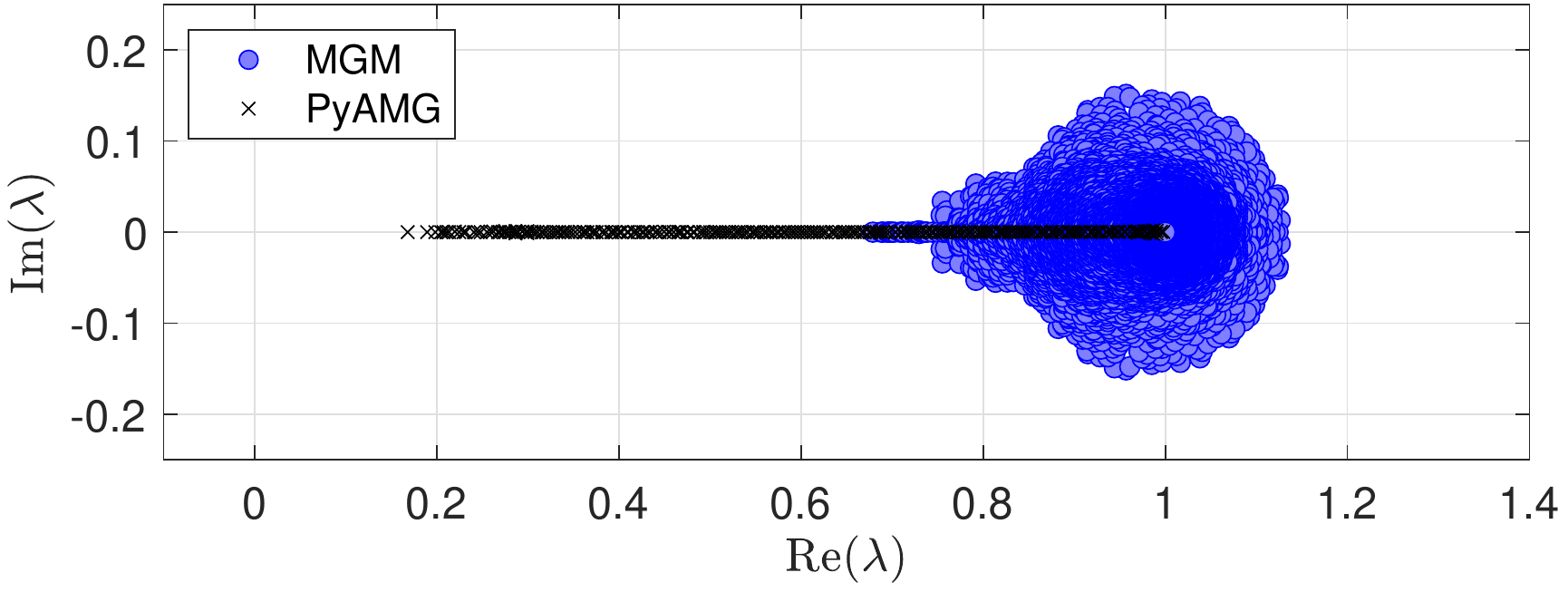} & \includegraphics[width=0.45\textwidth]{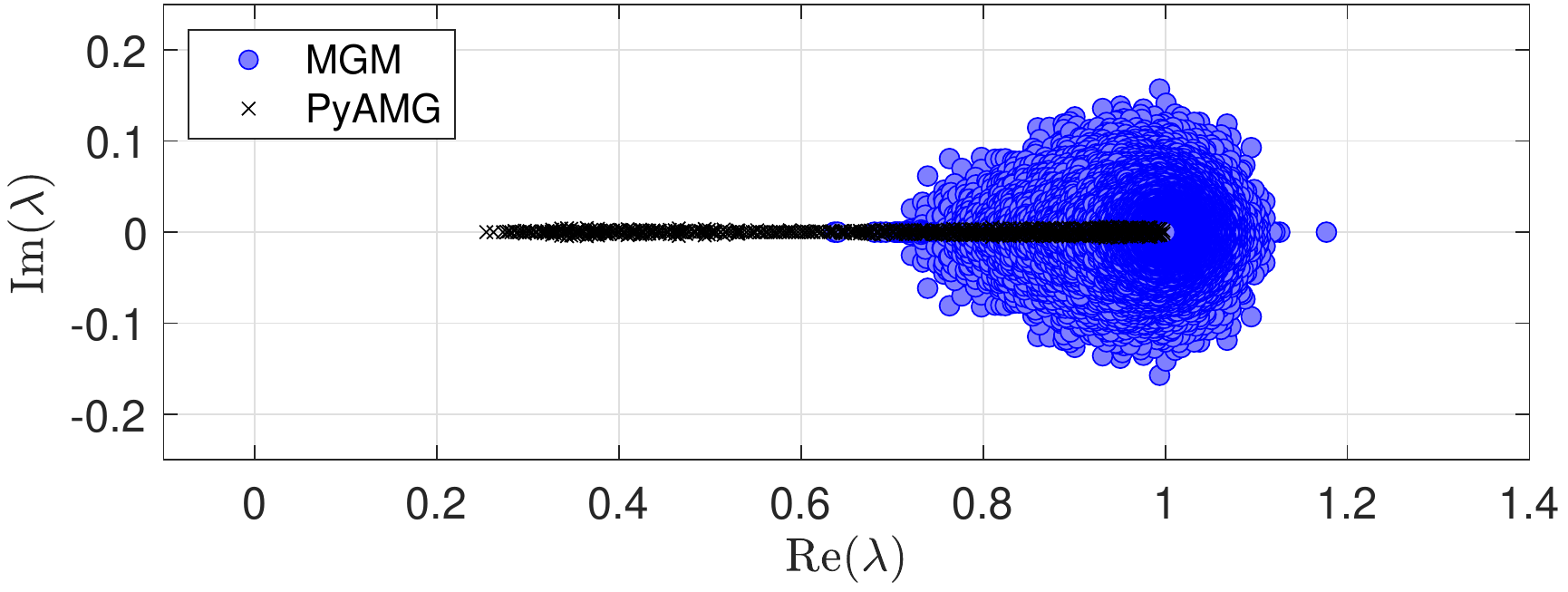} \\ 
\rotatebox{90}{\hspace{0.08\textwidth}\small $\ell=5$} & 
\includegraphics[width=0.45\textwidth]{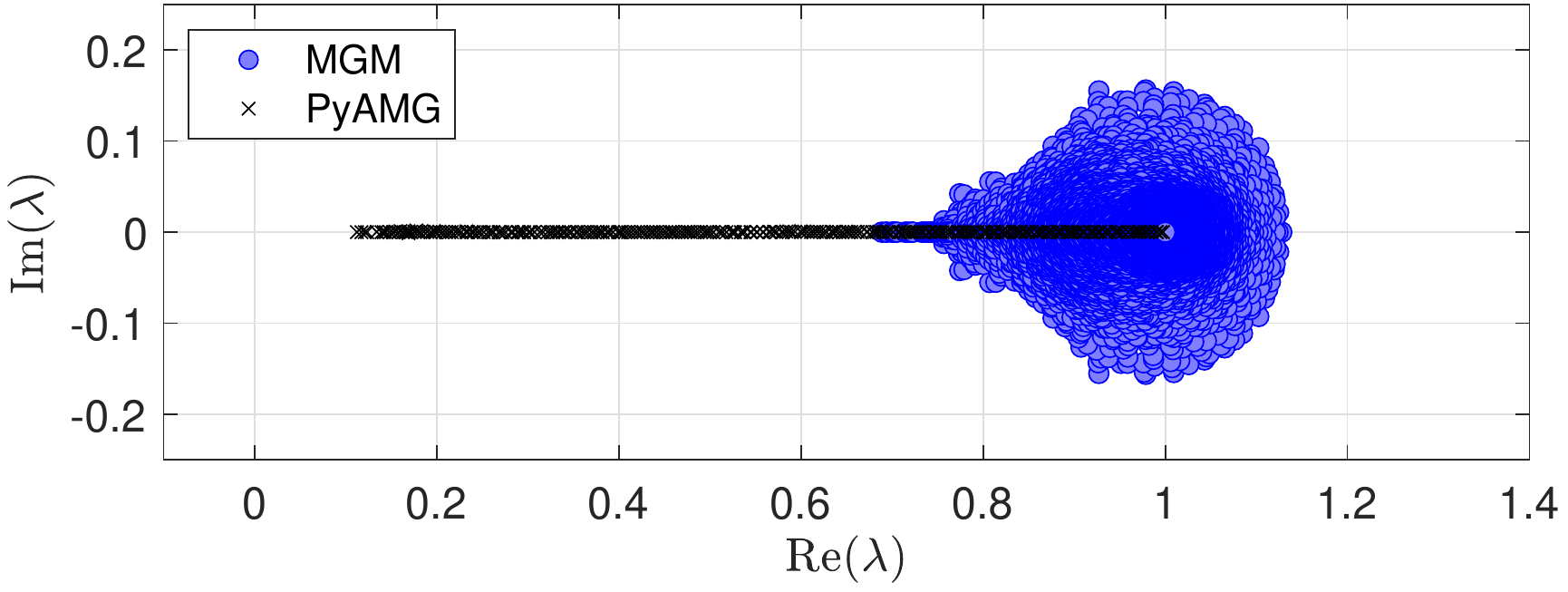} & \includegraphics[width=0.45\textwidth]{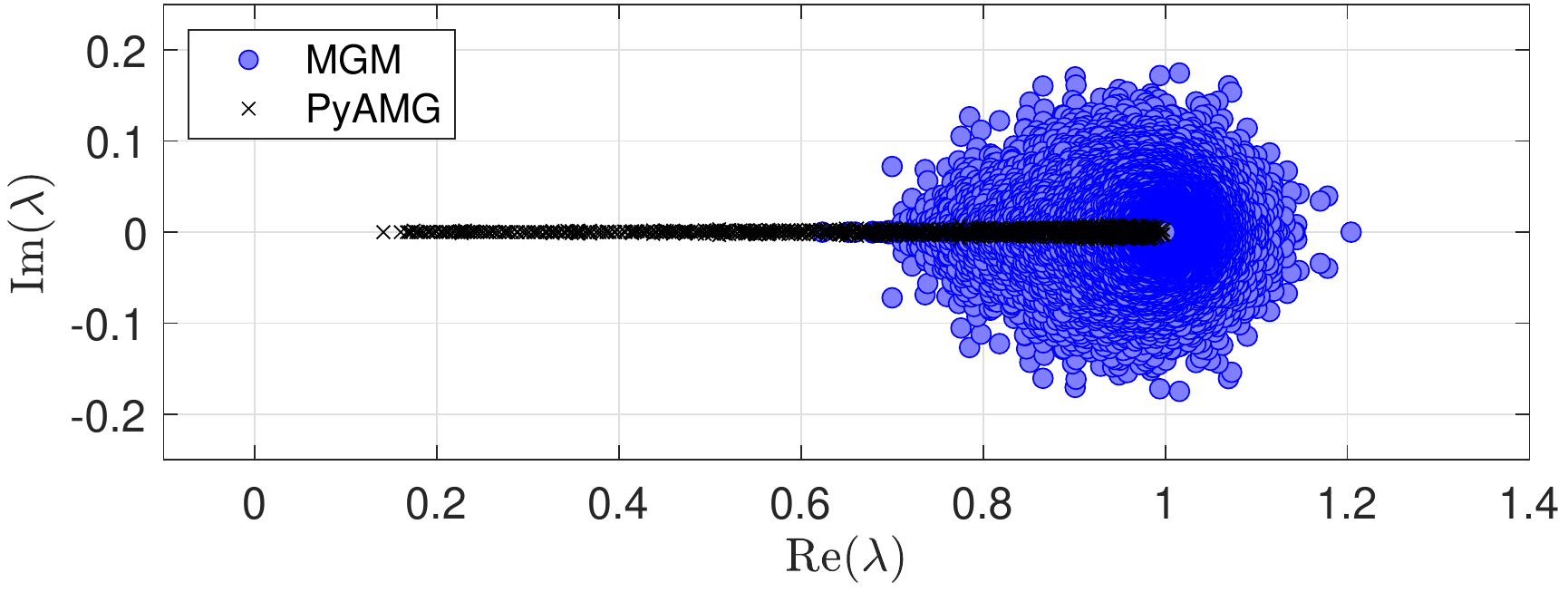} \\ 
\rotatebox{90}{\hspace{0.08\textwidth}\small $\ell=7$} & 
\includegraphics[width=0.45\textwidth]{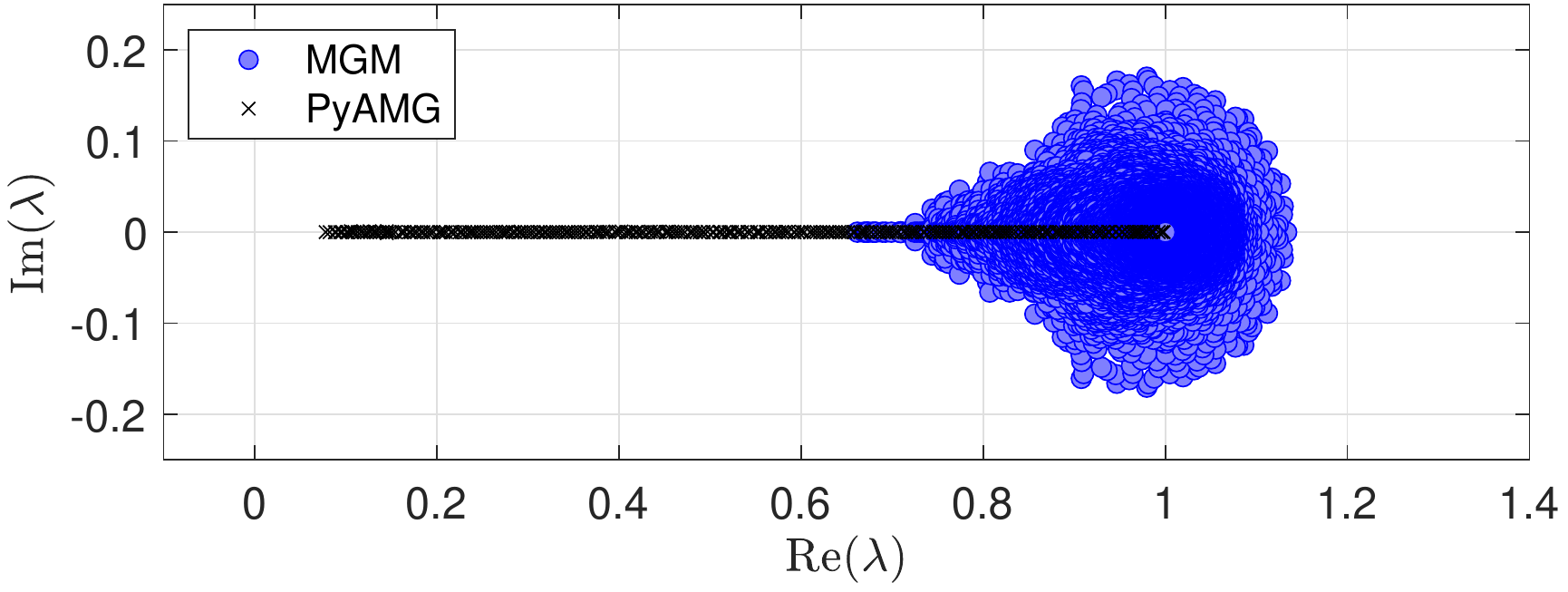} & \includegraphics[width=0.45\textwidth]{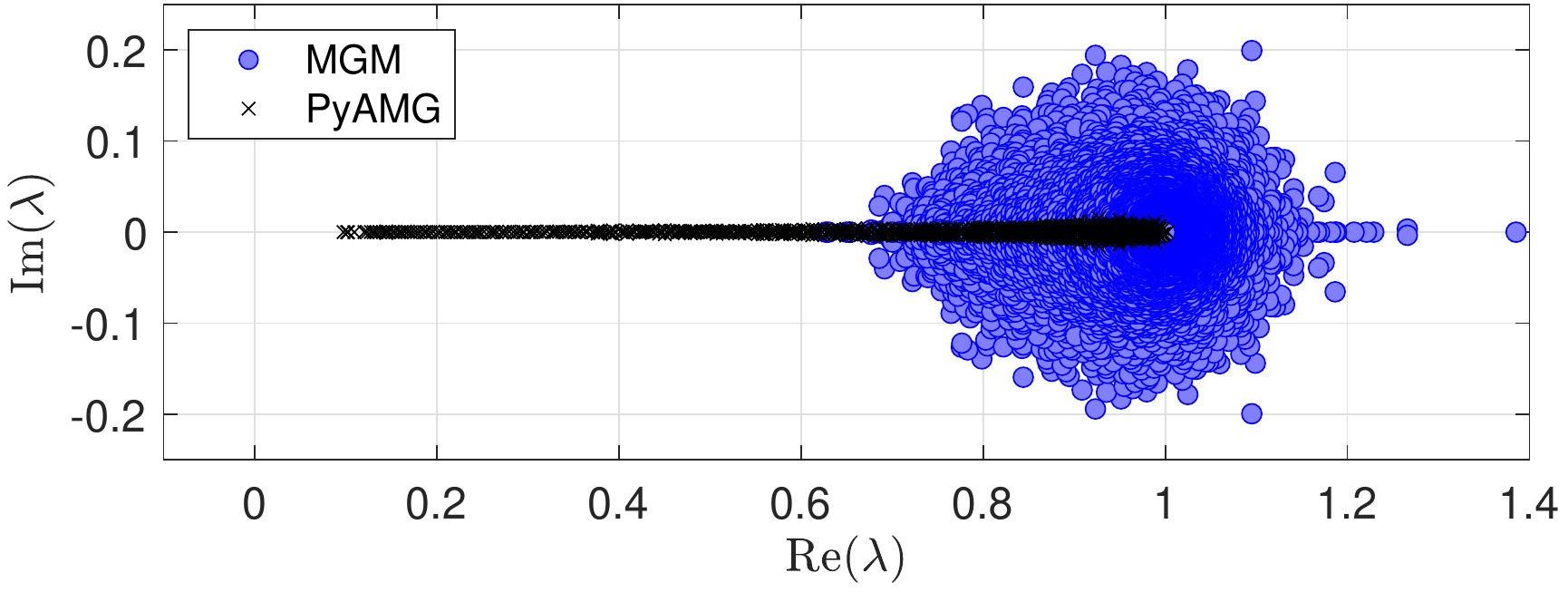}
\end{tabular}
\caption{\label{fig:eigsph} The spectra of the preconditioned matrix for shifted Poisson problem discretized with RBF-FD.}
\end{figure}
The previous two sections showed the preconditioned MGM methods outperforming the PyAMG methods.  To better understand these results, we investigate the spectrum (eigenvalues) of the preconditioned matrices from both methods.  Letting $M_h$ denote the matrix representation for applying one V-cycle of either MGM or PyAMG, we can write the (right) preconditioned system as $L_h M_h z^h = f^h$, where $z^h = (M_h)^{-1} u^h$.  The convergence behavior of Krylov methods can be understood by analyzing the spectrum $L_h M_h$.  As discussed in, for example~\cite{oosterlee1998evaluation}, the more clustered this spectrum is to one, the faster the Krylov methods will converge.  In Figure \ref{fig:eigsph} we display the complete spectrum of the preconditioned matrix $L_h M_h$ of both MGM and PyAMG for the RBF-FD discretizations on the sphere and cyclide.  Due to the cost of this eigenvalue computation, we were only able to compute the results for with $N_h=10242$ and $8192$, respectively.  We see from the figure that spectra for MGM are more clustered around one than PyAMG for both surfaces and increasing $\ell$, which explains the better iteration counts in Table \ref{tbl:iters_vs_N_rbffd}.  We omit the results for GFD, but note that the spectra were similar to RBF-FD, but were even more clustered near one.

%The major drawback of iterative methods is slow convergence for more complicated systems  arising in applications such as fluid dynamics, chemical reactions and transport, and electromagnetics \cite{Saad}. The use of preconditioning for iterative methods can add more efficiency and robustness to the method by stabilizing and accelerating the convergence. The quality of the preconditioner can be determined from the spectrum of preconditioning operator.
%In this section, we present the preconditioner spectra for solving a shifted Poisson problem on the sphere and the cyclide.
%\begin{figure}
%		\centering\includegraphics[width=0.45\textwidth,scale=0.35]{figures/spectra_N10242_lvl3_geosphere_lT0_methodRBF.png}
%		\centering\includegraphics[width=0.45\textwidth,scale=0.35]{figures/spectra_N8192_lvl3_geocyclide_poisson_lT0_methodRBF.png}
%		\caption{\label{fig:eigsph} The eigenvalue spectra of the preconditioning matrix for shifted Poisson problem on the sphere $N_h = 10242$ (left) and Dupin's cyclide $N_h =11884$ (right).}
%\end{figure}

%%%%%%%%%%%%%%%%%%%%%%%%%%%%%%%%%%%%%%%%%%%%%%%%%%%%%%%%%%%%%%%%%%%%%%55
\subsection{Iteration vs.\ accuracy\label{sec:accuracy}}
\begin{figure}[htb]\centering
\begin{tabular}{ccc}
\rotatebox{90}{\hspace{0.1\textwidth}\small $\ell=3$} & 
\includegraphics[width=0.4\textwidth]{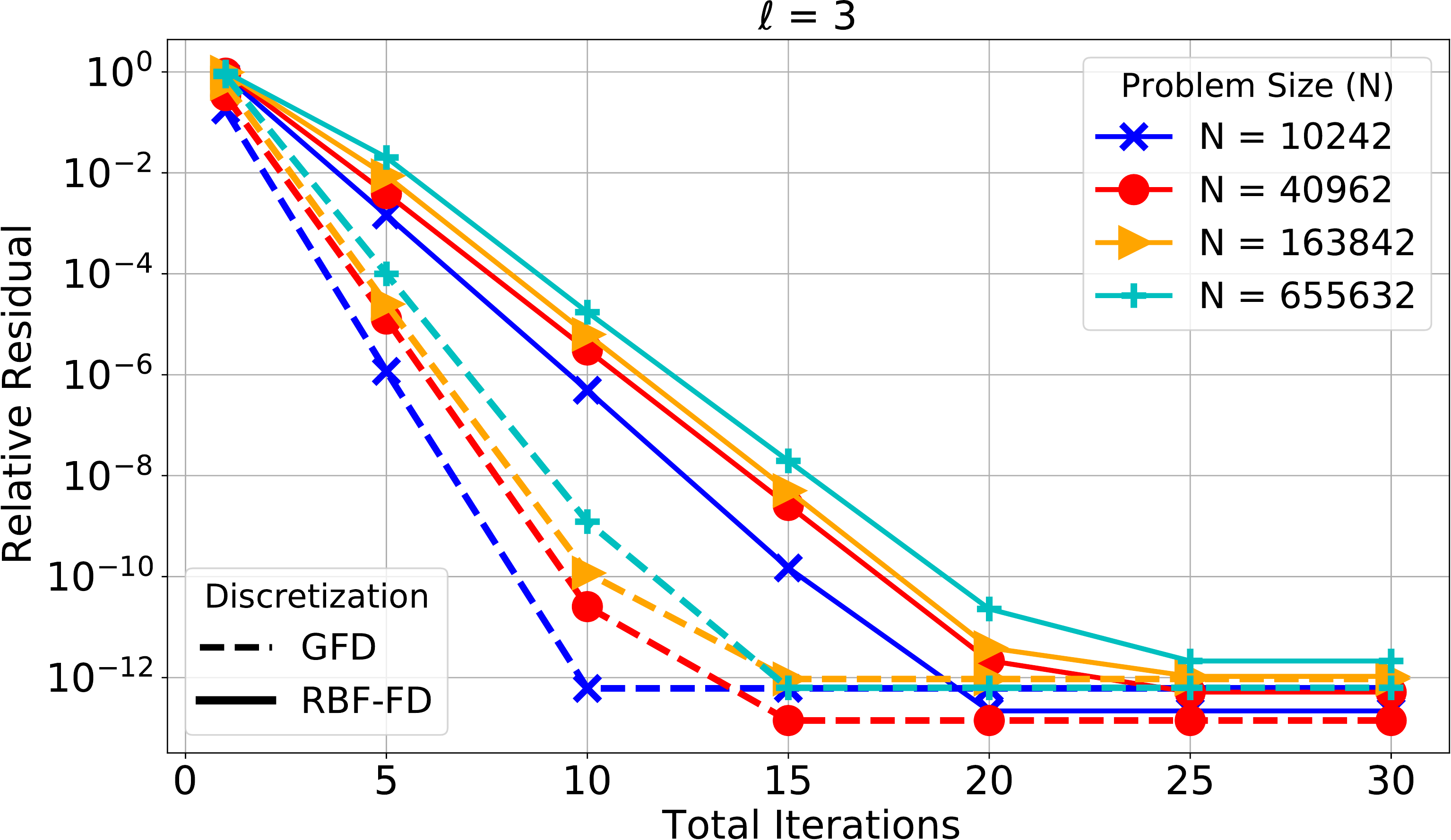} & 
\includegraphics[width=0.4\textwidth]{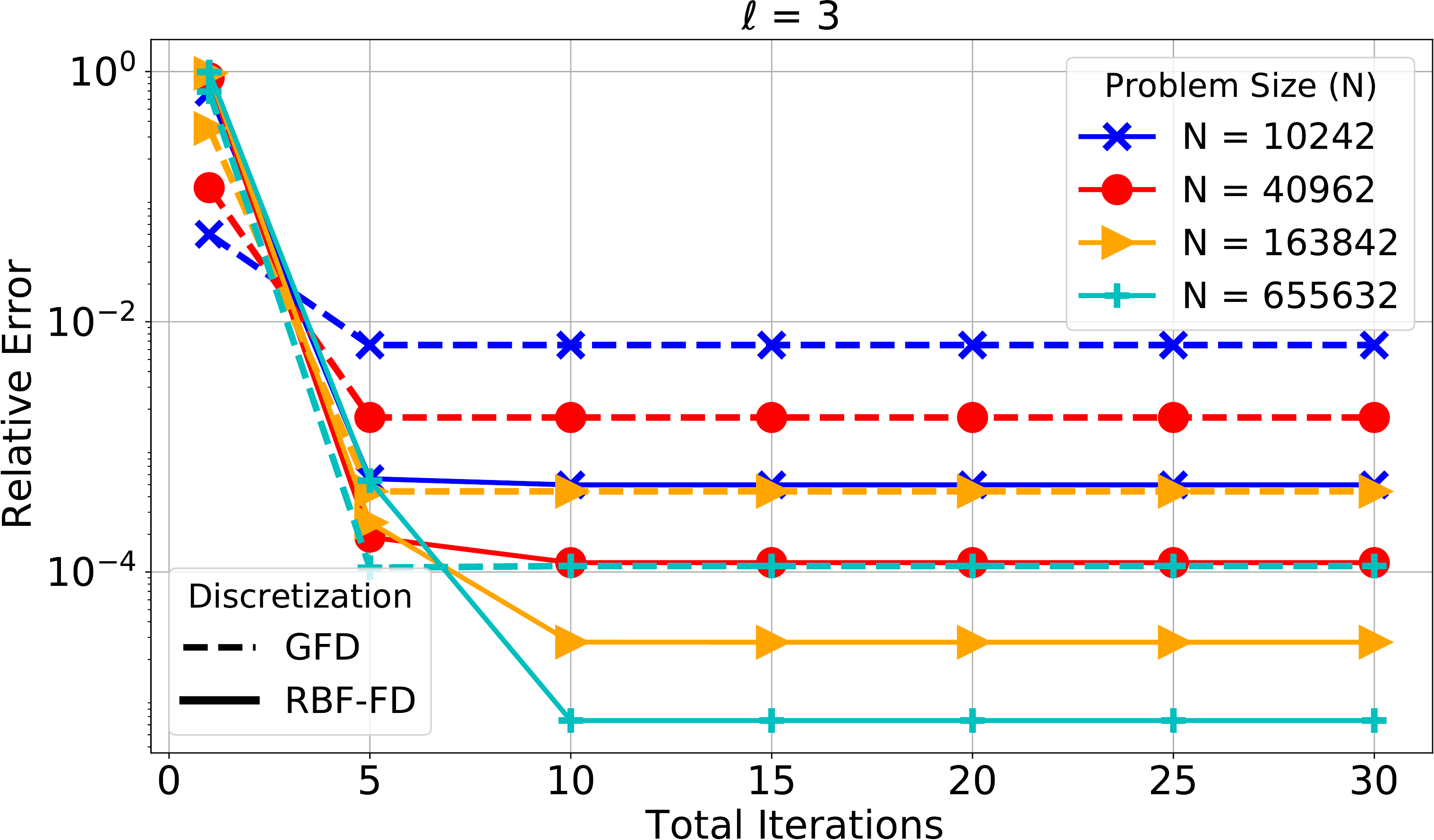} \\ 
\rotatebox{90}{\hspace{0.1\textwidth}\small $\ell=5$} & 
\includegraphics[width=0.4\textwidth]{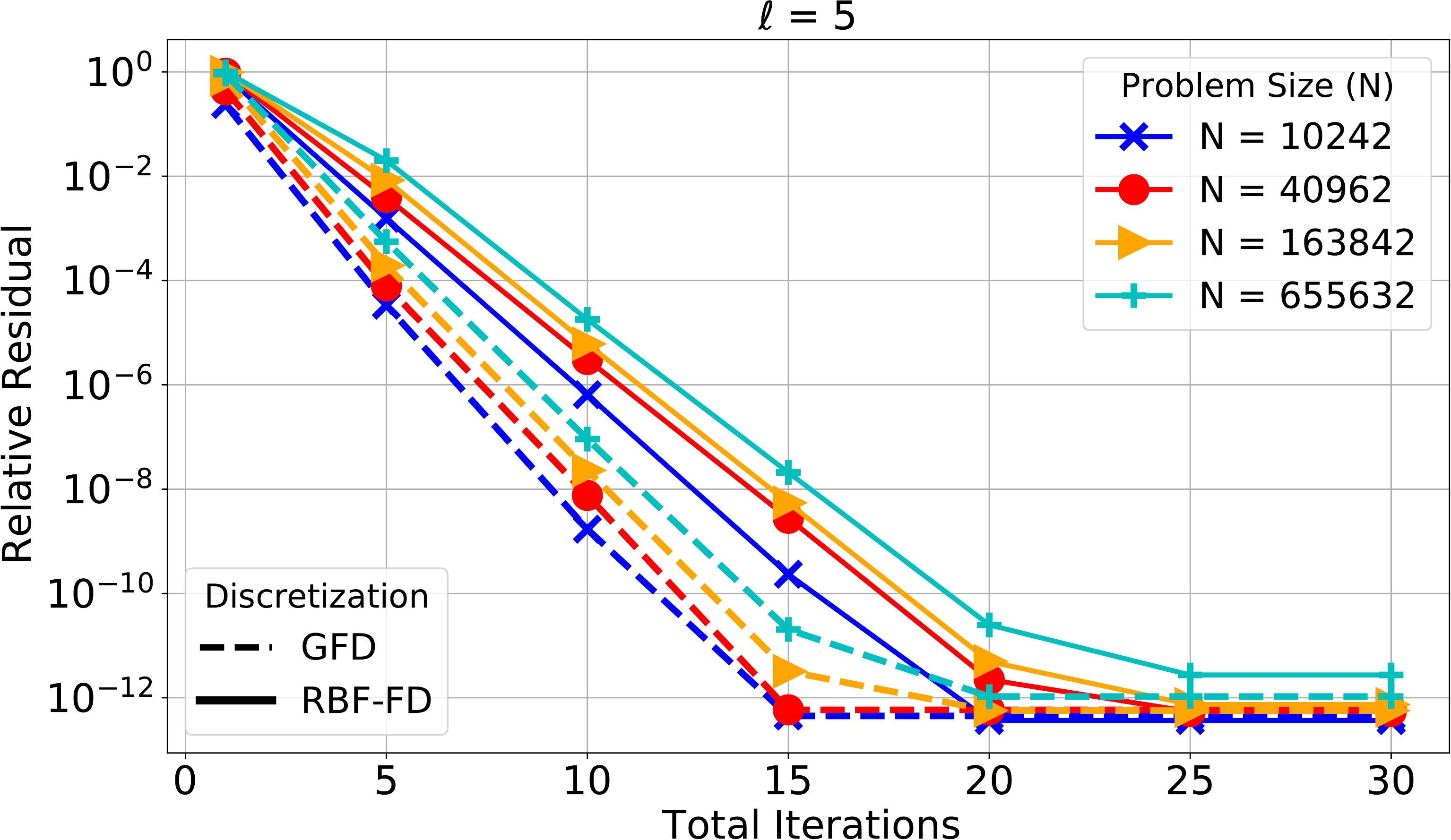} & 
\includegraphics[width=0.4\textwidth]{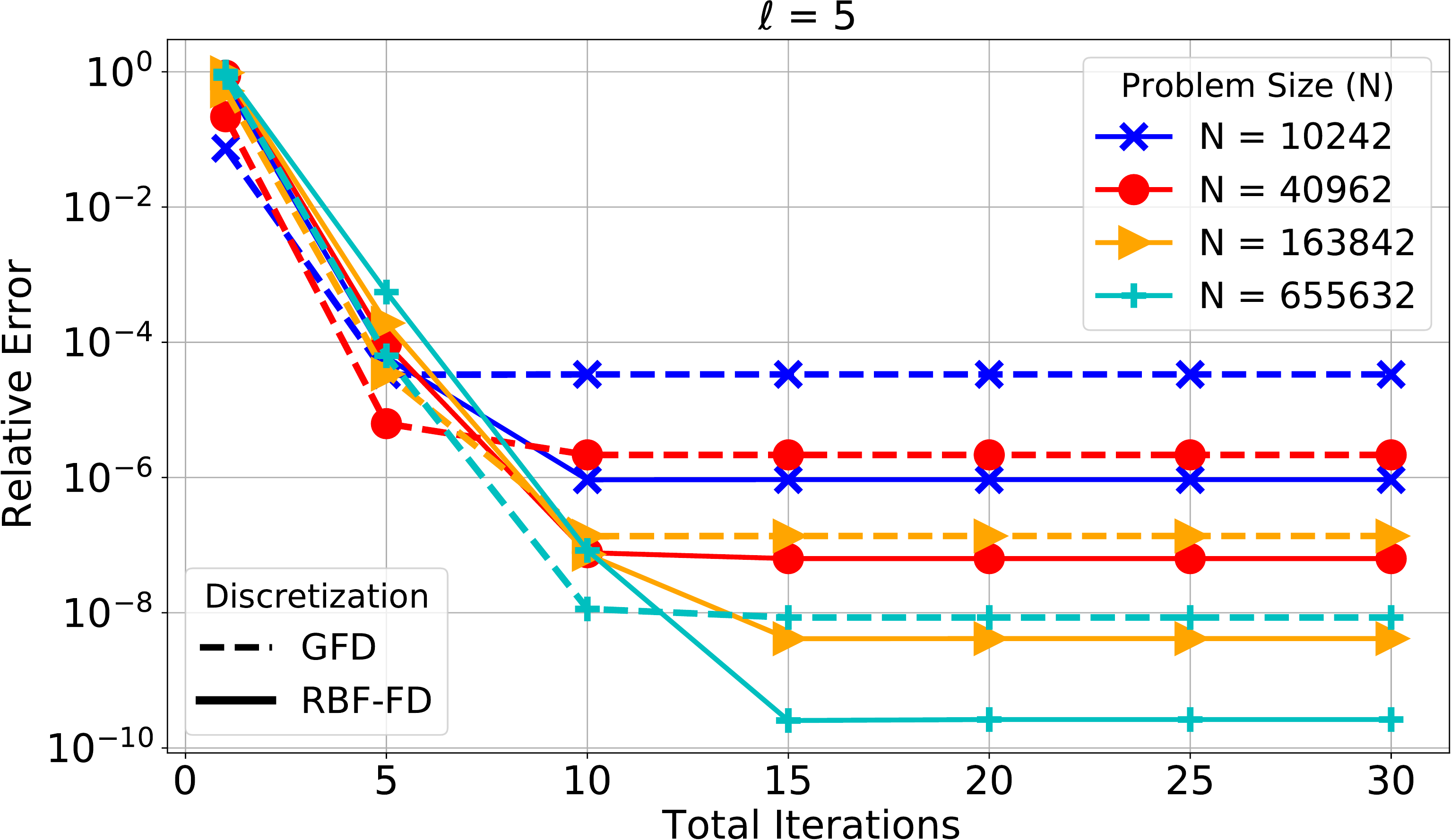} \\ 
\rotatebox{90}{\hspace{0.1\textwidth}\small $\ell=7$} & 
\includegraphics[width=0.4\textwidth]{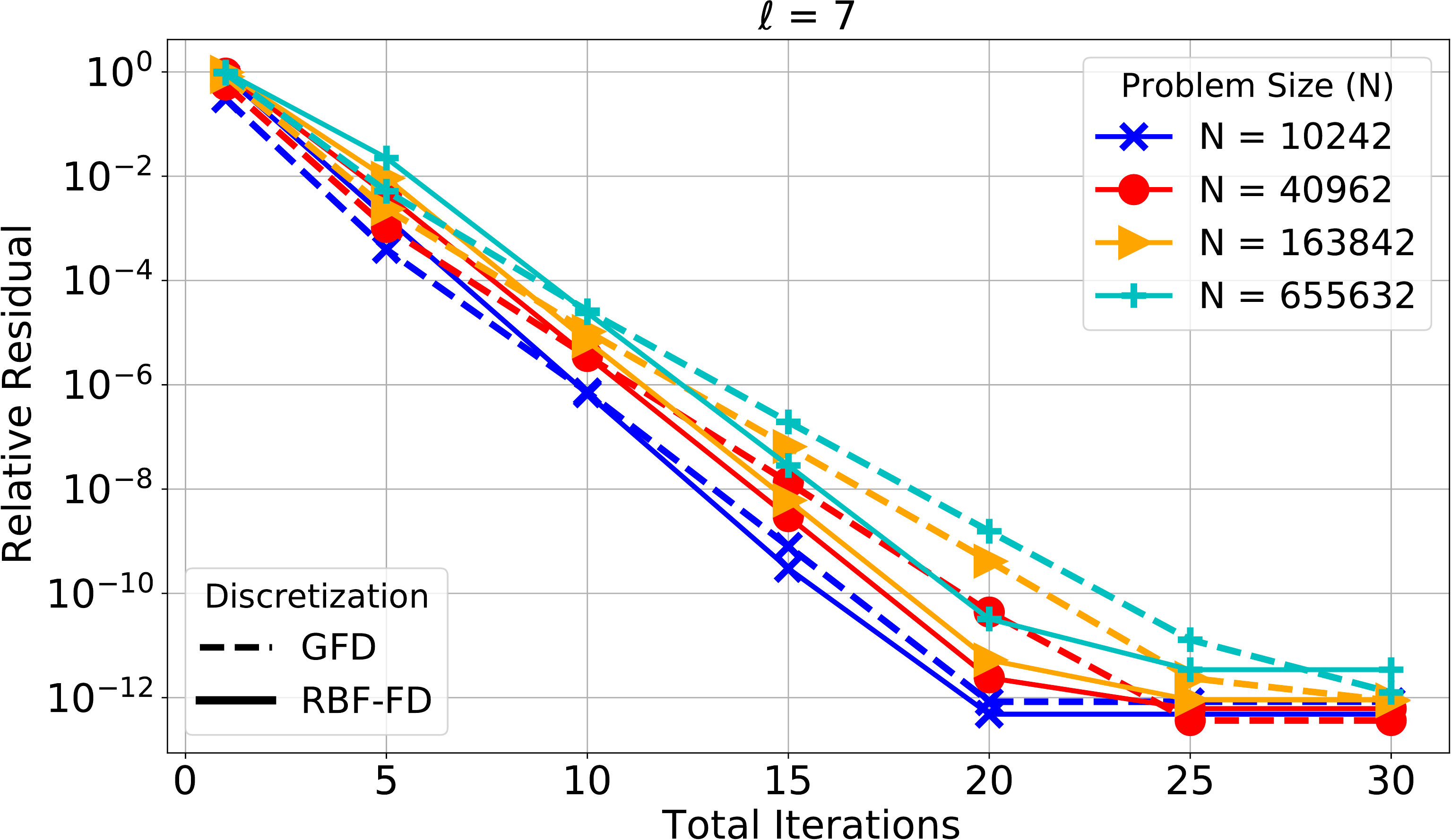} & 
\includegraphics[width=0.4\textwidth]{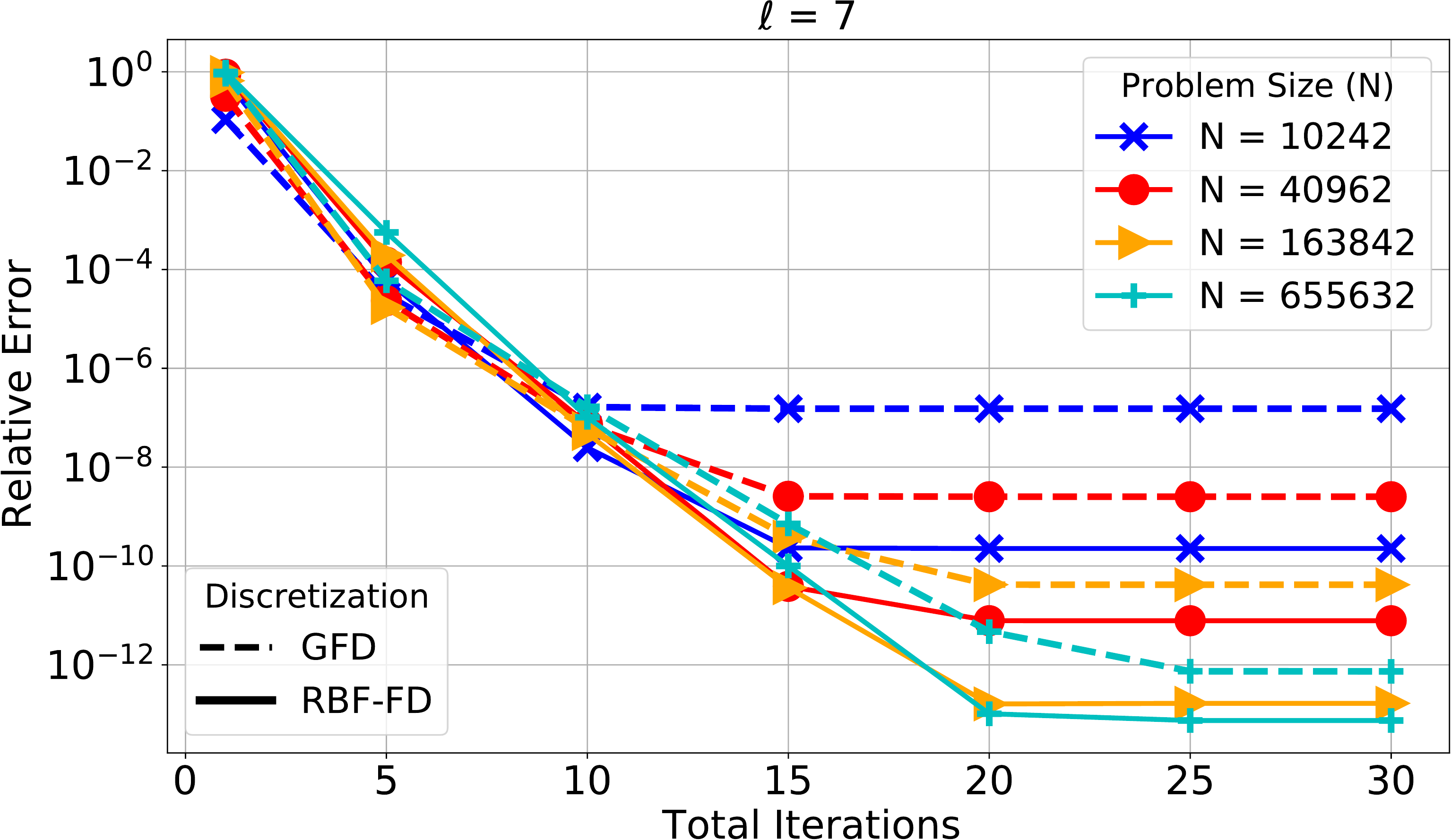} 
\end{tabular}
\caption{Relative residuals (left) and relative 2-norm errors (right) for solving a Poisson problem on the sphere with MGM GMRES.  Solid lines correspond to RBF-FD discretizations, while dashed lines correspond to GFD. \label{fig:accuracy}}
\end{figure}
%The previous results demonstrated that preconditioned MGM is more effective than preconditioned PyAMG for the surfaces and discretizations under consideration.  
In the final set of tests we focus on solving the (discretized) Poisson problem with MGM GMRES and examine how the accuracy of the RBF-FD and GFD discretizations depend on the iteration count for increasing $N_h$ and $\ell$.  We restrict our attention to the sphere, for which it is easy to construct test problems with exact solutions based on spherical harmonics.  For the test problem in the experiments, we use the $Y_5^{4}$ spherical harmonic, which can be written in Cartesian coordinates as $Y_5^4(x,y,z) = z(x^4-6x^2y^2+y^4)$. We fix the number of iterations of MGM GMRES for solving the discretized systems to $1, 5, 10, 15, 20, 25, 30$, and compute both the relative residual and relative errors (in the 2-norm) in the approximate solutions.  Figure \ref{fig:accuracy} displays the results from these experiments.  We see from the figure that in almost all cases the minimum error for either the RBF-FD and GFD is reached before the minimal residual is reached. Additionally, the results indicate that while the residuals for GFD converge faster than RBF-FD, the errors for a given $N_h$ and $\ell$ are smaller for RBF-FD.  So the cost per error for both methods is much more comparable than the previous experiments indicated and favor RBF-FD.

\section{Applications\label{sec:applications}}
\begin{figure}[htb]
    \centering
    \begin{tabular}{ccc}
    \includegraphics[width=0.18\textwidth]{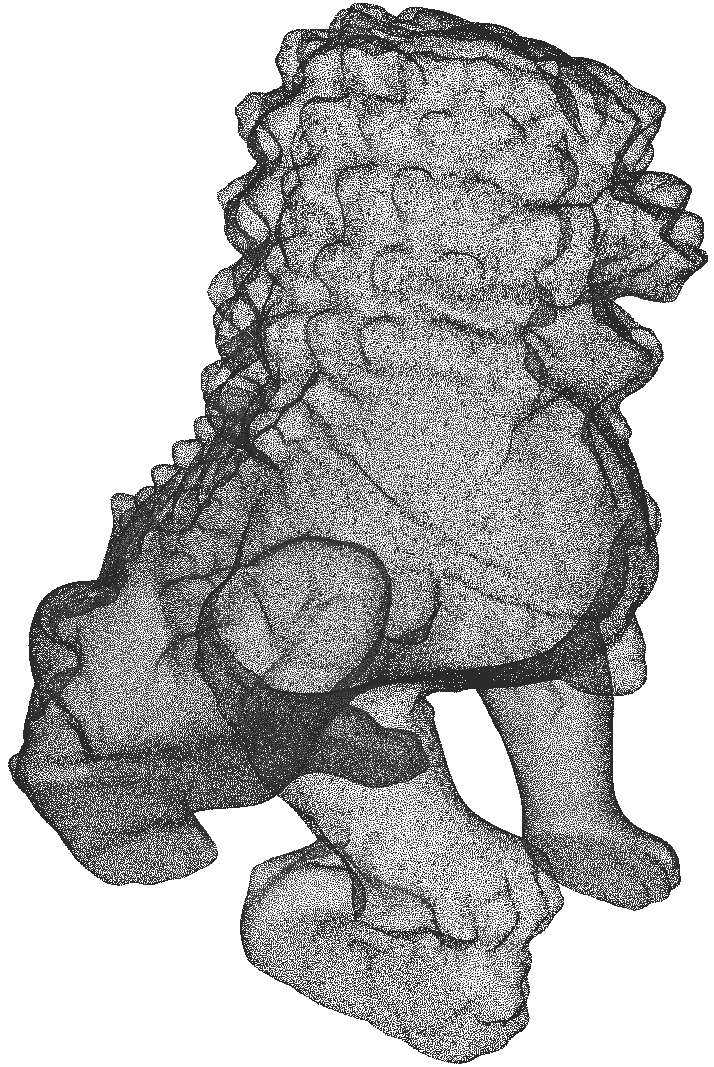}  & 
    \includegraphics[width=0.23\textwidth]{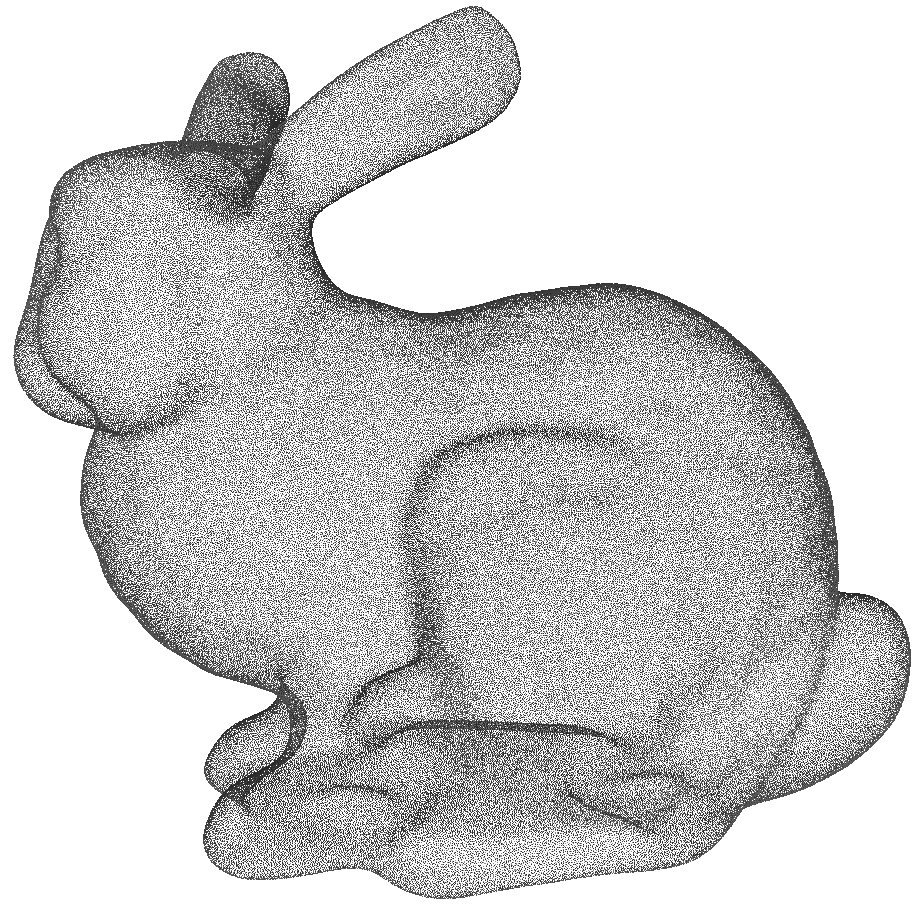} &
    \includegraphics[width=0.28\textwidth]{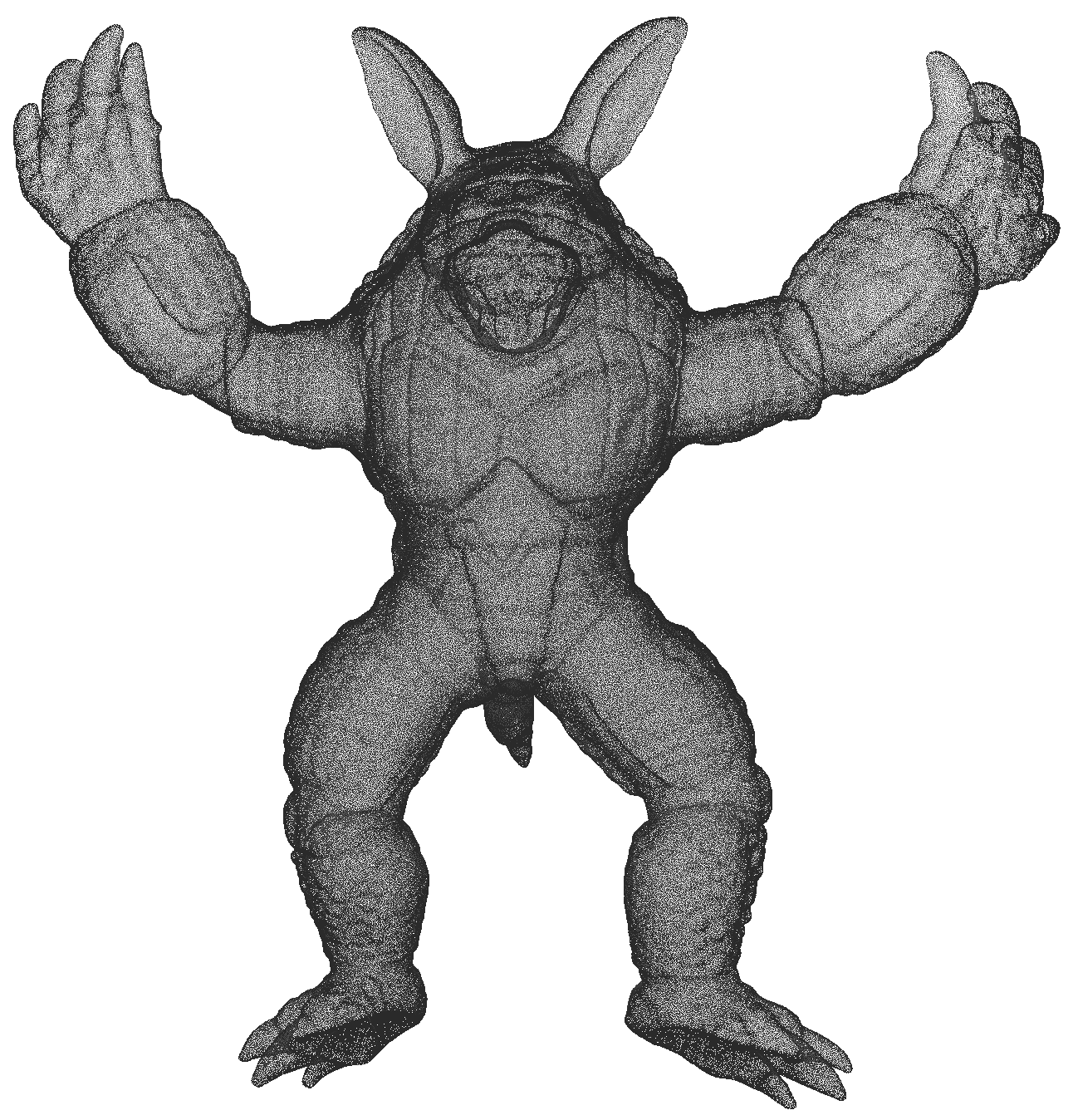}
    \end{tabular}
  \caption{Point clouds for the surfaces considered in the applications: Chinese Guardian Lion ($N_h = 436605$), Stanford Bunny ($N_h = 291804$), and Armadillo ($N_h = 872773$).\label{fig:point_clouds}}
\end{figure}

%%%%%%%%%%%%%%%%%%%%%%%%%%%%%%%%%%%%%%%%%%%%%%%%%%%%%%%%%%%%%%%%%%%%%%
In this section we demonstrate the performance of MGM on three different applications involving complicated surfaces represented by relatively large point clouds; see Figure \ref{fig:point_clouds}.  All these applications involve solving discrete (shifted) surface Poisson problems, for which we use the RBF-FD method to approximate the LBO and MGM GMRES to solve the resulting linear systems.

\subsection{Surface harmonics}
We first consider approximating the first several eigenvalues and eigenfunctions of the LBO on the Chinese Guardian Lion model.  The eigenfunctions of the LBO or  the ``surface harmonics'' have been used in various applications in data analysis.  For example, Reuter et.\ al.~\cite{reuter2009discrete} used the low frequency surface harmonics for shape segmentation and registration.  

\begin{figure}[htb]
    \centering
    \begin{tabular}{ccccc}
    \includegraphics[width=0.17\textwidth]{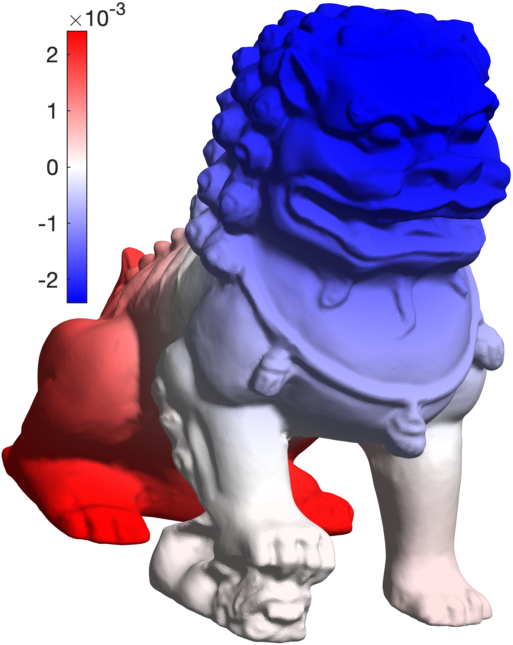}  & 
    \includegraphics[width=0.17\textwidth]{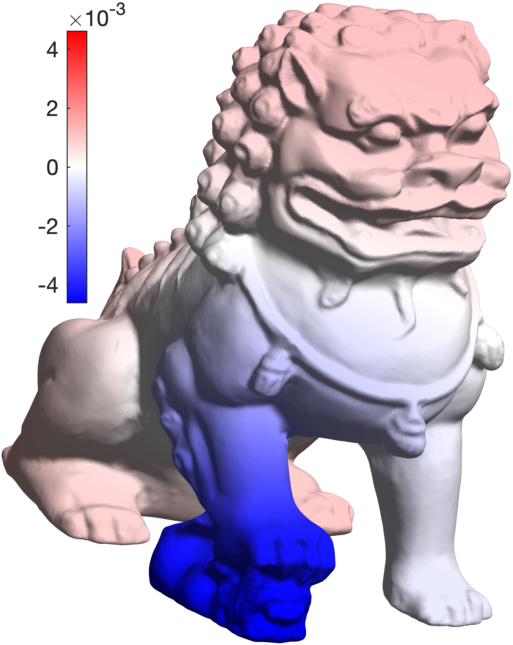}  & 
    \includegraphics[width=0.17\textwidth]{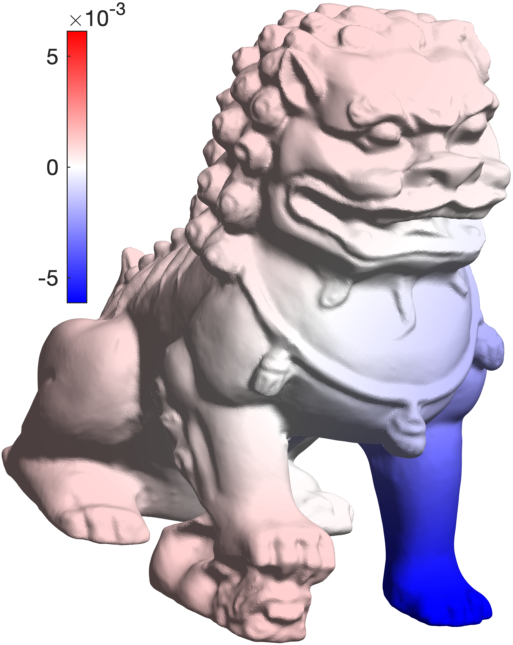}  & 
    \includegraphics[width=0.17\textwidth]{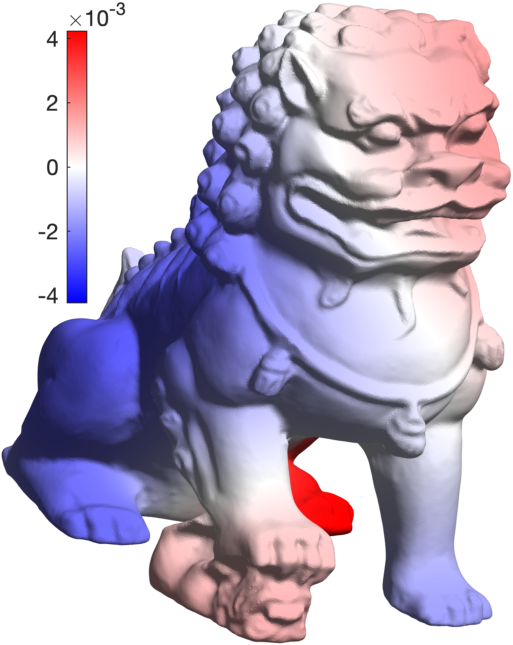}  & 
    \includegraphics[width=0.17\textwidth]{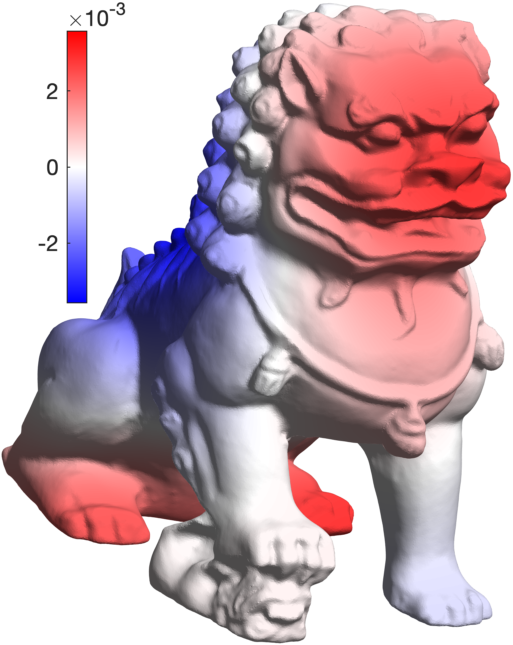}  \\
    \includegraphics[width=0.17\textwidth]{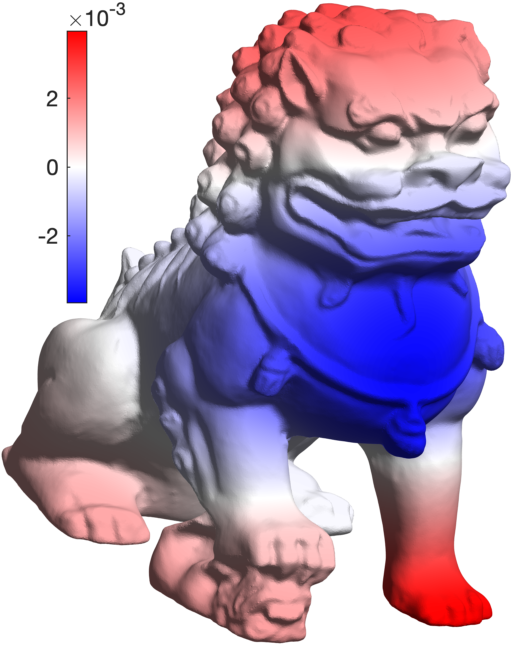}  & 
    \includegraphics[width=0.17\textwidth]{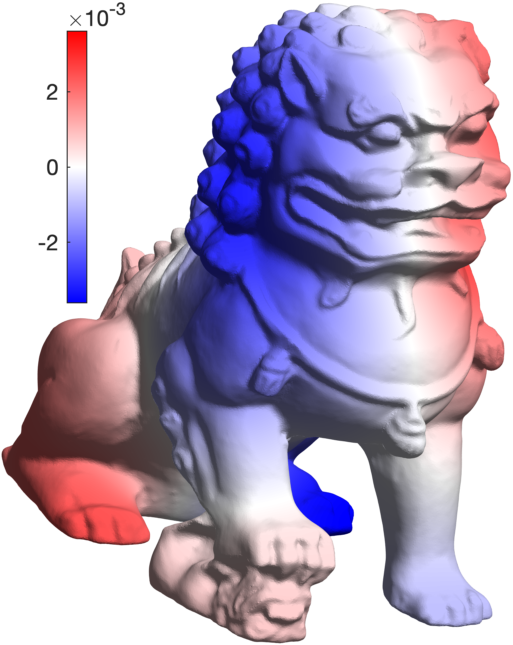}  & 
    \includegraphics[width=0.17\textwidth]{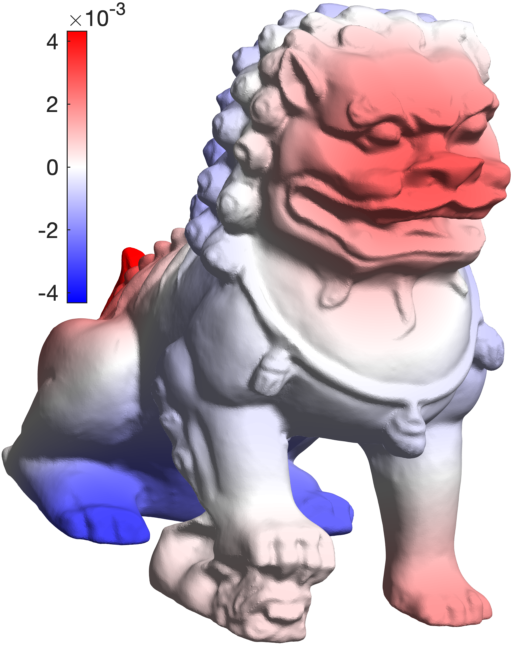}  & 
    \includegraphics[width=0.17\textwidth]{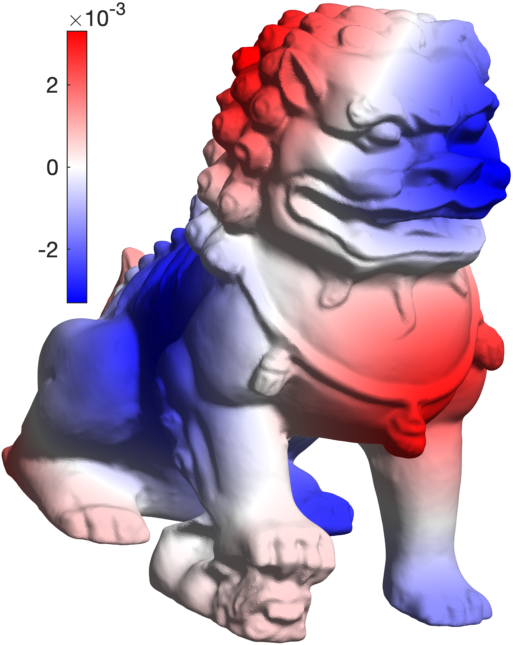}  & 
    \includegraphics[width=0.17\textwidth]{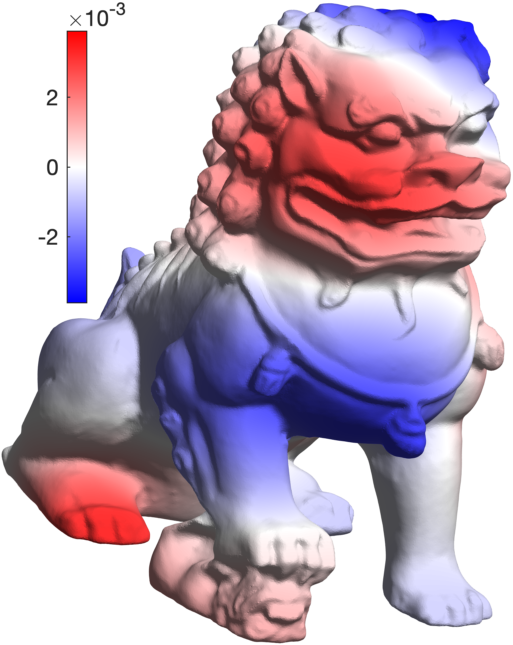}  \\    
    \end{tabular}
  \caption{Left: pseudocolor map of the first 10 non-zero surface harmonics of the Chinese Guardian Lion model.}\label{fig:dragon_harmonics}
\end{figure}

The LBO eigenvalue problem is given as $\laps u = \lambda u$.
%\begin{align*}
%\laps u = \lambda u.
%\end{align*}
To approximate the solutions of this problem we use the RBF-FD method with $\ell=5$ to approximate the LBO and ARPACK~\cite{lehoucq1998arpack} (accessed through the \texttt{eigs} function in MATLAB) to solve the discrete system for the first several eigenpairs that are smallest in magnitude.  ARPACK uses the Arnoldi method on the shifted inverse of a matrix to find the eigenpairs closest to the shift $\sigma$, which, for the surface problem, requires a routine for repeatedly solving systems of the form $(L_h - \sigma I_h)v^h = f^h$, for different $f^h$.  We use MGM GMRES to solve these linear systems with $\sigma=-1$ and set the tolerance to $10^{-10}$.  Figure \ref{fig:dragon_harmonics} displays the first 10 non-zero harmonics computed with this technique.
The ARPACK routine used 49 linear system solves to determine the eigenpairs; the median number of MGM GMRES iterations required to solve these systems was only 15 and the max was 16.

%%%%%%%%%%%%%%%%%%%%%%%%%%%%%%%%%%%%%%%%%%%%%%%%%%%%%%%%%%%%%%%%%%%%%%
\subsection{Pattern formation}
We next consider solving two coupled reaction-diffusion (RD) equations on the Stanford Bunny model.  These types of equations arise, for example, in phenomenological models of color patterns in animal coats~\cite{miyazawa2010blending}.  We consider the Gierer-Meinhardt two-species RD system~\cite{gierer1972theory} given as follows:
\begin{subequations}
\begin{align}
\frac{\partial u}{\partial t} =& D_u\laps\, u + A - Bu + \frac{u^2}{v(1+Cu^2)}, \label{eq:turingu} \\
\frac{\partial v}{\partial t} =& D_v\laps\, v + u^2 - v. \label{eq:turingv}
\end{align}
\label{eq:turing}
\end{subequations}
By altering the parameters $A$, $B$, $C$, $D_u$, and $D_v$ appropriately, this system can produce solutions that converge to spot or labyrinth patterns at ``steady-state''~\cite{miyazawa2010blending}.  For the bunny model, we set $A=0.08$, $B=1.5$, $C=0.45$, $D_u = 5\times 10^{-5}$, and $D_v = 10^{-3}$ to produce the labyrinth pattern.  We use a random initial condition, where at each point in $X_h$ the values of $u$ and $v$ are selected from a uniformly random distribution in the interval $[0,1]$.
%\begin{figure}[htb]
%    \centering
%    \begin{tabular}{ccc}
%    \includegraphics[width=0.23\textwidth]{figures/bunny_spots_N291804.png}  & \includegraphics[width=0.23\textwidth]{figures/bunny_stripes_N291804.png} & \includegraphics[width=0.45\textwidth]{figures/stripe_spot_iterations.pdf}
%    \end{tabular}
%  \caption{The illustration shows the formation of Turing spots on the Standford Bunny. }\label{fig:bunny_spots}
%\end{figure}
\begin{figure}[htb]
    \centering
    \begin{tabular}{cc}
    \includegraphics[width=0.27\textwidth]{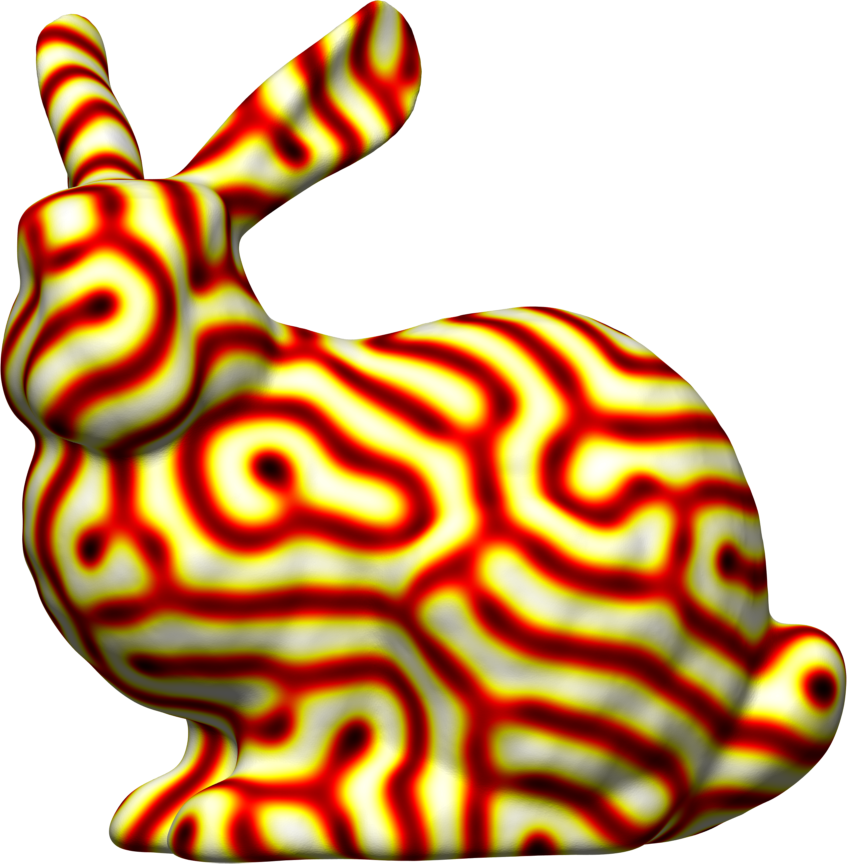}  & 
    \includegraphics[width=0.55\textwidth]{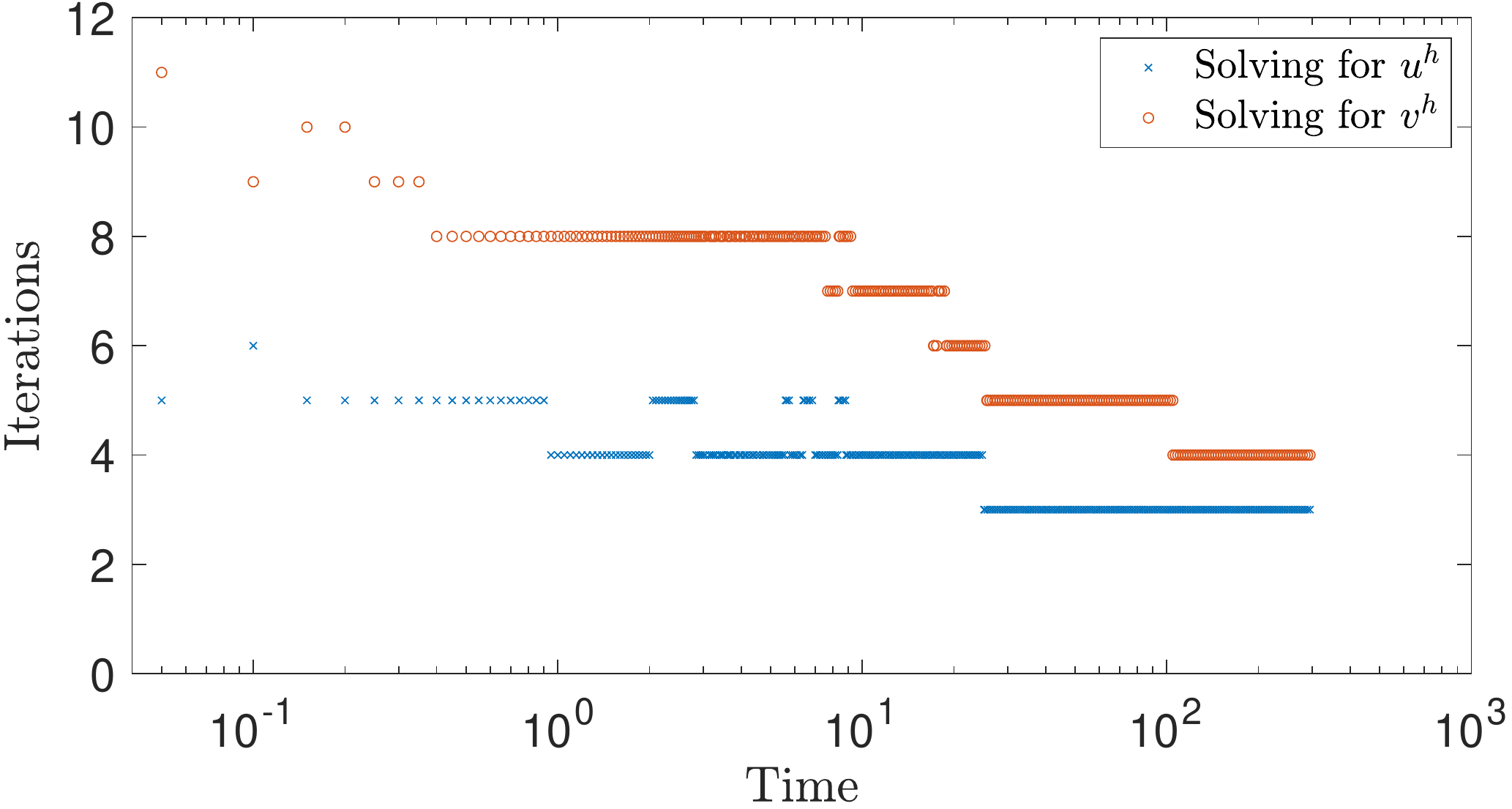}
    \end{tabular}
  \caption{Left: pseudocolor map of the $u$ variable in the numerical solution of \eqref{eq:turing} on the Stanford Bunny model; the colors transition from white to yellow to red to black, with white corresponding to $u=0$ and black to $u=1$. Right: iteration count of MGM GMRES for solving the linear systems associated with $u$ and $v$ variables at each time-step in the semi-implicit scheme for \eqref{eq:turing}.}\label{fig:bunny_labyrinth}
\end{figure}

To approximate the solution of \eqref{eq:turing} we use the RBF-FD method with $\ell=3$ to approximate the LBO and apply the second-order accurate semi-implicit backward difference scheme (SBDF2)~\cite{Ascher97} as the time-stepping method that treats the diffusion implicitly and reactions explicitly.  We set the time-step to $\Delta t = 0.05$.  The temporal discretization results in two decoupled (discrete) screened Poisson problems that need to be solved at each time-step for which we use GMRES preconditioned with MGM.  For the GMRES method we set the tolerance on the relative residual to $10^{-8}$ and use the previous time-step as the initial guess.  We set the final integration to $300$ time units, which resulted in a near steady-state pattern.  Figure \ref{fig:bunny_labyrinth} displays the results of the simulations.  Included in the figure are the iterations required by the preconditioned GMRES method as a function of time.  We see from the figure that the maximum iteration count is 6 for the $u$ variable and 11 for the $v$ variable, and decreases to 2 and 3, respectively as the solutions approach steady-state.  The larger iteration count for the $v$ variable is expected since the diffusion coefficient is larger in  \eqref{eq:turingv}.

\subsection{Geodesic distance}
Lastly, we consider the classic problem of approximating the geodesic distance from a given point on a surface to all other points.  We use the \textit{heat method} introduced by Crane et.\ al.~\cite{crane2013geodesics} to solve this problem.  This method transforms the non-linear geodesic distance problem, typically formulated in terms of the eikonal equation, into solving a pair of linear parabolic and elliptic problems.  The heat method is comprised of the three steps:
\begin{enumerate}
\item Solve $u_t = \laps u$, with $u_0 = \delta(\vx^*)$, to some time $t_{\rm final} > 0$
\item Compute the vector field $\boldsymbol{\eta} = -\nabla_{\M}u/|\nabla_{\M} u|$
\item Solve the Poisson problem $\laps \varphi = \nabla_{\M}\cdot \boldsymbol{\eta}$
\end{enumerate}
Here $\vx^*\in$ denotes the target point on the surface $\M$ to compute the distance from, $\nabla_{\M}$ denotes the surface gradient, $\nabla_{\M} \cdot$ is the surface divergence, and $\delta$ denotes the Dirac delta function.  As discussed in~\cite{crane2013geodesics}, the function $\varphi$ approximates the geodesic distance and converges to the exact distance as $t_{\rm final}\rightarrow 0$.

\begin{figure}[htb]
    \centering
    \begin{tabular}{ccc}
    \includegraphics[width=0.28\textwidth]{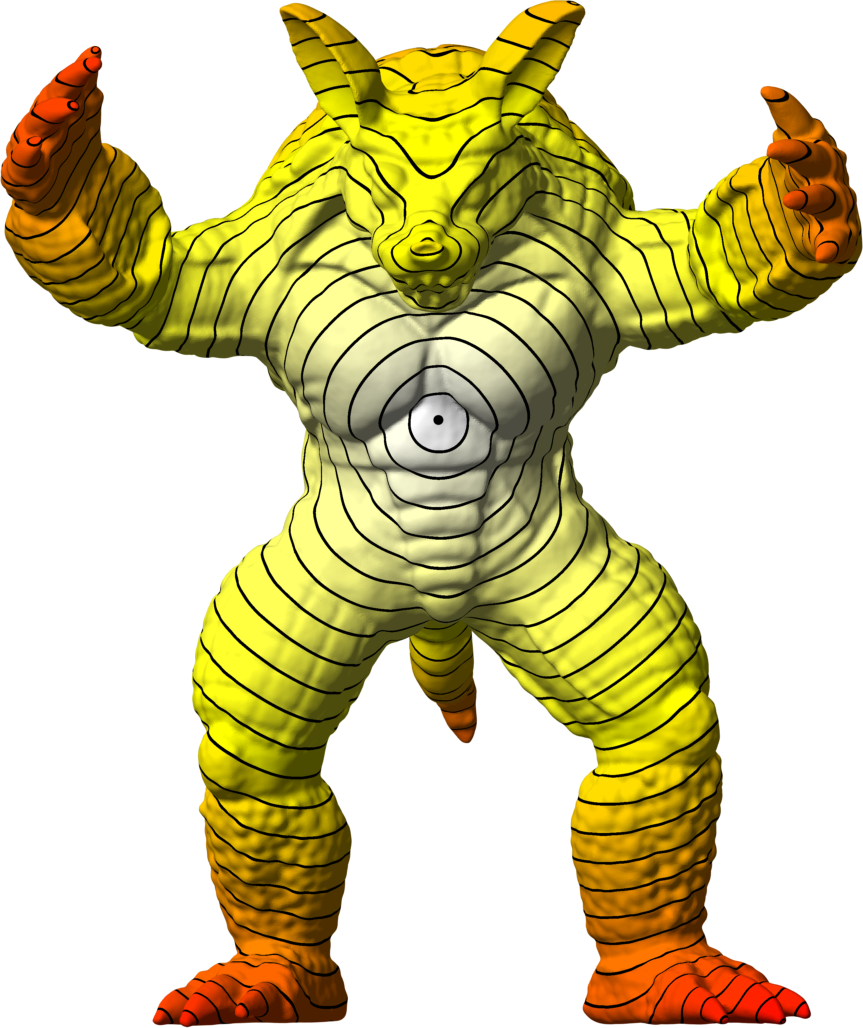}  & 
    \includegraphics[width=0.28\textwidth]{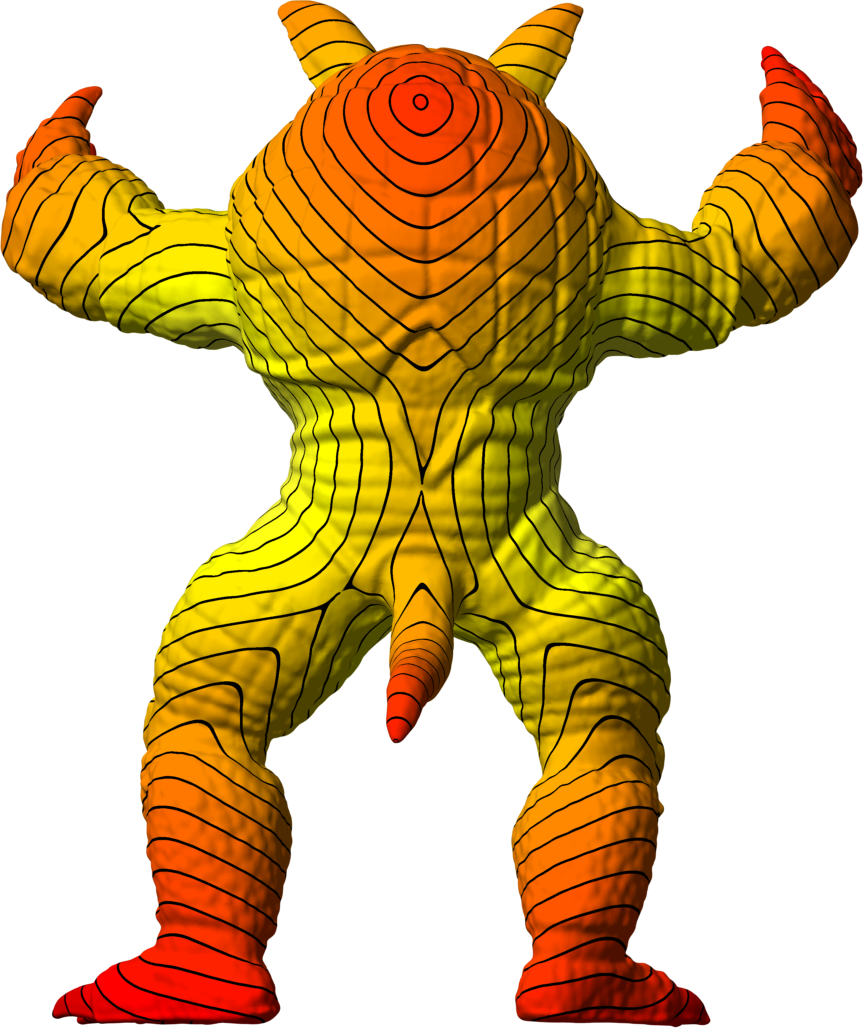} &
    \includegraphics[width=0.28\textwidth]{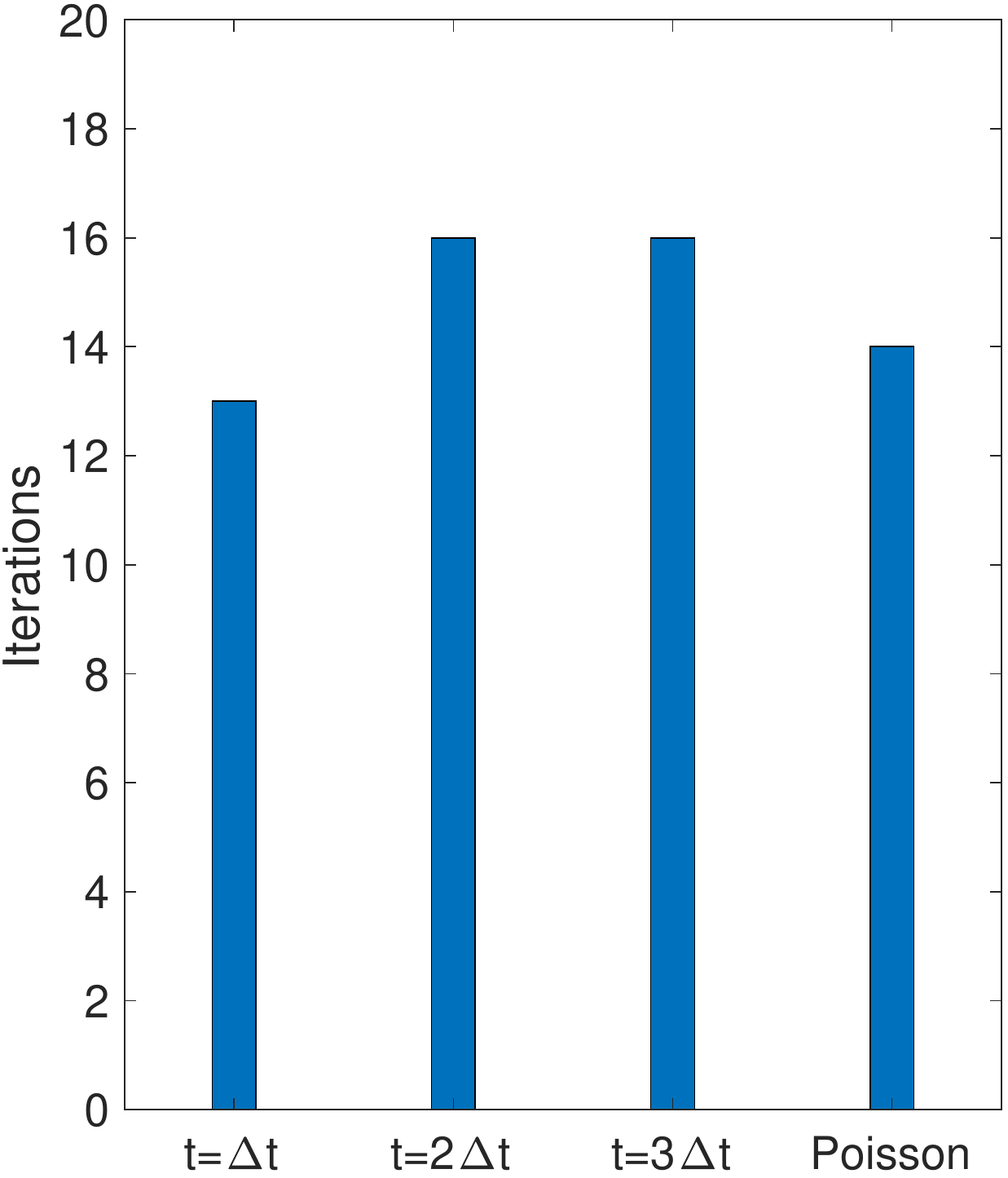}
    \end{tabular}
  \caption{Left: pseudocolor map of the approximate geodesic distance from the solid circle on the chest of the Armadillo model (viewed from the front and backside) computed with the heat method.  Solid black lines mark the contours of the distance field and the colors transition from white to yellow to red with increasing distance from the solid circle.  Right: iteration count of GMRES preconditioned with MGM for solving the linear systems associated with the heat method.\label{fig:geodesic}}
\end{figure}

We apply the heat method on the Armadillo model.  We again use the RBF-FD method with $\ell=3$ to approximate the LBO in steps 1 and 3 above.  To approximate the surface gradient and divergence, we also use the RBF-FD method formulated in the tangent plane similar to the method described in~\cite{SUCHDE20192789} for GFD.  For these approximations, we use $\ell=2$, which result in a second-order approximation.  We discretize the heat equation in the first step with backward Euler in time with a time-step of $\Delta t=10^{-3}$ and set $t_{\rm final} = 3\Delta t$.  To solve the linear systems associated with this implicit discretization and the system from the discretized Poisson equation in step 3, we use GMRES preconditioned with MGM, setting the tolerance to $10^{-8}$.  The results for a point $\vx^{*}$ on the chest of the Armadillo are displayed in the first to images of Figure \ref{fig:geodesic}.  The last image in this figure displays the iterations of the preconditioned GMRES method for solving the systems from the heat equation discretization for three time-steps and the Poisson system to determine $\varphi$. We see that the iteration count remains low for all these systems.

\section{Concluding remarks}\label{sec:discussion}
We have presented a new geometric multilevel method, MGM, for solving linear systems associated with discretizations of elliptic PDEs on point clouds.  The method is entirely meshfree and uses the WSE algorithm for coarsening the point clouds, interpolation/restrictions operators based on polyharmonic spline RBFs, Galerkin coarsening of the operator, and standard smoothers.  All of these choices make MGM particularly straightforward to implement.  We numerically analyzed the method as a standalone solver and preconditioner on test problems for the sphere and cyclide discretized using RBF-FD and GFD methods, and found that it compares favorably to AMG methods in terms of convergence rates and wall-clock time.  When using MGM as a preconditioner, we also found that it scaled well as both the problem size and accuracy of the discretizations increased.  Finally, we demonstrated that the method can be used in three challenging applications involving large systems of equations.

There are several extensions of MGM that we plan to pursue in the future.  One is to test the method on other discretizations.  MGM is agnostic to the underlying discretization and could be used even for (nodal) mesh-based discretizations.  Here the nodal points of the mesh could be treated as a point cloud and WSE could be applied, or if there is a natural way to coarsen the mesh, then this could be used instead.  A second idea we plan to pursue is extending MGM to domains with boundaries, which in principle should be straightforward.  Finally, we plan to look into parallel implementations of the method to further improve the performance.

\section*{Acknowledgements} 
AMJ and GBW's work was partially supported by US NSF grants CCF-1717556 and DMS-1952674.  VS's work was partially supported by US NSF grants CCF-1714844. We benefitted from discussions with Drs. Cem Yuksel (on the WSE algorithm) and Hari Sundar (on multigrid methods), both from the University of Utah. The Stanford Bunny and Armadillo data were obtained from the Stanford University 3D Scanning Repository.  The Chinese Guardian Lion data was obtained from the AIM@SHAPE-VISIONAIR Shape Repository.

\bibliographystyle{siamplain}
\bibliography{biber}

\end{document}